\documentclass{article}

\usepackage{arxiv}
\usepackage{tikz}
\usepackage[utf8]{inputenc} 
\usepackage[T1]{fontenc}    
\usepackage{hyperref}       
\usepackage{url}            
\usepackage{booktabs}       
\usepackage{amsfonts}       
\usepackage{nicefrac}       
\usepackage{microtype}      
\usepackage{lipsum}

\usepackage{amssymb}
\usepackage{lineno,hyperref}
\usepackage{amsmath}
\usepackage{bm}
\usepackage{subcaption}
\usepackage{graphicx}
\usepackage{svg}
\usepackage{rotating}
\usepackage{lscape}
\usepackage{array,multirow}
\usepackage[title]{appendix}
\usepackage{float}
\usepackage{stmaryrd}
\usepackage{textcomp}
\usepackage{algorithm}
\usepackage{algpseudocode}
\usepackage{amsthm}
\usepackage{listings}
\usepackage{xcolor}
\usepackage{siunitx}
\usepackage{symbols}
\usepackage{pifont}
\usepackage{lineno}

\DeclareMathOperator{\tr}{tr}

\title{Non-intrusive reduced order modeling of poroelasticity of heterogeneous media based on a discontinuous Galerkin approximation}

\author{
  Teeratorn Kadeethum \\
  Sibley School of Mechanical and Aerospace Engineering\\
  Cornell University, New York, USA\
   Ithaca, New York, USA \\
  \texttt{tk658@cornell.edu} \\
   \And
 Francesco Ballarin\\
  mathLab, Mathematics Area\\
   SISSA, Italy \\
  \texttt{francesco.ballarin@sissa.it} \\
   \And
 Nikolaos Bouklas \\
  Sibley School of Mechanical and Aerospace Engineering\\
  Cornell University, New York, USA\
   Ithaca, New York, USA \\
  \texttt{nb589@cornell.edu} \\
}

\begin{document}
\maketitle

\begin{abstract}
A simulation tool capable of speeding up the calculation for linear poroelasticity problems in heterogeneous porous media is of large practical interest for engineers, in particular, to effectively perform sensitivity analyses, uncertainty quantification, optimization, or control operations on the fluid pressure and bulk deformation fields. Towards this goal, we present here a non-intrusive model reduction framework using proper orthogonal decomposition (POD) and neural networks, based on the usual offline-online paradigm. As the conductivity of porous media can be highly heterogeneous and span several orders of magnitude, we utilize the interior penalty discontinuous Galerkin (DG) method as a full order solver to handle discontinuity and ensure local mass conservation during the offline stage.
We then use POD as a data compression tool and compare the nested POD technique, in which time and uncertain parameter domains are compressed consecutively, to the classical POD method in which all domains are compressed simultaneously. The neural networks are finally trained to map the set of uncertain parameters, which could correspond to material properties, boundary conditions, or geometric characteristics, to the collection of coefficients calculated from an $L^2$ projection over the reduced basis. We then perform a non-intrusive evaluation of the neural networks to obtain coefficients corresponding to new values of the uncertain parameters during the online stage. We show that our framework provides reasonable approximations of the DG solution, but it is significantly faster. Moreover, the reduced order framework can capture sharp discontinuities of both displacement and pressure fields resulting from the heterogeneity in the media conductivity, which is generally challenging for intrusive reduced order methods. The sources of error are presented, showing that the nested POD technique is computationally advantageous and still provides comparable accuracy to the classical POD method. We also explore the effect of different choices of the hyperparameters of the neural network on the framework performance.
\end{abstract}

\keywords{poroelasticity \and reduced order modeling \and neural networks \and discontinuous Galerkin \and heterogeneity \and finite element}

\section{Introduction}

Coupled hydro-mechanical (HM) processes in porous media, in which fluid flow and solid deformation tightly interact, are involved in various problems ranging from groundwater and contaminant hydrology to biomedical engineering \cite{bisdom2016geometrically, juanes2016were, lee2016pressure,nick2013reactive, choo2018cracking, kadeethum2020enriched, yu2020poroelastic, vinje2018fluid, kadeethum2019investigation, bouklas2015nonlinear, salimzadeh2019effect}. Porous media are generally anisotropic and heterogeneous. Additionally, the sharp contrast of phases of the porous microstructure leads to discontinuous material properties that can span several orders of magnitude \cite{matthai2009upscaling,kadeethum2020finite,baker2015practical,jia2017comprehensive,kadeethum2018investigation,bisdom2016geometrically,muljadi2016impact,nicolaides2015impact}. The sharp discontinuities of the strongly heterogeneous microstructure significantly influence the time-dependent solid-fluid interactions. For instance, spatially dependent volumetric deformation of a porous medium caused by fluid pressure changes in pore spaces may impact the hydraulic storability and permeability of porous material, which enhances the complexity of the fluid flow field \cite{chen2007reservoir,Du2007,kadeethum2019investigation,kadeethum2020well}. Therefore, accurate numerical modeling of coupled HM processes requires a computational method that can robustly handle substantial heterogeneity in porous media and ensure local mass conservation.

Solution of coupled HM processes can be handled analytically for simple geometries \cite{terzaghi1951theoretical,wang2017theory}. However, for complex geometries and non-homogeneous boundary conditions, we could use numerical approximations such as finite volume discretization \cite{nordbotten2014cell,sokolova2019multiscale,honorio2018stabilized} and finite element methods \cite{choo2018enriched,deng2017locally,li2020variational,Haga2012,liu2018lowest,murad2013new,wheeler2014coupling,bouklas2015effect}. Recently, the possibility of solving linear elasticity problems and coupled HM processes using physics-informed neural networks is also presented in \cite{kadeethum2020pinn,Kadeethum2020ARMA,guo2020energy,haghighat2020deep}. In the context of finite element methods, the coupled HM processes can be cast as two-, three-, four-, or five-field formulations \cite{Haga2012,phillips2007coupling1,phillips2007coupling2,kumar2020conservative,zdunek2016five}. In this study, we use the two-field formulation, in which fluid pressure and  displacement fields are the primary variables, and following adapt the discontinuous Galerkin (DG) finite element solver used in \cite{Kadeethum2019ARMA,kadeethum2019comparison,LeeKadNick_ML_2019,kadeethum2020finite,kadeethum2020enriched}. The DG solver is selected because it can handle strong heterogeneity of permeability, and it is locally mass conservative by construction \cite{riviere2008discontinuous,phillips2008coupling,liu2009coupled,choo2018enriched}.

The DG solver, also referenced as a full order model (FOM) in the following, requires substantial computational resources \cite{hansen2010discrete,hesthaven2016certified}. Hence, it is not directly suitable to handle large-scale inverse problems, optimization, or control, in which an extensive set of simulations are needed to be explored \cite{ballarin2019pod,hijazi2020effort,strazzullo2018model,hansen2010discrete,hesthaven2016certified}. Consequently, reduced order modeling (ROM) is a framework that can be employed towards handling large-scale inverse problems, optimization, or control, since ROM aims to produce a low-dimensional representation of FOM, which requires much less computational time than the FOM while maintaining computational accuracy \cite{schilders2008model,schilders2008introduction}.
The ROM methodology is applicable to a parametrized problem (i.e., to repeated evaluations of a problem depending on a parameter $\bm{\mu}$, which could correspond to physical properties, geometric characteristics, or boundary conditions \cite{ballarin2019pod,venturi2019weighted,hesthaven2016certified}).
ROM is generally composed of two stages, the offline and online stages \cite{hesthaven2016certified}. The offline stage begins with the initialization of uncertain input parameters,  which we call a training set. Then the FOM is solved corresponding to each member in the training set, and in the following, we will refer to the corresponding solutions as snapshots. Data compression techniques are then used to compress FOM snapshots to produce basis functions that span a reduced space of very low dimensionality but still guarantee accurate reproduction of the snapshots \cite{decaria2020artificial,cleary1984data,hijazi2020effort}. ROM can then be solved during the online stage for any new value of $\bm{\mu}$ by seeking an approximated solution in the reduced space.

This study focuses on utilizing a non-intrusive ROM approach \cite{wang2019non,xiao2019domain,hesthaven2016certified} to tackle the solution of strongly heterogeneous linear poroelasticity problems. This approach has not yet been extended to coupled HM processes in porous media. The non-intrusive approach is especially attractive because it does not require any cumbersome modifications of FOM source codes, and it can capture sharp changes present in the FOM approximation, which is a challenge in classical intrusive ROM models \cite{hesthaven2016certified}. This characteristic alleviates several code compatibility complications since many of the source codes used to build FOMs may not be available or easily accessible, especially in legacy codes or commercial software \cite{xiao2015non1}. Hence, the non-intrusive variant of ROM can provide more flexibility in coupling to existing FOM platforms \cite{mignolet2013review,hesthaven2018non,xiao2015non2}. The non-intrusive ROM framework presented in this manuscript is based on the development of \cite{hesthaven2018non}, which has been adapted and applied to a wide range of problems \cite{girfoglio2020nonA,girfoglio2020nonB,ortali2020gaussian,hijazi2020data,demo2020efficient,gadalla2020comparison}. In this approach, a set of reduced basis functions are extracted from high-fidelity FOM solutions through proper orthogonal decomposition (POD). The coefficients of the reduced basis functions are in turn approximated by artificial neural networks (ANN) using $\bm{\mu}$ as an input. This work has three main novelties

\begin{enumerate}

    \item The non-intrusive ROM framework has been extended to the coupled HM problem, allowing the exploration of cases where media permeability is strongly heterogeneous.

    \item We illustrate that the non-intrusive ROM solution can mimic the sharp discontinuities of both displacement and pressure fields approximated by the DG solver. This problem is challenging for classical intrusive ROM approaches.

    \item We distinguish between the error due to the truncation of reduced bases and the error introduced by the prediction of ANN. Besides, we illustrate that the nested POD technique, in which time and uncertain parameter domains are compressed consecutively, could provide comparable accuracy to the classical POD method, in which all domains are compressed simultaneously.

\end{enumerate}

The rest of the manuscript is summarized as follows. We begin with the model description and corresponding governing equations (section 2) and follow with the presentation of the discretization of the coupled HM problem and the DG solver structure (section 3). Subsequently, the ROM framework and its components are discussed (section 4). We validate our developed ROM framework using a series of benchmark problems (section 5). We then analyze the ROM performance concerning different hyperparameters and data compression strategies. We provide discussions on the sources of error, the ROM performance, and circumstances in which the ROM is more suitable than the FOM (section 6).

\section{Governing equations}\label{sec:governing_equations}

This section briefly describes governing equations following Biot's formulation for linear poroelasticity \cite{biot1941general,biot1957elastic}. Let $\Omega \subset \mathbb{R}^d$ ($d \in \{1,2,3\}$) denote the physical domain and $\partial \Omega$ denote its boundary. The time domain is denoted by $\mathbb{T} = \left(0,\mathrm{T}\right]$ with $\mathrm{T}>0$. Primary variables used in this paper are $p : \Omega \times  \mathbb{T} \to \mathbb{R}$, which is a scalar-valued fluid pressure (\si{Pa}) and $\bm{u} : \Omega \times \mathbb{T} \to \mathbb{R}^d$, which is a vector-valued displacement (\si{m}).

Although linear poroelasticity theory may oversimplify deformations in soft porous materials such as soils \cite{Choo2016,Borja2016,Macminn2016,Zhao2020}, this description is reasonably good for stiff materials such as rocks, which are the focus of this work. The theory of linear poroelasticity theory describes the HM problem through two coupled governing equations, namely linear momentum and mass balance equations.

The infinitesimal strain tensor $\bm{\varepsilon}$ is defined as
\begin{equation}
\bm{\varepsilon}(\bm{u}) :=\frac{1}{2}\left(\nabla \bm{u}+(\nabla \bm{u})^{\intercal}\right),
\end{equation}
where $\lambda_{l}$ and $\mu_{l}$ are the Lam\'e constants, which are related to the bulk modulus $K$ and the Poisson ratio $\nu$ of the porous solid as

\begin{equation}\label{eq:lambda_l}
\lambda_{l}=\frac{3 K \nu}{1+\nu}, \quad\text{ and }\quad \mu_{l}=\frac{3 K(1-2 \nu)}{2(1+\nu)}.
\end{equation}
\noindent

Further, $\bm{\sigma}$ is the total Cauchy stress tensor, which may be related to the effective Cauchy stress tensor $\bm{\sigma}^{\prime}$ and the pore pressure $p$ as
\begin{equation}
\bm{\sigma} (\bm{u},p) = \bm{\sigma}^{\prime}(\bm{u}) - \alpha p \mathbf{I}.
\end{equation}
Here, $\mathbf{I}$ is the second-order identity tensor, and
$\alpha$ is the Biot coefficient defined as \cite{jaeger2009fundamentals}:
\begin{equation} \label{eq:biot_coeff}
\alpha = 1-\frac{K}{K_{{s}}},
\end{equation}
\noindent
with {$K$} and {$K_s$} being the bulk moduli of the bulk porous material and the solid matrix, respectively.
According to linear elasticity, the effective stress tensor relates to the infinitesimal strain tensor, and therefore to the displacement, through the following constitutive relationship, which can be written as
\begin{equation}
\bm{\sigma}^{\prime}(\bm{u}) =
\lambda_{l}  \tr(\bm{\varepsilon}(\bm{u})) \mathbf{I}
+ 2 \mu_{l} \bm{\varepsilon}{(\bm{u})}.
\end{equation} 
\noindent

Under quasi-static conditions, the linear momentum balance equation can be written as
\begin{equation}
\nabla \cdot \bm{\sigma} (\bm{u},p) + \bm{f} = \bm{0},
\end{equation}
where $\bm{f}$ is the body force term defined as $\rho \phi \mathbf{g}+\rho_{s}(1-\phi) \mathbf{g}$, where $\rho$ is the fluid density, $\rho_s$ is the solid density, $\phi$ is the porosity, and $\mathbf{g}$ is the gravitational acceleration vector.
The gravitational force will be neglected in this study, but the body force term will be kept in the succeeding formulations for a more general case.

For this solid deformation problem, the domain boundary $\partial \Omega$ is assumed to be suitably decomposed into displacement and traction boundaries, $\partial \Omega_u$ and $\partial \Omega_{t}$, respectively.
Then the linear momentum balance equation is supplemented by the boundary and initial conditions as:
\begin{equation} \label{eq:linear_balance}
\begin{split}
\nabla \cdot \bm{\sigma}^{\prime}(\bm{u}) -\alpha \nabla \cdot \left(p \mathbf{I}\right)
+ \bm{f} = \bm{0}  &\text { \: in \: } \Omega \times \mathbb{T}, \\
\bm{u} =\bm{u}_{D} &\text { \: on \: } \partial \Omega_{u} \times \mathbb{T},\\
\bm{\sigma} {(\bm{u})} \cdot  \mathbf{n}=\bm{t_D} &\text { \: on \: } \partial \Omega_{t} \times \mathbb{T}, \\
\bm{u}=\bm{u}_{0}  &\text { \: in \: } \Omega \text { at } t = 0,
\end{split}
\end{equation}

\noindent
where $\bm{u}_D$ and ${\bm{t}_D}$ are prescribed displacement and traction values at the boundaries, respectively, and $\mathbf{n}$ is the unit normal vector to the boundary.

Next, the mass balance equation is given as \cite{coussy2004poromechanics,kim2011stability}:
\begin{equation} \label{eq:mass_balance_old}
\frac{1}{M} \dfrac{\partial p}{\partial t} +
\alpha  \frac{\partial {\varepsilon_{v}}}{\partial t}
+ \nabla \cdot  \bm{q} = g \text { in } \Omega \times \mathbb{T},
\end{equation}

\noindent
where
\begin{equation} \label{eq:1/M}
\frac{1}{M}  = \left(\phi c_{f}+\dfrac{\alpha-\phi}{K_{s}}\right)
\end{equation}
is the Biot modulus.
\noindent
Here, $c_f$ is the fluid compressibility, ${\varepsilon_{v}}$ := $\operatorname{tr}(\bm{\varepsilon}) = \nabla \cdot \bm{u}$ is the volumetric strain, and $g$ is a sink/source term. The superficial velocity vector $\bm{q}$ is given by Darcy's law as
\begin{equation} \label{eq:darcy}
\bm{q} =- \bm{\kappa}(\nabla p-\rho \mathbf{g}).
\end{equation}

\noindent
Here $\bm{\kappa}=\frac{\bm{k}}{\mu_f}$ is the porous media conductivity, $\mu_f$ is the fluid viscosity. Again, the gravitational force, $\rho \mathbf{g}$, will be neglected in this work, without loss of generality. In addition, $\bm{k}$ is the matrix permeability tensor defined as

\begin{equation} \label{eq:permeability_matrix}
\bm{k} :=
\begin{cases}
 \left[ \begin{array}{lll}{{k}_{xx}} & {{k}_{xy}} & {{k}_{xz}} \\ {{k}_{yx}} & {{k}_{yy}} & {{k}_{yz}} \\ {{k}_{zx}} & {{k}_{zy}} & {k}_{zz}\end{array}\right] & \text{if} \ d = 3, \ \\ \\
 \left[ \begin{array}{ll}{{k}_{xx}} & {{k}_{xy}}  \\ {{k}_{yx}} & {{k}_{yy}} \\ \end{array}\right]  & \text{if} \ d = 2, \ \\ \\
 \ {k}_{xx}  & \text{if} \ d = 1,
\end{cases}
\end{equation}

\noindent

For the fluid flow problem, the domain boundary $\partial \Omega$ is also suitably decomposed into the pressure and flux boundaries,  $\partial \Omega_p$ and $\partial \Omega_q$, respectively.
In what follows, we apply the fixed stress split scheme \cite{kim2011stability,mikelic2013convergence}, assuming $\left(\sigma_{v}-\sigma_{v, 0}\right)+\alpha \left(p-p_{0}\right)=K \varepsilon_{v}$.
Then we write the fluid flow problem with boundary and initial conditions as

\begin{equation} \label{eq:mass_balance}
\begin{split}
\left(\frac{1}{M}+\frac{\alpha^{2}}{K}\right) \frac{\partial p}{\partial t}+\frac{\alpha}{K} \frac{\partial \sigma_{v}}{\partial t}- \bm{\kappa} \nabla p = g  &\text { \: in \: } \Omega \times \mathbb{T}, \\
p=p_{D} &\text { \: on \: } \partial \Omega_{p} \times \mathbb{T}, \\
- \bm{\kappa} \nabla p \cdot \mathbf{n}=q_{D} &\text { \: on \:} \partial \Omega_{q} \times \mathbb{T}, \\
p=p_{0} &\text { \: in \: } \Omega \text { at } t = 0,
\end{split}
\end{equation}

\noindent
where $\sigma_{v}:=\frac{1}{3} \tr(\bm{\sigma})$ is the volumetric stress, and $p_D$ and $q_D$ are the given boundary pressure and flux, respectively.

\section{Finite element discretization and solution}\label{sec:fem}

We use the discontinuous Galerkin solver from \cite{Kadeethum2019ARMA,kadeethum2019comparison,LeeKadNick_ML_2019,kadeethum2020finite,kadeethum2020enriched}, and we briefly revisit the discretization in this section. We begin by introducing the necessary notation.
Let $\mathcal{T}_h$ be a shape-regular triangulation obtained by a partition of $\Omega$ into $d$-simplices (segments in $d=1$, triangles in $d=2$, tetrahedra in $d=3$). For each cell $T \in \mathcal{T}_h$, we denote by $h_{T}$ the diameter of $T$, and we set $h=\max_{T \in \mathcal{T}_h} h_{T}$ and $h_{l}=\min_{T \in \mathcal{T}_h} h_{T}$.
We further denote by $\mathcal{E}_h$ the set of all facets (i.e., $d - 1$ dimensional entities connected to at least a $T \in \mathcal{T}_h$) and by $\mathcal{E}_h^{I}$ and $\mathcal{E}_h^{\partial}$ the collection of all interior and boundary facets, respectively.
The boundary set $\mathcal{E}_h^{\partial}$ is decomposed into two disjoint subsets associated with the Dirichlet boundary facets, and the Neumann boundary facets for each of Eqs. \eqref{eq:linear_balance} and \eqref{eq:mass_balance}.
In particular, $\mathcal{E}_{h}^{D,u}$ and $\mathcal{E}_{h}^{N,u}$ correspond to the facets on $\partial \Omega_u$ and $\partial \Omega_{tr}$, respectively,  for Eq. \eqref{eq:linear_balance}.
On the other hand, for Eq. \eqref{eq:mass_balance},  $\mathcal{E}_{h}^{D,m}$ and $\mathcal{E}_{h}^{N,m}$ conform to $\partial \Omega_p$ and $\partial \Omega_{q}$, respectively.

We also define
$$
e = \partial T^{+}\cap \partial T^{-}, \ \ e \in \mathcal{E}_h^I,
$$
\noindent
where  $T^{+}$ and $T^{-}$ are the two neighboring elements to $e$. We denote by $h_e$ the characteristic length of $e$ calculated as
\begin{equation}
h_{e} :=\frac{\operatorname{meas}\left(T^{+}\right)+\operatorname{meas}\left(T^{-}\right)}{2 \operatorname{meas}(e)},
\end{equation}

\noindent
depending on the argument, meas($\cdot$) represents the measure of a cell or of a facet.

Let $\mathbf{n}^{+}$ and $\mathbf{n}^{-}$ be the outward unit normal vectors to  $\partial T^+$ and $\partial T^-$, respectively.
For any given scalar function $\zeta: \mathcal{T}_h \to \mathbb{R}$ and vector function $\bm{\tau}: \mathcal{T}_h \to \mathbb{R}^d$, we denote by $\zeta^{\pm}$ and $\bm{\tau}^{\pm}$ the restrictions of $\zeta$ and $\bm{\tau}$ to $T^\pm$, respectively.
Subsequently, we define the weighted average operator as

\begin{equation}
\{\zeta\}_{\delta e}=\delta_{e} \zeta^{+}+\left(1-\delta_{e}\right) \zeta^{-}, \ \text{ on } e \in \mathcal{E}_h^I,
\end{equation}

\noindent
and

\begin{equation}
\{\bm{\tau}\}_{\delta e}=\delta_{e} \bm{\tau}^{+}+\left(1-\delta_{e}\right) \bm{\tau}^{-},
\ \text{ on } e \in \mathcal{E}_h^I,
\end{equation}

\noindent
where $\delta_{e}$ is calculated by \cite{ErnA_StephansenA_ZuninoP-2009aa,ern2008posteriori}:

\begin{equation}
\delta_{e} :=\frac{{k}^{-}_e}{{k}^{+}_e+{k}^{-}_e}.
\end{equation}
Here,
\begin{equation}
{k}^{+}_e :=\left(\mathbf{n}^{+}\right)^{\intercal} \cdot \bm{k}^{+}  \mathbf{n}^{+}, \ \text{ and }
{k}^{-}_e :=\left(\mathbf{n}^{-}\right)^{\intercal} \cdot \bm{k}^{-}  \mathbf{n}^{-},
\end{equation}
where  ${k_e}$ is a harmonic average of $k^{+}_e$ and ${k}^{-}_e$ which reads
\begin{equation}
{k_{e}}:= \frac{2{k}^{+}_e {k}^{-}_e}{{k}^{+}_e+{k}^{-}_e},
\end{equation}
and $\bm{k}$ is defined as in Eq. \eqref{eq:permeability_matrix}.
\noindent
The jump across an interior edge will be defined as
\begin{align*}
\jump{\zeta} = \zeta^+\mathbf{n}^++\zeta^-\mathbf{n}^- \quad \mbox{ and } \quad \jtau = \bm{\tau}^+\cdot\mathbf{n}^+ + \bm{\tau}^-\cdot\mathbf{n}^- \quad \mbox{on } e\in \mathcal{E}_h^I.
\end{align*}

Finally, for $e \in \mathcal{E}^{\partial}_h$, we set $\av{\zeta}_{\delta_e} :=   \zeta$ and $\av{\bm{\tau}}_{\delta_e} :=  \bm{\tau}$ for what concerns the definition of the weighted average operator, and $\jump{\zeta} :=  \zeta \mathbf{n}$ and $\jump{\bm{\tau}} :=  \bm{\tau} \cdot \mathbf{n}$ as definition of the jump operator.

\subsection{Temporal Discretization}

We adapt the Biot's system solver from \cite{kadeethum2020enriched,kadeethum2020finite}. The time domain, $\mathbb{T} = \left(0, \mathrm{T}\right]$, is partitioned into $N^t$ open intervals such that, $0=: t^{0}<t^{1}<\cdots<t^{ N^t} := \mathrm{T}$. The length of the interval, $\Delta t^n$, is defined as $\Delta t^n=t^{n}-t^{n-1}$ where $n$ represents the current time step. $\Delta t^0$ is an initial $\Delta t$, which is defined as $t^{1}-t^{0}$, while the other time steps, $\Delta t^n$, are calculated as follows

\begin{equation} \label{eq:time_mult}
\Delta t^n :=
\begin{cases}
 \Delta t_{mult}\times \Delta t^{n-1} & \text{if} \ \Delta t^n \le \Delta t_{max} \ \\
 \Delta t_{max}  & \text{if} \ \Delta t^n > \Delta t_{max},
\end{cases}
\end{equation}

\noindent
where $\Delta t_{mult}$ is a positive constant multiplier, and $\Delta t_{max}$ is the maximum allowable time step. Then, let $\varphi(\cdot, t)$ be a scalar function and $\varphi^{n}$ be its approximation at time $t^n$, i.e. $\varphi^{n} \approx \varphi\left(t^{n}\right)$. We employ the following backward differentiation formula for time discretization of all primary variables \cite{ibrahim2007implicit,akinfenwa2013continuous,lee2018phase,kadeethum2020locally}

\begin{equation} \label{eq:bdf_gen}
\mathrm{BDF}_{1}\left(\varphi^{n}\right):=
\frac{1}{\Delta t^n}\left(\varphi^{n}-\varphi^{n-1}\right).
\end{equation}

\subsection{Full Discretization}

Following \cite{liu2009coupled,kadeethum2020enriched,kadeethum2020finite,Kadeethum2019ARMA}, in this study, the displacement field is approximated by the classical continuous Galerkin method (CG) method, and the fluid pressure field is discretized by discontinuous Galerkin (DG) method to ensure local mass conservation and provide a better flux approximation \cite{kadeethum2020enriched,kadeethum2020finite}.

We begin with defining the finite element space for the continuous Galerkin (CG) method for a vector-valued function

\begin{equation}
 \mathcal{U}_{h}^{\mathrm{CG}_{k}}\left(\mathcal{T}_{h}\right) :=\left\{\bm{\psi_u} \in \mathbb{C}^{0}(\Omega{; \mathbb{R}^d}) :\left.\bm{\psi_u}\right|_{T} \in \mathbb{Q}_{k}(T{; \mathbb{R}^d}), \forall T \in \mathcal{T}_{h}\right\},
\label{eq:CG_U}
\end{equation}

\noindent
where $k$ indicates the order of polynomial that can be approximated in this space, $\mathbb{C}^0(\Omega{; \mathbb{R}^d})$ denotes the space of vector-valued piece-wise continuous polynomials, $\mathbb{Q}_{k}(T{; \mathbb{R}^d})$ is the space of polynomials of degree at most $k$ over each element $T$, and $\mathbb{R}$ is a set of real numbers. We will denote in the following by $N_h^u$ the dimension of the space $\mathcal{U}_{h}^{\mathrm{CG}_{k}}\left(\mathcal{T}_{h}\right)$, i.e. the number of degrees of freedom for the displacement approximation.

Next, the DG space for scalar-valued functions is defined as

\begin{equation}
\mathcal{P}_{h}^{\mathrm{DG}_{k}}\left(\mathcal{T}_{h}\right) :=\left\{\psi_p \in L^{2}(\Omega) :\left.\psi_p\right|_{T} \in \mathbb{Q}_{k}(T), \forall T \in \mathcal{T}_{h}\right\},
\label{eq:DG_space}
\end{equation}

\noindent
where $L^{2}(\Omega)$ is the space of square integrable functions. This non conforming finite element space allows us to consider discontinuous coefficients and preserves the local mass conservation. We will denote in the following by $N_h^p$ the dimension of $\mathcal{P}_{h}^{\mathrm{DG}_{k}}\left(\mathcal{T}_{h}\right)$.

We seek the approximated displacement ($\bm{u}_{h}$) and pressure (${p}_{h}$) solutions by discretizing the linear momentum balance equation Eq. \eqref{eq:linear_balance} employing the above CG finite element spaces for $\bm{u}_{h}$  and the DG spaces for ${p}_{h}$.  The fully discretized linear momentum balance equation Eq. \eqref{eq:linear_balance} can be defined using the following forms

\begin{equation}
\mathcal{A}_u\left((\bm{u}_{h}^{n}, p_{h}^{n}), \bm{\psi_u} \right) = \mathcal{L}_u\left(\bm{\psi_u} \right), \quad\forall \bm{\psi_u} \in \mathcal{U}_{h}^{\mathrm{CG}_{2}}\left(\mathcal{T}_{h}\right),
\label{eq:A_u_L_u}
\end{equation}

\noindent
at each time step $t^n$, where

\begin{equation} \label{eq:linear_balance_dis_1}
\mathcal{A}_u\left((\bm{u}_{h}^{n}, p_{h}^{n}), \bm{\psi_u} \right) = \sum_{T \in \mathcal{T}_{h}} \int_{T} \bm{\sigma}^{\prime}\left(\bm{u}_{h}\right) : \nabla^{s} \bm{\psi_u} \: d V  -  \sum_{T \in \mathcal{T}_{h}} \int_{T} \alpha  p_{h} \mathbf{I}  : \nabla^{s} \bm{\psi_u} \: d V,
\end{equation}

\noindent
and

\begin{equation} \label{eq:linear_balance_dis_2}
 \mathcal{L}_u\left(\bm{\psi_u} \right) =\sum_{T \in \mathcal{T}_{h}} \int_{T} \bm{f} \bm{\psi_u} \: d V+\sum_{e \in \mathcal{E}_{h}^{N}} \int_{e} \bm{t_D} \bm{\psi_u} \: d S.
\end{equation}

\noindent
Here, $\nabla^{s}$ is a symmetric gradient operator. We then discretize Eq. \eqref{eq:mass_balance} as



\begin{equation}
\mathcal{A}_p\left((\bm{u}_{h}^{n}, p_{h}^{n}), \psi_p \right) = \mathcal{L}_p\left(\psi_p \right), \quad\forall \psi_p \in \mathcal{P}_{h}^{\mathrm{DG}_{1}}\left(\mathcal{T}_{h}\right),
\label{eq:A_p_L_p}
\end{equation}

\noindent
for each time step $t^n$, where

\begin{equation} \label{eq:mass_discretization_dg_1}
\begin{split}
\mathcal{A}_p\left((\bm{u}_{h}^{n}, p_{h}^{n}), \psi_p \right) & =
\sum_{T \in \mathcal{T}_{h}} \int_{T} \frac{\alpha}{K} \mathrm{BDF}_1(\sigma_{v}) \psi_p \: d V + \sum_{T \in \mathcal{T}_{h}} \int_{T}  \left(\frac{1}{M}+\frac{\alpha^{2}}{K}\right)
\mathrm{BDF}_{1}\left( p_{h}^{n} \right) \psi_p \: d V \\
& + \sum_{T \in \mathcal{T}_{h}} \int_{T} \bm{\kappa} \nabla p_{h}^n \cdot \nabla \psi_{p} \: d V  - \sum_{e \in \mathcal{E}_h^{I} \cup \mathcal{E}_{h}^{D}}  \int_{e}\left\{\bm{\kappa}\nabla p_{h}^n\right\}_{\delta_{e}} \cdot \llbracket \psi_p \rrbracket \: d S \\ & - \sum_{e \in \mathcal{E}_h^{I} \cup \mathcal{E}_{h}^{D}}  \int_{e}\left\{\bm{\kappa} \nabla \psi_{p}\right\}_{\delta_{e}} \cdot \llbracket p_h^n \rrbracket \: d S  +  \sum_{e \in \mathcal{E}_h^{I} \cup \mathcal{E}_{h}^{D}}  \int_{e} \frac{\beta}{h_{e}} {\bm{\kappa}}_{{e}}  \llbracket p_h^n \rrbracket \cdot \llbracket \psi_p \rrbracket \: d S,
\end{split}
\end{equation}

\noindent
and

\begin{equation}
\begin{split}
\mathcal{L}_p\left(\psi_p  \right) & :=  \sum_{T \in \mathcal{T}_{h}} \int_{T} g \psi_{p} \: d V+ \sum_{e \in \mathcal{E}_{h}^{N}} \int_{e} q_{D} \psi_{p} \: d S \\ &  -\sum_{e \in \mathcal{E}_{h}^{D}} \int_{e} \bm{\kappa} \nabla \psi_{p} \cdot p_{D} \mathbf{n} \: d S + \sum_{e \in \mathcal{E}_{h}^{D}} \int_{e} \frac{\beta}{h_{e}} {\bm{\kappa}}_{{e}} \llbracket \psi_p \rrbracket \cdot p_D \mathbf{n} \: d S.
\end{split}
\label{eq:L_p}
\end{equation}

More details regarding block structure and solver algorithm could be found in \cite{kadeethum2020finite,kadeethum2020enriched}. For all the computations, matrices and vectors are built using the FEniCS form compiler \cite{AlnaesBlechta2015a}. The block structure is assembled by using the multiphenics toolbox \cite{Ballarin2019}. Solvers are employed from PETSc package \cite{petsc-user-ref}. All simulations are computed on AMD Ryzen Threadripper 3970X with a single thread.







\section{Reduced order modeling}

\begin{figure}[!ht]
   \centering

         \includegraphics[keepaspectratio, height=6.0cm]{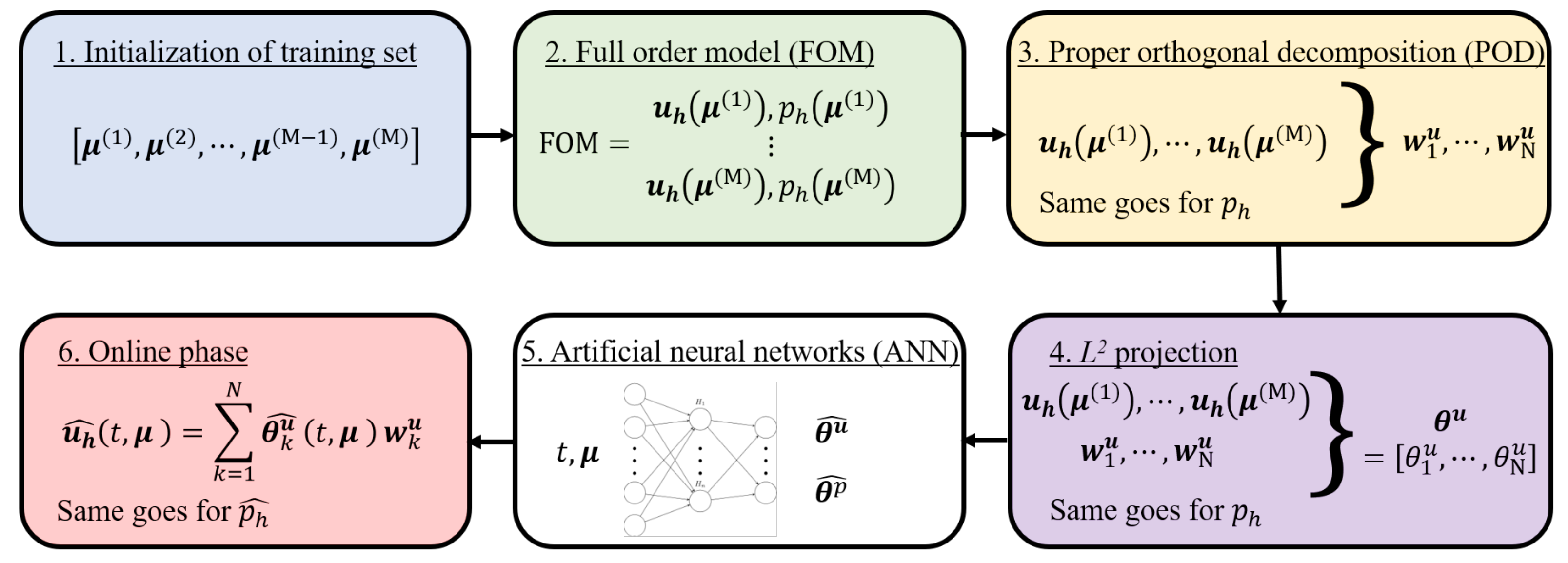}

   \caption{Summary of non-intrusive model reduction framework for Biot's system. We note that, for the sake of simplicity, we only show $\bm{\theta}$ and $\widehat{\bm{\theta}}$ in this figure. They; however, represent $\bm{\theta} = [\bm{\theta}^u(t, \bm{\mu}), \bm{\theta}^p(t, \bm{\mu})]$ and $\widehat{\bm{\theta}} = [\widehat{\bm{\theta}}^u(t, \bm{\mu}), \widehat{\bm{\theta}}^p(t, \bm{\mu})]$ }
   \label{fig:pod_nn_explain}
   \end{figure}

The DG solution scheme introduced in the previous section is typically a time-consuming operation, making it impractical to query such a solver in a real-time context whenever parametric studies are carried out. Such parametric studies are often of interest to account for uncertain material properties of porous media. Therefore in this work, we propose to employ a reduced order modeling strategy, based on the developments in \cite{hesthaven2018non}. A graphical summary of the reduced order modeling paradigm is presented in Figure \ref{fig:pod_nn_explain}.
The idea of this framework has been adapted and applied to a wide range of problems \cite{girfoglio2020nonA,girfoglio2020nonB,ortali2020gaussian,hijazi2020data,demo2020efficient,gadalla2020comparison}. The computations are divided into an offline phase for the ROM construction, which we will present through five consecutive steps, and (single-step) online stage for the ROM evaluation, described as follows.

The first step of the offline stage (colored in blue in the Figure) represents an initialization of a training set of parameters used to train the framework, of cardinality $\mathrm{M}$.
Then, in the second step (green), we query the full order model (FOM), based on the DG finite element solver discussed in the previous section, for each parameter $\bm{\mu}$ in the training set. At this point, we have $\mathrm{M}$ snapshots (FOM results) associated with the different parametric configurations in the training set. Each snapshot contains approximations of the primary variables ($\bm{u}_{h}$ and $p_h$) at each time step of the partition of the time domain $\mathbb{T}$ as introduced in the previous section.
The third step (yellow) aims to compress the information provided by the snapshots through the proper orthogonal decomposition (POD) technique \cite{hesthaven2016certified,wang2020recurrent,hesthaven2018non,hijazi2020data,stabile2017pod,willcox2002balanced}. The POD is used to determine characteristic spatial modes based on relative energy content criteria \cite{chatterjee2000introduction,liang2002proper,wang2018model}. In order to carry out a compression, only the first $\mathrm{N}$ spatial modes are retained \cite{hesthaven2016certified}, and employed as basis functions for the reduced basis spaces $\mathcal{U}_\mathrm{N}$ and $\mathcal{P}_\mathrm{N}$, used for approximating the displacement and pressure fields respectively. The typical goal is to achieve $\mathrm{N} \ll \mathrm{M} N^t$ (compression of the snapshots data), but also $\mathrm{N} \ll N_h^u$ and $\mathrm{N} \ll N_h^p$ (dimensionality reduction for the model discretization).
Next, in the fourth step (purple) we obtain the optimal representation of each snapshot in the reduced basis spaces by means of an $L^2$ projection \cite{girfoglio2020nonA,girfoglio2020nonB,ortali2020gaussian,hijazi2020data,demo2020efficient,gadalla2020comparison}.
This operation defines a map between each pair $(t, \bm{\mu})$, with $t \in \{t^0, \hdots, t^{N^t}\}$ and $\bm{\mu}$ in the training set, and a vector of coefficients $\bm{\theta}^u(t, \bm{\mu})$ that characterize the best approximation in the reduced space $\mathcal{U}_\mathrm{N}$ for the displacement field $\bm{u}_h(\bm{\mu})$ at time $t$. A similar map can be defined for the pressure field, denoted in the following by $\bm{\theta}^p(t, \bm{\mu})$.
Finally, in the fifth step (white) we aim to define a map $\widehat{\bm{\theta}}^u(t, \bm{\mu})$ and $\widehat{\bm{\theta}}^p(t, \bm{\mu})$ for any time $t$ in the time interval $\mathbb{T}$ and any value of the parameter $\bm{\mu}$ by training artificial neural networks (ANN) to approximate $\bm{\theta}^u(t, \bm{\mu})$ and $\bm{\theta}^p(t, \bm{\mu})$ based on the training data points obtained at the fourth step \cite{hesthaven2018non}. We note that the wall time used to perform the $L^2$ projection during the fourth step is relatively much smaller compared to the wall time used to train the ANN in the fifth step. This concludes the offline stage.

Finally, during the online phase (red), for given values of the parameter $\bm{\mu}$ and time instance $t$ we aim to recover the online approximation to our primary variables by querying the ANN evaluation for $\widehat{\bm{\theta}}^u(t, \bm{\mu})$ and $\widehat{\bm{\theta}}^p(t, \bm{\mu})$ and reconstructing the resulting finite element representation by means of the reduced basis functions spanning $\mathcal{U}_\mathrm{N}$ and $\mathcal{P}_\mathrm{N}$, respectively \cite{hesthaven2018non}. The details of each phase will be further discussed in the following paragraphs.

\subsection{Initialization of the training set} \label{sec:init}

Let $\mathbb{P} \subset \mathbb{R}^P$, $P \in \mathbb{N}$, be a compact set representing the range of variation of the parameters $\bm{\mu} \in \mathbb{P}$.
For the sake of notation we will denote by $\mu_p$, $p = 1, \hdots, P$, the $p$-th component of $\bm{\mu}$.
To explore the parametric dependence of the phenomena, we define a discrete training set of $\mathrm{M}$ parameter instances. Each parameter instance in the training set will be indicated with the notation $\bm{\mu}^{(i)}$, for $i = 1, \hdots, \mathrm{M}$. Thus, the $p$-th component of the $i$-th parameter instance in the training set is denoted by $\mu_p^{(i)}$ in the following. The choice of the value of $\mathrm{M}$, as well as the sampling procedure from the range $\mathbb{P}$ is typically user- and problem-dependent. In this work, we use an equispaced distribution for the training set, and we will also briefly discuss the robustness of numerical results for varying $\mathrm{M}$.
We note that adaptive sampling approaches could be employed and might result in a better model accuracy with a lower number of training instances \cite{paul2015adaptive,vasile2013adaptive}. Still, we prefer here to use a straightforward equispaced sampling as a thorough discussion of the parameter space exploration is not among the main goals of this work.

\subsection{Full order model (FOM)} \label{sec:FOM_in_ROM}

We employ the DG scheme discussed in Section \ref{sec:fem} as our full order model (FOM). The FOM finite element solver is operated $\mathrm{M}$ times, corresponding to each parameter instance of $\bm{\mu}^{(i)}$. As we develop our solver on top of the FEniCS platform \cite{AlnaesBlechta2015a}, this part could be performed in parallel using any desirable processor numbers. In this study, however, we perform our FOM simulation using only a single core and run each FOM snapshot sequentially to provide a clear comparison of wall time used to perform the FOM and the wall time used to construct the ROM solution in the online phase. Since the problem formulation is time-dependent, the output of the FOM solver for each parameter instance $\bm{\mu}^{(i)}$ collects the time series representing the time evolution of the primary variables for each time step $t$. Therefore, based on the training set cardinality $\mathrm{M}$ and the number $N^t$ of time steps, we have a total of $N^t \mathrm{M}$ training data to be employed in the subsequent steps.

\subsection{Proper orthogonal decomposition (POD)}

In this work, we utilize POD as a data compression tool, i.e. we seek a reduced order approximation in an optimal linear subspace \cite{hesthaven2016certified,wang2020recurrent,hesthaven2018non,hijazi2020data,stabile2017pod,willcox2002balanced}. If the problem does not allow such representation, nonlinear variants (e.g., autoencoders) could be considered as data compression tools \cite{hinton1994autoencoders,phillips2020autoencoder,o2019learning,goh2019solving}. Here, we prefer to employ POD because it is generally faster than the nonlinear variants. As the numerical results will show, POD spaces provide ROM results of sufficient accuracy for the problem at hand.


Let $\bm{\mu}^{(i)}$ be a parameter instance in the training set, $i = 1, \hdots, \mathrm{M}$. The corresponding displacement field snapshot contains
\begin{equation}
{\mathbb{S}_u^{(i)}}=\left[{\bm{u}}_{h}\left({\cdot; t^{0}, \bm{\mu}^{(i)}}\right), \cdots, {\bm{u}}_{h}\left({\cdot; t^{{N^t}}, \bm{\mu}^{(i)}}\right)\right] \in \mathbb{R}^{{N_{h}^u} \times N^t}. 
\end{equation}
where ${\bm{u}}_{h}\left({\cdot; t^{n}, \bm{\mu}^{(i)}}\right)$ represent the displacement field at time $t^{n}$ and parameter instance $\bm{\mu}^{(i)}$. We recall that $N_{h}^u$ is the number of DOFs in the displacement finite element space, and $N^t$ is the total number of time steps.
In this study, $N_{h}^u$ is constant (i.e., the mesh and finite element function space remain the same), and $N^t$ is fixed (i.e., each snapshot utilizes the same initial and final time, and time step.).

We compare in the rest of the paper two variants of POD-based compression for the set of snapshots ${\mathbb{S}_u^{(i)}}$, $i = 1, \hdots, \mathrm{M}$. For compactness of exposition in the rest of this section, we will focus on the displacement field $\bm{u}_{h}$, but a very similar procedure is indeed carried out for the pressure field $p_h$ as well. Apart from primary variables $\bm{u}_{h}$ and $p_h$, compression of any other quantity of interest (e.g., the fluid flux at internal and external boundaries or the maximum total stress) could also be carried out using ROM. While this may be of great interest in applications, we are just interested in reducing the primary variables as our focus is to validate the methodology.

The first choice is based on a \emph{standard} POD algorithm where all snapshots are compressed in a single procedure. We first collect all snapshots in a matrix
\begin{equation}
{\mathbb{S}_u}=\left[{\mathbb{S}_u^{(1)}, \cdots, \mathbb{S}_u^{(\mathrm{M})}}\right] \in \mathbb{R}^{{N_{h}^u} \times N^t\mathrm{M}},
\end{equation}
by horizontally stacking all matrices $\mathbb{S}_u^{(i)}$, $i = 1, \hdots, \mathrm{M}$. We then perform the singular value decomposition (SVD) of ${\mathbb{S}_u}$ as

\begin{equation}
{\mathbb{S}_u}=\mathbb{W}\left[\begin{array}{cc}
\mathbb{D} & 0 \\
0 & 0
\end{array}\right] \mathbb{Z}^{\top}
\end{equation}

\noindent
where $\mathbb{W}=\left[\mathbf{w}_{1},\cdots, \mathbf{w}_{{N_{h}^u}}\right] \in \mathbb{R}^{{N_{h}^u} \times {N_{h}^u}}$ and $\mathbb{Z}=\left[\mathbf{z}_{1},\cdots, \mathbf{z}_{N^t \mathrm{M}}\right] \in \mathbb{R}^{N^t \mathrm{M}\times N^t \mathrm{M}}$ are orthogonal matrices, $\mathbb{D}=\operatorname{diag}\left(\sigma_{1}, \cdots, \sigma_{r}\right) \in \mathbb{R}^{r \times r}$
is a diagonal matrix, with singular values $\sigma_{1} \geq \sigma_{2} \geq \cdots \geq \sigma_{r}>0 .$ Here, $r$ is the number of non-zero singular values and $r \leq \min \left\{{N_{h}^u}, N^t \mathrm{M}\right\}$. The columns of $\mathbb{W}$ are called left singular vectors of $\mathbb{S},$ and the columns of $\mathbb{Z}$ are called right singular vectors of $\mathbb{S}$.
To carry out a dimensionality reduction, the POD basis of rank $\mathrm{N} \ll r$ consisting of the first $\mathrm{N}$ left singular vectors of $\mathbb{S}$, and it has the property of minimizing the projection error defined by

\begin{equation}
\left\{\mathbf{w}_{1}, \cdots, \mathbf{w}_{\mathrm{N}}\right\} = \arg\min \left\{
\varepsilon\left(\tilde{\mathbf{w}}_{1}, \cdots, \tilde{\mathbf{w}}_{\mathrm{N}}\right)={\sum_{i=1}^{\mathrm{M}}\sum_{k=0}^{N^t} \left\|\bm{u}_{h}\left(\cdot; t^{k}, \bm{\mu}^{(i)}\right)-\sum_{n=1}^{\mathrm{N}}\left(\bm{u}_{h}\left(\cdot; t^{k}, \bm{\mu}^{(i)}\right), \tilde{\mathbf{w}}_{n}\right)_{u} \tilde{\mathbf{w}}_{n}\right\|_{u}^{2}}
\right\}
\end{equation}

\noindent
among all the orthonormal bases $\left\{\tilde{\mathbf{w}}_{1}, \cdots, \tilde{\mathbf{w}}_{\mathrm{N}}\right\} {\subset \mathbb{R}^{N_{h}^u}}$. Here $(\cdot, \cdot)_u$ denotes an inner product for the displacement space, while $\left\|\cdot\right\|_u$ its induced norm.
The reduced basis space $\mathcal{U}_{\mathrm{N}}$ is then defined as the span of $\left\{\mathbf{w}_{1}, \cdots, \mathbf{w}_{\mathrm{N}}\right\}$. The effect of the dimension $\mathrm{N}$ on the accuracy of the resulting ROM will be discussed in the numerical examples section (see Section \ref{sec:results}).

The second choice is instead based on a \emph{nested} POD algorithm. The primary rationale for this choice is that the SVD computation of the standard POD case may become unfeasible when $\mathrm{M}$ is large (i.e., finer sampling of the parameter space) and $N^t$ is large as well (i.e., small-time steps). Indeed, the SVD of a matrix with a large number $N^t\mathrm{M}$ of columns may require a large amount of resources, both in terms of CPU time and memory storage. The bottleneck is due to the simultaneous compression in parameter space and time. The nested POD algorithm, instead, aims at decoupling the compression in consecutive stages, operating only either on the time interval or on the parameter space. Similar algorithms are often employed by practitioners in the reduced order modeling community and can be found in the literature with various names, such as two-level POD, or hierarchical approximate POD \cite{Audouze1,rapun2010reduced,Audouze2,ballarin2016fast,himpe2018hierarchical,wang2019non,jacquier2020non}. The nested POD algorithm can be summarized in the two following sequential stages:
\begin{enumerate}
\item[1)] \emph{compression on the temporal evolution}: for each parameter instance $\bm{\mu}^{(i)}$ in the training set compress the temporal evolution stored in $\mathbb{S}_u^{(i)} \in \mathbb{R}^{N_{h}^u \times N^t}$ by means of a POD, retaining only the first $\mathrm{N}_{\mathrm{int}} \ll N^t$ modes. A compressed matrix $\widetilde{\mathbb{S}}_u^{(i)} \in \mathbb{R}^{N_{h}^u \times \mathrm{N}_{\mathrm{int}}}$ is the assembled by storing by column the first $\mathrm{N}_{\mathrm{int}}$ modes, scaled by the respective singular values. The value of $\mathrm{N}_{\mathrm{int}}$ can be chosen according to energy criteria (and thus it will be, in general, depending on the index $i$) or can be fixed a priori (as we do in this study), and is typically considerably smaller than the number of time steps $N^t$.
\item[2)] \emph{compression on the parameter space}: after the temporal evolution of each parameter instance has been compressed, one can assemble the following matrix
\begin{equation}
{\widetilde{\mathbb{S}}_u}=\left[{\widetilde{\mathbb{S}}_u^{(1)}, \cdots, \widetilde{\mathbb{S}}_u^{(\mathrm{M})}}\right] \in \mathbb{R}^{{N_{h}^u} \times \mathrm{N}_{\mathrm{int}}\mathrm{M}}.
\end{equation}
One can proceed as in the standard POD and define the reduced basis space $\mathcal{U}_{\mathrm{N}}$ obtained after compression of ${\widetilde{\mathbb{S}}_u}$. Note that the final goal of the nested POD is still to obtain a reduced basis space of dimension $\mathrm{N}$, which is computed from an SVD of a matrix with $\mathrm{N}_{\mathrm{int}}\mathrm{M} \ll N^t\mathrm{M}$ columns, thus overcoming the bottleneck of the first algorithm.
\end{enumerate}

\begin{figure}[!ht]
   \centering

         \includegraphics[keepaspectratio, height=17.0cm]{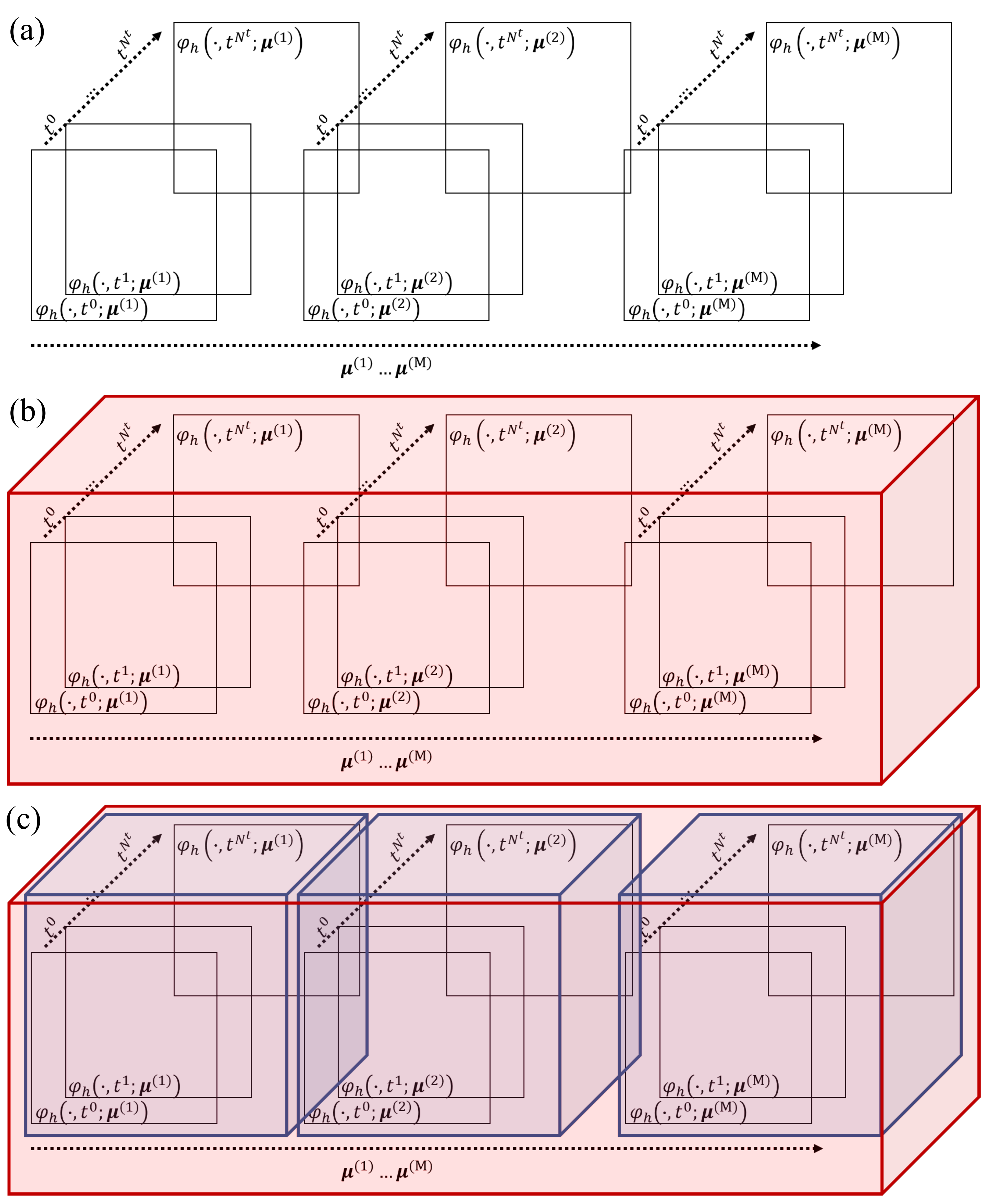}

   \caption{Proper orthogonal decomposition (POD) variants used in this study: (a) input data based on $\mathrm{M}$ training instances and ${N^t}$ time steps, (b) standard POD, and (c) nested POD. Each colored box represents a compression by SVD.}
   \label{fig:pod_explain}
   \end{figure}

We summarize the computations required by each of the two algorithms in Figure \ref{fig:pod_explain}.
Figure \ref{fig:pod_explain}a reports the input data to the two algorithms, namely parameters $\bm{\mu}^{(i)}$ in the training set, time steps $t^{n}$, and corresponding displacement or pressure fields obtained querying the DG solver. A generic field $\varphi_h$ is shown to serve as a reminder that the compression is carried out for both primary variables to obtain reduced spaces $\mathcal{U}_{\mathrm{N}}$ and $\mathcal{P}_{\mathrm{N}}$.
In the first approach (see Figure \ref{fig:pod_explain}b), we perform a compression over the whole matrix $\mathbb{S}$; a reduction by a SVD is represented in the picture by means of a colored box. In the second approach (see Figure \ref{fig:pod_explain}c), we utilize a nested POD method instead; compressions on the temporal evolution are represented by a blue box, while the final compression on the parameter space is depicted by a red box.
We finally note that, due to the adopted scaling in 1), the standard POD is formally equivalent to a nested POD algorithm with $\mathrm{N_{int}} = N^t$. However, it would be impractical to carry out the standard POD in such a manner because it would require $\mathrm{M}$  intermediate compressions without resolving the underlying bottleneck. Still, this formal equivalence motivates us to present numerical results for the standard POD with the label $\mathrm{N_{int}} = \infty$, where the symbol $\infty$ (instead of the actual value $N^t$) serves us as a reminder that no intermediate compressions are carried out.

\subsection{$L^2$ projection}

Again, for the sake of compactness, we will focus on $\bm{u}_{h}$, but a similar procedure is carried out for $p_h$. Let $\left\{\mathbf{w}_{1}, \cdots, \mathbf{w}_{\mathrm{N}}\right\}$ denote the basis functions spanning $\mathcal{U}_{\mathrm{N}}$. Given a time $t^{n}$ in the discretization of the time interval $\mathbb{T}$ and a parameter instance $\bm{\mu}^{(i)}$ in the training set we can define the best approximation $\widetilde{\bm{u}}_{h}\left(\cdot; t^{n}, \bm{\mu}^{(i)}\right)$ to $\bm{u}_{h}\left(\cdot; t^{n}, \bm{\mu}^{(i)}\right)$ in $\mathcal{U}_{\mathrm{N}}$ as
\begin{equation}
\widetilde{\bm{u}}_{h}\left(\cdot; t^{n}, \bm{\mu}^{(i)}\right) = \sum_{k=1}^{\mathrm{N}} {\theta}_k^u(t^n, \bm{\mu}^{(i)}) \mathbf{w}_{k}
\end{equation}
Here we collect in the vector $\bm{\theta}^u(t^{n}, \bm{\mu}^{(i)})=\left[\theta_{1}^u(t^{n}, \bm{\mu}^{(i)}), \cdots, \theta_{\mathrm{N}}^u(t^{n}, \bm{\mu}^{(i)})\right] \in \mathbb{R}^{\mathrm{N}}$, where the coefficients $\theta_j^u$ are solutions to the $L^2$ projection problem, which can be stated as: Given $\bm{u}_{h}\left(\cdot; t^{n}, \bm{\mu}^{(i)}\right)$, find $\bm{\theta}^u(t^{n}, \bm{\mu}^{(i)})$ such that: $\sum_{j=1}^{\mathrm{N}} {\theta}_j^u(t, \bm{\mu}) (\mathbf{w}_{j}, \mathbf{w}_{k})_u =
(\bm{u}_{h}\left(\cdot; t^{n}, \bm{\mu}^{(i)}\right), \mathbf{w}_{k})_u,  k = 1, \hdots, \mathrm{N}.$ We note that this results in a linear system, which left-hand side $(\mathbf{w}_{j}, \mathbf{w}_{k})_u$ can be easily precomputed and stored in a $\mathrm{N} \times \mathrm{N}$ matrix. However, the right-hand side can only be computed once the DG solutions are available for the training set and corresponding time steps. The goal of the next subsection is to generalize the computation of the coefficients of the ROM expansion for any (time, parameter) pair using artificial neural networks trained on the available data $\bm{\theta}^u(t^{n}, \bm{\mu}^{(i)})$.



\subsection{Artificial neural networks (ANN)}\label{sec:ann}

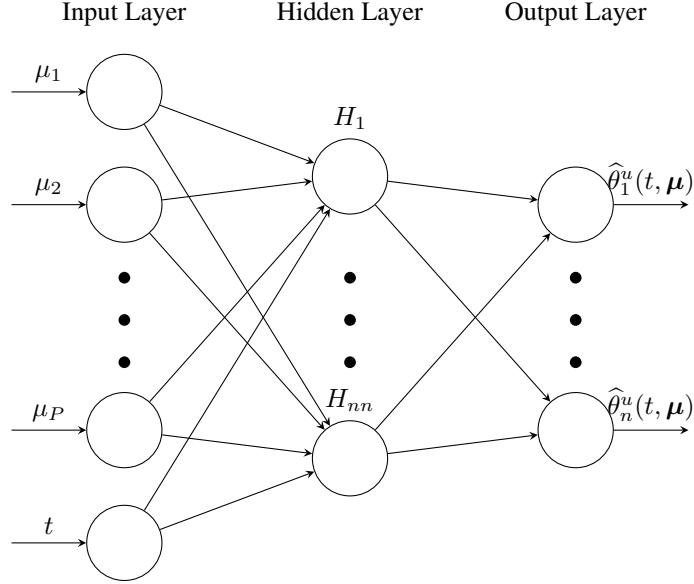
\begin{figure}[!ht]
   \centering
        \tikzset{%
  every neuron/.style={
    circle,
    draw,
    minimum size=1cm
  },
  neuron missing/.style={
    draw=none,
    scale=4,
    text height=0.333cm,
    execute at begin node=\color{black}$\vdots$
  },
}
\begin{tikzpicture}[x=1.5cm, y=1.5cm, >=stealth]

\foreach \m/\l [count=\y] in {1,2,missing,3,4}
  \node [every neuron/.try, neuron \m/.try] (input-\m) at (0,2.5-\y) {};

\foreach \m [count=\y] in {1,missing,2}
  \node [every neuron/.try, neuron \m/.try ] (hidden-\m) at (2,2-\y*1.25) {};

\foreach \m [count=\y] in {1,missing,2}
  \node [every neuron/.try, neuron \m/.try ] (output-\m) at (4,1.5-\y) {};


\draw [<-] (input-1) -- ++(-1,0)
    node [above, midway] {$\mu_1$};

\draw [<-] (input-2) -- ++(-1,0)
    node [above, midway] {$\mu_2$};

\draw [<-] (input-3) -- ++(-1,0)
    node [above, midway] {$\mu_P$};

\draw [<-] (input-4) -- ++(-1,0)
    node [above, midway] {$t$};

\foreach \l [count=\i] in {1,nn}
  \node [above] at (hidden-\i.north) {$H_{\l}$};

\foreach \l [count=\i] in {1,n}
  \draw [->] (output-\i) -- ++(1,0)
    node [above, midway] {$\widehat{\theta}_\l^u(t, \bm{\mu})$};

\foreach \i in {1,...,4}
  \foreach \j in {1,...,2}
    \draw [->] (input-\i) -- (hidden-\j);

\foreach \i in {1,...,2}
  \foreach \j in {1,...,2}
    \draw [->] (hidden-\i) -- (output-\j);

\node [align=center, above] at (0,2) {Input Layer};
\node [align=center, above] at (2,2) {Hidden Layer};
\node [align=center, above] at (4,2) {Output Layer};

\end{tikzpicture}
   \caption{{
Artificial neural network architecture used in this study. Input layer contains up to $P+1$ input nodes ($\mu_1$ to $\mu_P$ and $t$), and output layer is composed of $1,...,\mathrm{N}$ output nodes ($\widehat{\theta}_1^u(t, \bm{\mu})$ to $\widehat{\theta}_\mathrm{N}^u(t, \bm{\mu})$). The number of hidden layers are denoted by $\mathrm{N_{hl}}$ and each hidden layer is composed of $\mathrm{N_{nn}}$ neurons ($H_1$ to $H_{nn}$).}}
   \label{fig:gerneral_nn}
\end{figure}

Following the determination of the correspondence between $(t^{n}, \bm{\mu}^{(i)})$ and ${\bm{\theta}^u(t^{n}, \bm{\mu}^{(i)})}$, we aim to construct artificial neural networks (ANN) to map an input space of $\mathbb{T} \times \mathbb{P} \ni (t, \bm{\mu})$ to a vector of coefficients $\widehat{\bm{\theta}}^u(t, \bm{\mu})$ that reproduce the training data. The network architecture used in this work is presented in Figure \ref{fig:gerneral_nn}. The number of hidden layers ($\mathrm{N_{hl}}$) and the number of neurons ($\mathrm{N_n}$) act as so-called  hyperparameters \cite{goodfellow2016deep}. Each neuron (e.g., ${H_{1, 1}}$ $...$ ${H_{1, \mathrm{N_n}}}$) is connected to the nodes of the previous layer with adjustable weights and also has an adjustable bias. We denote the set of weights and biases as ($\mathrm{W}$) and ($\mathrm{b}$), respectively. These variables are learned during a training phase \cite{hinton2006reducing, goodfellow2016deep}. The neural networks are built on the PyTorch platform \cite{NEURIPS2019_9015}. The results produced using either the rectified linear unit (ReLU) or the hyperbolic tangent ($\tanh$) were comparable. Hence, we only present the results using the $\tanh$ activation function in this paper.

\noindent
Here we use a mean squared error ($\mathrm{MSE}^{\theta, u}$) as the metric of our network loss function, defined as follows

\begin{equation}\label{eq:loss_mse}
{\mathrm{MSE}^{\theta, u}}=\frac{1}{\mathrm{M} N^t} \sum_{i=1}^{\mathrm{M}}\sum_{k=0}^{N^t}\left|\widehat{\bm{\theta}}^u\left(t^k, \bm{\mu}^{(i)}\right)-{\bm{\theta}^u}\left(t^k, \bm{\mu}^{(i)}\right)\right|^{2}.
\end{equation}

\noindent
To minimize Eq. \eqref{eq:loss_mse}, we train the neural network using the adaptive moment estimation (ADAM) algorithm \cite{kingma2014adam}. Throughout this study, we use a batch size of 32, a learning rate of 0.001, a number of epoch of 20,000, and we normalize both our input and output to $[0, 1]$. To prevent our networks from overfitting behavior, we follow early stopping and generalized cross-validation criteria \cite{hesthaven2018non,prechelt1998early,prechelt1998automatic}. Note that instead of literally stopping our training cycle, we only save the set of trained $\mathrm{W}$ and $\mathrm{b}$ to be used in the online phase when the current validation loss is lower than the lowest validation from all the previous training cycle. This procedure ensures we compare our ANN training time with a fixed number of epochs. As already noted in the two previous subsections, we train the ANN specifically for each primal variable.

\subsection{Online phase}

During the online phase, for each inquiry (i.e., a novel value of $\bm{\mu}$), we evaluate the ANN to obtain $\widehat{\bm{\theta}}^u(t, \bm{\mu})$ for each $t \in \{t^0, \cdots, t^{{N^t}}\}$). Subsequently, we reconstruct the displacement as

\begin{equation}
\widehat{\bm{u}}_{h}\left(\cdot; t, \bm{\mu}\right) = \sum_{k=1}^{\mathrm{N}} \widehat{\theta}_k^u(t, \bm{\mu}) \mathbf{w}_{k},
\label{eq:online_solution}
\end{equation}

\noindent
and similarly for the pressure.
We note that the reduced basis $\{\mathbf{w}_{k}\}_{k=1}^{\mathrm{N}}$ is already constructed during the POD phase; hence, recovering the online solutions, requires to evaluate $\widehat{\theta}_k^u(t, \bm{\mu})$ from the trained ANN (which is typically extremely fast), and subsequently, perform a reconstruction using Eq. \eqref{eq:online_solution} (which only requires a linear combination of finite element functions). As a result, one typically enjoys an inexpensive online phase for each inquiry.

\section{Numerical Examples}\label{sec:results}

Throughout this section we take $\Omega = \left(0, 1\right)^2$ corresponding to a square domain of $1\mathrm{m}^2$ area, and decompose its boundary $\partial \Omega$ with the following labels

\begin{equation}
\begin{aligned}
\mathrm{Left} & = \{0\} \times [0,1]  \\
\mathrm{Top} &= [0,1] \times \{1\}   \\
\mathrm{Right} &= \{1\} \times[0,1]   \\
\mathrm{Bottom} &= [0,1] \times \{0\}.
\end{aligned}
\end{equation}

\begin{figure}[!ht]
   \centering
    \includegraphics[width=6.0cm,keepaspectratio]{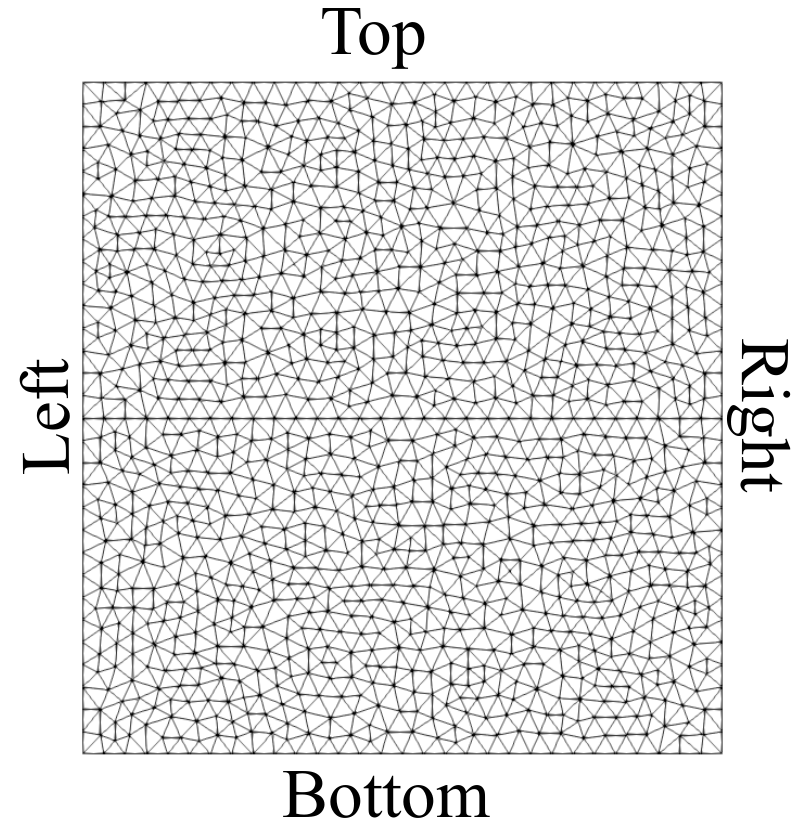} 
   \caption{Domain, its boundaries, and mesh used for all numerical examples}
   \label{fig:mesh}
\end{figure}

\noindent
A plot of the domain, its boundary labels, and the mesh we utilized for its discretization is shown in Figure \ref{fig:mesh}. The mesh contains 2370 elements, and its maximum size $h$ is $0.047$ $m$. We note that mesh is split into two conforming subdomains $(0, 1) \times (0, 0.5)$ and $(0, 1) \times (0.5, 1)$ because one of the test cases presented in this section will employ different material properties in the two subdomains. The degrees of freedom associated with this mesh are $9722$ for the continuous approximation of the displacement field $\bm{u}$, and $7110$ for the discontinuous approximation of the pressure $p$. For what concerns the time discretization, we choose $\Delta t^0= 20.0\,\mathrm{s}$, $\Delta t_{mult}= 1.0\,\mathrm{s}$, and $\Delta t_{max} = 20\,\mathrm{s}$. In each of the following subsections, we will specify the input parameters and boundary and initial conditions for each considered test case.

\subsection{Model validation}

In this subsection, we verify the developed reduced order modeling framework through a series of benchmark problems. Here we fix $\mathrm{N} = 5$, $\mathrm{N_{int}} = 5$, $\mathrm{N_{hl}} = 3$, and $\mathrm{N_{nn}} = 7$ to simplify the presentation, since the goal of this subsection is to showcase the versatility of the proposed framework for test cases of increasing difficulty. Furthermore, the effects (e.g., in terms of training time and model accuracy) of each of the aforementioned hyperparameters will be discussed later in Section \ref{sec:analysis}.

\subsubsection{Example 1: Terzaghi's consolidation problem} \label{sec:1d_1l}

\begin{figure}[!ht]
   \centering
    \includegraphics[width=6.0cm,keepaspectratio]{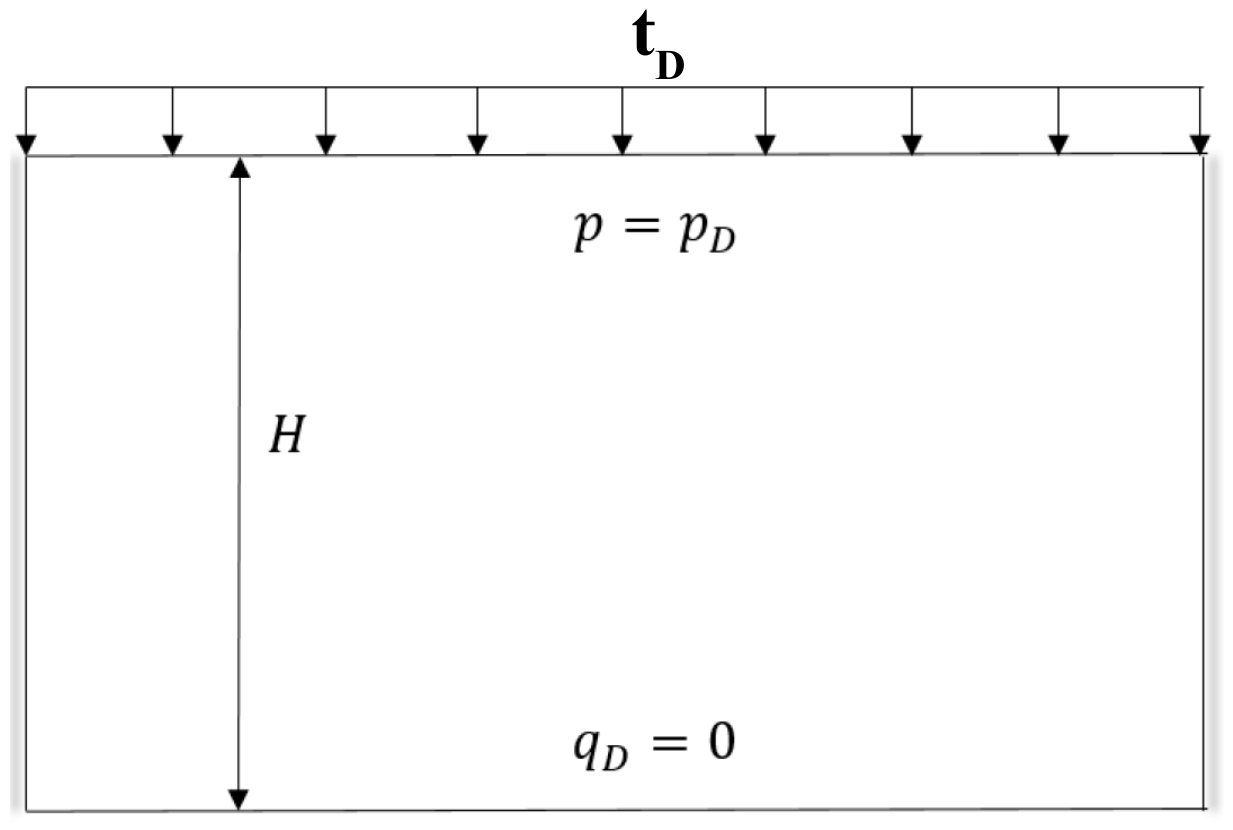} 
   \caption{Example 1: Setup for the 1-dimensional consolidation problem in a homogeneous material. Here $H = 1, {p_D} = 0$ and $\bm{t_D} =[0, -1]  \quad  \si{kPa}$. This figure is adapted from \cite{kadeethum2020finite}.}
   \label{fig:validation_geo}
\end{figure}

We first verify the presented reduced order model by means of the test case that we have already employed in the validation of the finite element solver in \cite{Kadeethum2019ARMA,kadeethum2020finite,kadeethum2020enriched}. This benchmark problem is built upon Terzaghi's 1-dimensional consolidation problem \cite{terzaghi1951theoretical}. We assume the domain is homogeneous, isotropic, and saturated with a single-phase fluid. The boundary conditions are described as follows

\begin{equation} \label{eq:ex1_bound_u}
\begin{split}
\bm{u}_{D} \cdot  \mathbf{n}= 0 \quad \si{m} &\text { \: on \: } \mathrm{Left} \times \mathbb{T},\\
\bm{t_D} =[0, -1]  \quad  \si{kPa} &\text { \: on \: } \mathrm{Top} \times \mathbb{T}, \\
\bm{u}_{D} \cdot  \mathbf{n}= 0 \quad  \si{m} &\text { \: on \: } \mathrm{Right} \times \mathbb{T},\\
\bm{u}_{D} \cdot  \mathbf{n}= 0 \quad \si{m} &\text { \: on \: } \mathrm{Bottom} \times \mathbb{T},
\end{split}
\end{equation}

\noindent
for Eq. \eqref{eq:linear_balance}, and

\begin{equation} \label{eq:ex1_bound_p}
\begin{split}
{q}_{D} = 0 \quad \si{m/s} &\text { \: on \: } \mathrm{Left} \times \mathbb{T},\\
{p_D} = 0  \quad  \si{Pa} &\text { \: on \: } \mathrm{Top} \times \mathbb{T}, \\
{q}_{D} = 0 \quad \si{m/s} &\text { \: on \: } \mathrm{Right} \times \mathbb{T},\\
{q}_{D} = 0 \quad \si{m/s} &\text { \: on \: } \mathrm{Bottom} \times \mathbb{T},
\end{split}
\end{equation}

\noindent
for Eq. \eqref{eq:mass_balance}.
A graphical summary of such boundary conditions is reported in Figure \ref{fig:validation_geo}.
The coefficients appearing in section \ref{sec:governing_equations} will either be considered as input parameters, or given fixed values. In particular, we fix $\alpha \approx 1$, as the porous matrix is characterized by $K=1000$ $\mathrm{kPa}$ while the bulk solid is modeled by $K_{\mathrm{s}} \to \infty$ $\mathrm{kPa}$, $c_f = 1.0 \times 10^{-9}$ $\mathrm{Pa^{-1}}$, $\phi = 0.3$, and fluid viscosity - $\mu_f=10^{-3}$ $\mathrm{Pa.s}$. The Poisson ratio $\nu$ and the permeability coefficient $k_{xx}$ are instead considered as input parameters $\bm{\mu} = (\nu, k_{xx})$. The admissible range of variation for $\nu$ is $[0.1, 0.4]$, while that for $k_{xx}$ is $[1.0 \times 10^{-15}, 1.0 \times 10^{-11}]$. For any parametric realization of $\nu$ one can then easily compute the corresponding Lam\'e constants $\lambda_l$ and $\mu_l$ by Eq. \eqref{eq:lambda_l}; for any parametric realization of $k_{xx}$ the matrix permeability tensor $\bm{k}$ is defined as

\begin{equation}  \label{eq:permeability_matrix_iso}
\bm{k}:=\left[ \begin{array}{ll} k_{xx} & k_{xy}  := 0.0 \\ k_{yx} := 0.0 & k_{yy} := k_{xx} \end{array}\right].
\end{equation}

For the validation we aim to carry out in this section we focus on a fixed realization of the uncertain parameter $\bm{\mu}$, and which is outside of the training set. Further discussion on the sensitivity of the reduced order model over the entire parametric range $\mathbb{P}$ will follow in section \ref{sec:sensitivity}. We use a mean squared error ($\mathrm{MSE}_{\varphi}(t, \bm{\mu})$) and maximum error ($\mathrm{ME}_{\varphi}(t, \bm{\mu})$) as the metrics to evaluate our developed framework. $\mathrm{MSE}_{\varphi}(t, \bm{\mu})$ and $\mathrm{ME}_{\varphi}(t, \bm{\mu})$ are defined as follows

\begin{equation}\label{eq:validation_mse}
    {\mathrm{MSE}_\varphi(t, \bm{\mu}) := \left\|\varphi_h(\cdot; t, \bm{\mu}) - \widehat{\varphi}_h(\cdot; t, \bm{\mu})\right\|_{\varphi}^2},
\end{equation}

\noindent
and

\begin{equation}
    {\mathrm{ME}_\varphi(t, \bm{\mu}) := \left\|\varphi_h(\cdot; t, \bm{\mu}) - \widehat{\varphi}_h(\cdot; t, \bm{\mu})\right\|_{\varphi}^{\infty}.}
\end{equation}

\noindent
where $\varphi$ stands for either the displacement $\bm{u}$ or the pressure $p$, $\varphi_h(\cdot; t, \bm{\mu})$ is the corresponding finite element solution at time $t$, and $\widehat{\varphi}_h(\cdot; t, \bm{\mu})$ is the corresponding ROM solution at time $t$. Furthermore, $\left\| \cdot \right\|_{\varphi}^2$ denotes the norm in the space of the primary variable $\varphi$, $\left\| \cdot \right\|_{\varphi}^{\infty}$ the infinity norm (i.e., the maximum pointwise absolute value of its function argument). We remark that, even though Eqs. \eqref{eq:loss_mse} and \eqref{eq:validation_mse} are both MSE errors; they play two fundamentally distinct roles: Eq. \eqref{eq:loss_mse} is employed during the training phase and provides a measurement of the error over the entire time and parametric range; in contrast, Eq. \eqref{eq:validation_mse} is employed during the testing phase and provides a measurement of the error for each time step and each parametric instance. Furthermore, Eq. \eqref{eq:loss_mse} accounts for errors introduced due to the ANN approximation and measures such errors in the space $\mathbb{R}^N$ of the reduced order coefficients; in contrast, Eq. \eqref{eq:validation_mse} accounts for errors introduced by both the POD basis truncation and ANN approximation, and measures such errors in the spatial norm associated to each primary variable. Consistently with this observation, it will then be quite natural to study further whether the error is primarily introduced by the basis truncation or the ANN evaluation, as we will do in section \ref{sec:source_err}.

\begin{figure}[!ht]
   \centering
        \includegraphics[keepaspectratio, height=6.0cm]{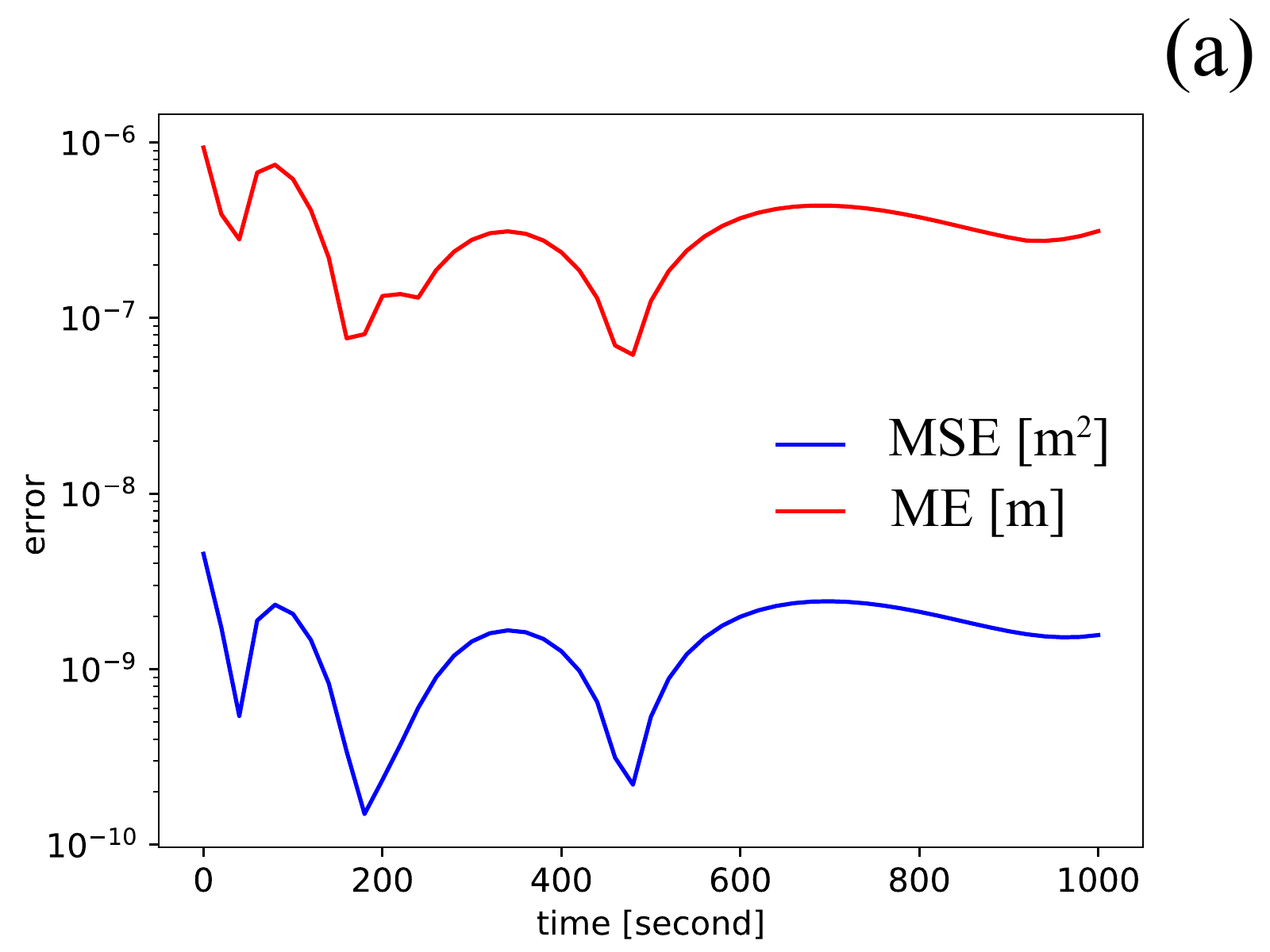}
         \includegraphics[keepaspectratio, height=6.0cm]{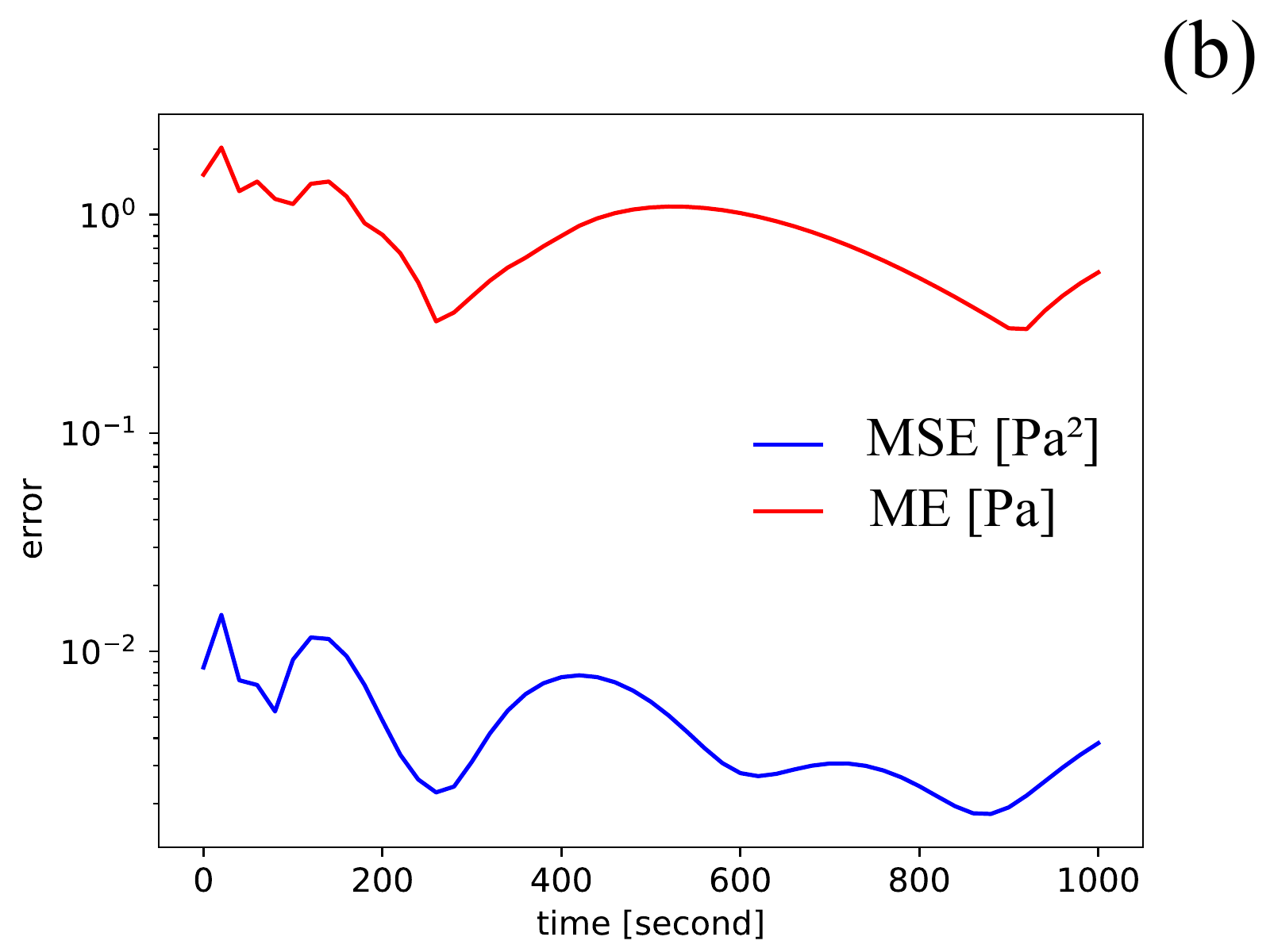}
   \caption{Example 1: mean squared error (MSE) and maximum error (ME) plots using $\bm{\mu} = (\nu, k_{xx}) = (0.2, 1.0 \times 10^{-12})$ - outside of the training snapshots: (a) displacement field ($\bm{u}$) and (b) fluid pressure field ($p$).}
   \label{fig:ex1_err}
\end{figure}

The $\mathrm{MSE}$ and $\mathrm{ME}$ results are presented in Figure \ref{fig:ex1_err}.
We observe that the error in both primary variables maintains the same order of magnitude in the entire time interval. In particular, the reduced order approximation of the displacement field is affected by an MSE of $O(10^{-7})$ \si{m^2} and ME of $O(10^{-9})$ \si{m}. Since the finite element displacement $\bm{u}_h$ is $O( 10^{-4})$ \si{m} for the current parameter value, the corresponding relative ME (i.e., the ratio between ME and the magnitude of the finite element displacement) is $O( 10^{-5})$, which makes the online evaluation an accurate surrogate for any practical engineering scenario. Similarly, the pressure field has MSE of $O(10^{0})$ \si{Pa^2}, ME $O(10^{-2})$ \si{Pa}; since the finite element pressure $p_h$ have values of $O(10^{3})$ \si{Pa}, also the online evaluation of the pressure results in an approximation with a relative ME of $O( 10^{-5})$.


\subsubsection{Example 2: Consolidation problem with anisotropic permeability} \label{sec:2d_ani}

\begin{figure}[!ht]
   \centering
        \includegraphics[keepaspectratio, height=6.0cm]{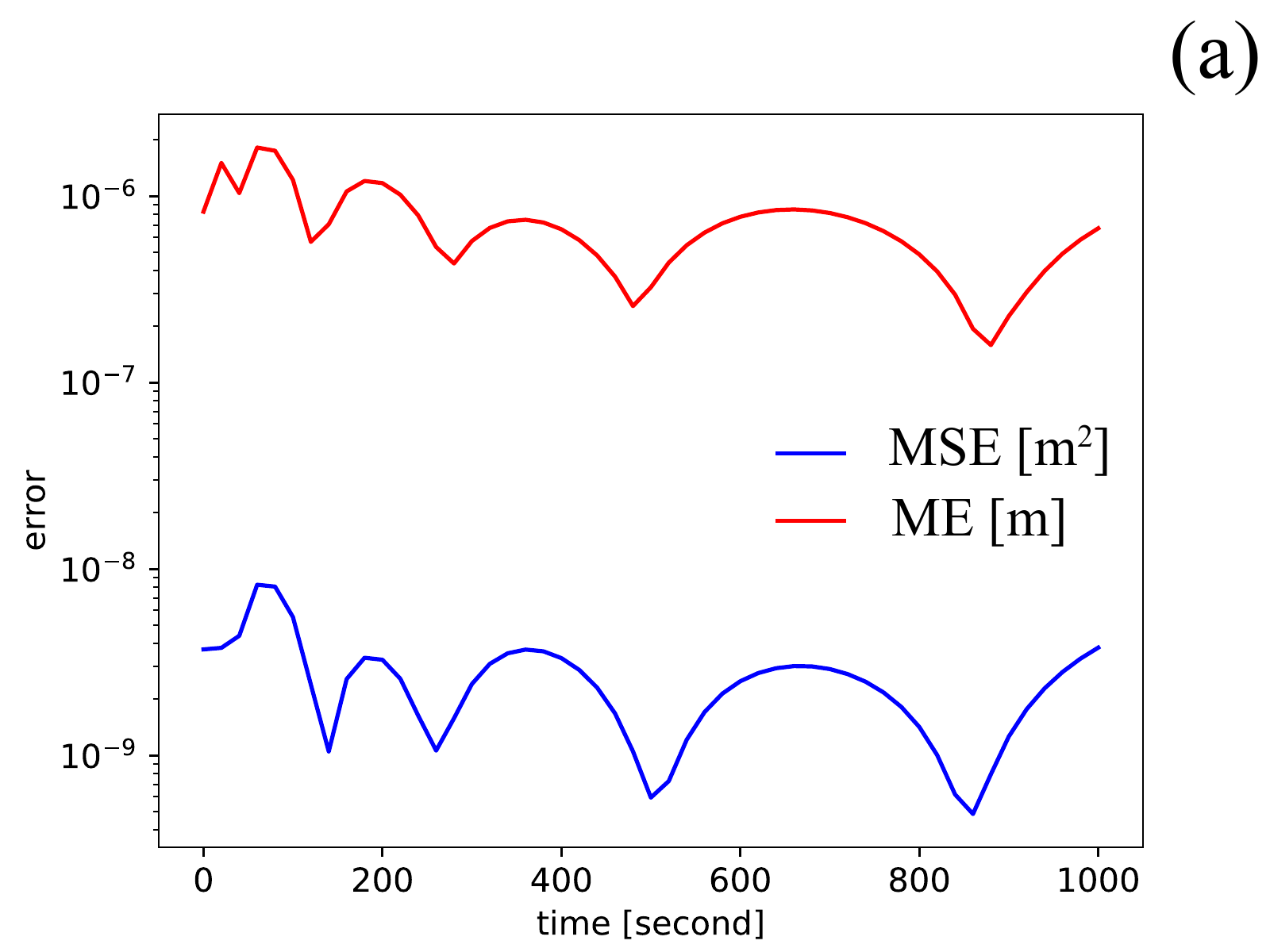}
         \includegraphics[keepaspectratio, height=6.0cm]{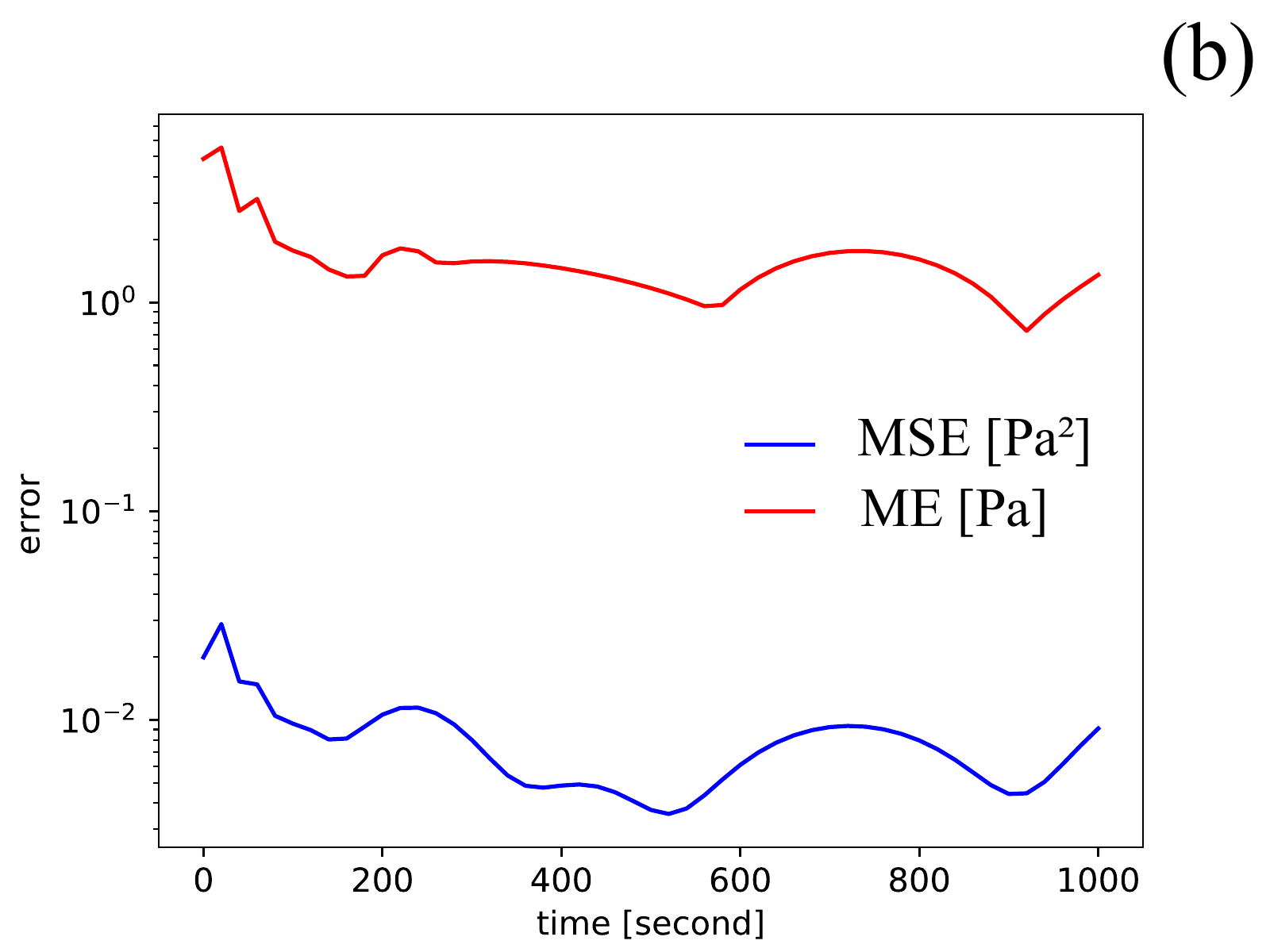}
   \caption{Example 2: mean squared error (MSE) and maximum error (ME) plots using ${\bm{\mu} = (\nu, k_{xx}) = (0.2, 1.0 \times 10^{-12})}$ - outside of the training snapshots: (a) displacement field ($\bm{u}$) and (b) fluid pressure field ($p$).}
   \label{fig:ex2_err}
\end{figure}

We then move to a case where medium permeability is anisotropic.
This benchmark case has been employed, e.g., in \cite{lipnikov2009local,choo2018cracking} and shows the advantages of the DG formulation that we use in this work over traditional finite volume methods, which use a standard two-point flux approximation scheme, as the latter requires the grid to be aligned with the principal directions of the permeability/diffusivity tensors. Boundary conditions, fixed coefficients, and input parameters are as in Section \ref{sec:1d_1l}, except for the anisotropic permeability tensor

\begin{equation}  \label{eq:permeability_matrix_ani}
\bm{k}:=\left[ \begin{array}{ll} k_{xx} & k_{xy} := 0.1 k_{xx} \\k_{yx} := 0.1 k_{xx} & k_{yy} := 0.1 k_{xx} \end{array}\right],
\end{equation}

\noindent
where $k_{xx}$ is the second input parameter. The $\mathrm{MSE}$ and $\mathrm{ME}$ results of this case are illustrated in Figure \ref{fig:ex2_err}. Similarly to the previous example, considering that $\bm{u}_h$ and $p_h$ at the initial state are $O( 10^{-4})$ \si{m} and $O( 10^{3})$ \si{Pa}, respectively, we get a relative ME of $O(10^{-5})$ for both displacement and pressure.


\subsubsection{Example 3: Consolidation problem with 2-layered material} \label{sec:1d_2l}

\begin{figure}[!ht]
   \centering
        \includegraphics[width=6.0cm,keepaspectratio]{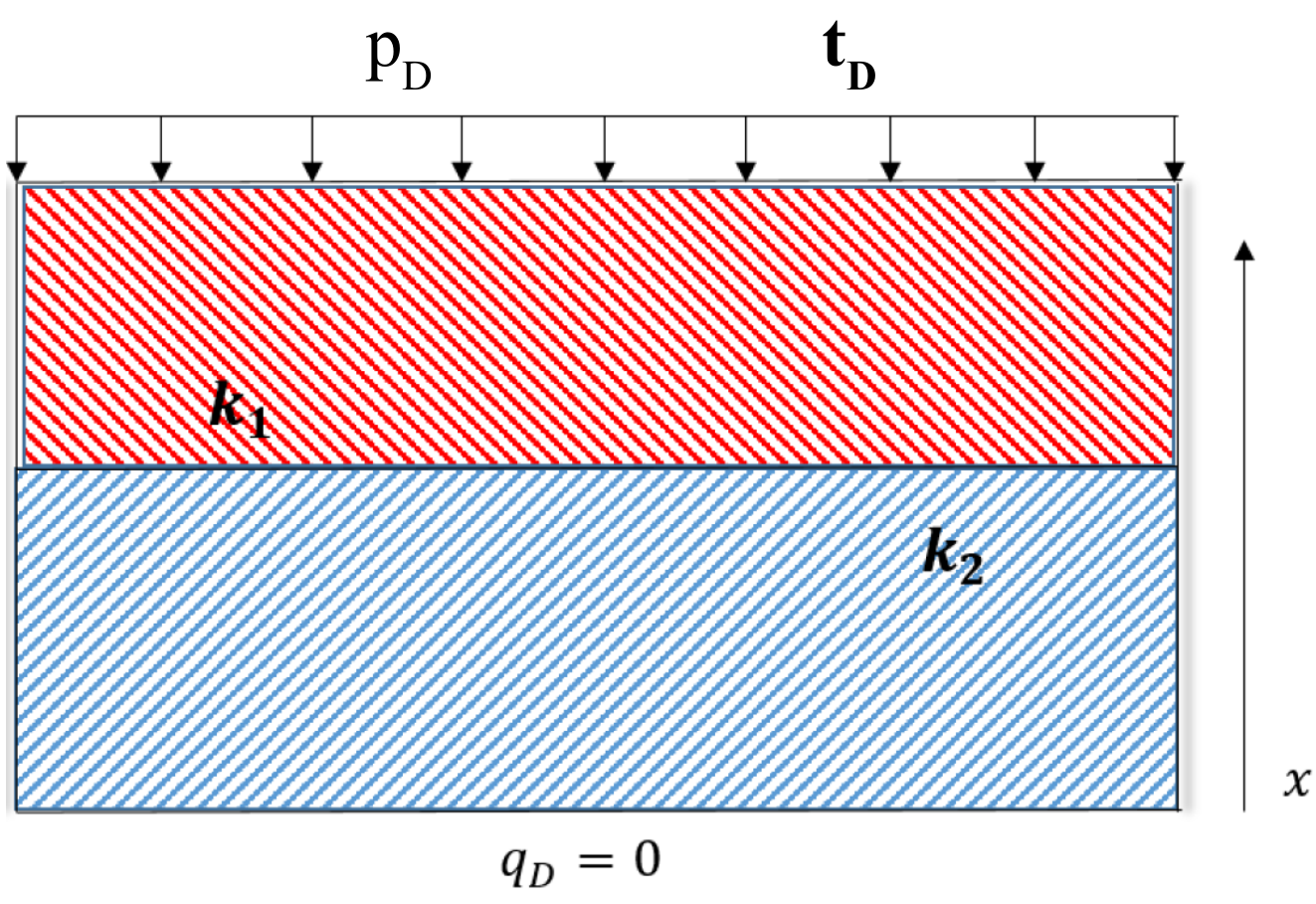}
   \caption{Example 3: Setup for the 1-dimensional consolidation problem in a 2-layered material. This figure is adapted from \cite{kadeethum2020finite}.}
   \label{fig:mat_two_case_geo}
\end{figure}

Finally, we evaluate the developed model reduction framework using a 2-layered material as presented in Figure \ref{fig:mat_two_case_geo}. Boundary conditions, fixed coefficients, and input parameters are as in Section \ref{sec:1d_1l}, except for medium permeability defined as

\begin{equation}  \label{eq:permeability_matrix_iso_1}
\bm{k}(x, y) = \begin{cases}
\bm{k}_1, y > 0.5\\
\bm{k}_2, y < 0.5\\
\end{cases},
\quad \text{where} \quad
\begin{aligned}
\bm{k}_1&:=\left[   \begin{array}{ll} 1.0 \times 10^{-12}  &  0.0 \\ 0.0 & 1.0 \times 10^{-12} \end{array}\right], \text{and}\\
\bm{k}_2&:=\left[ \begin{array}{ll} k_{xx} & 0.0 \\ 0.0 & k_{xx} \end{array}\right].
\end{aligned}
\end{equation}

\noindent
The second input parameter thus affects the isotropic permeability $\bm{k}_2$ in the bottom subdomain depicted in Figure \ref{fig:mat_two_case_geo}; instead, the permeability tensor $\bm{k}_1$ is parameter independent in the top subdomain. We restrict the range for the second parameter $k_{xx}$ to the interval $[1.0 \times 10^{-16}, 1.0 \times 10^{-15}]$ to simulate parametric configurations in which the two layers have very different material properties.

\begin{figure}[!ht]
   \centering
        \includegraphics[keepaspectratio, height=6.0cm]{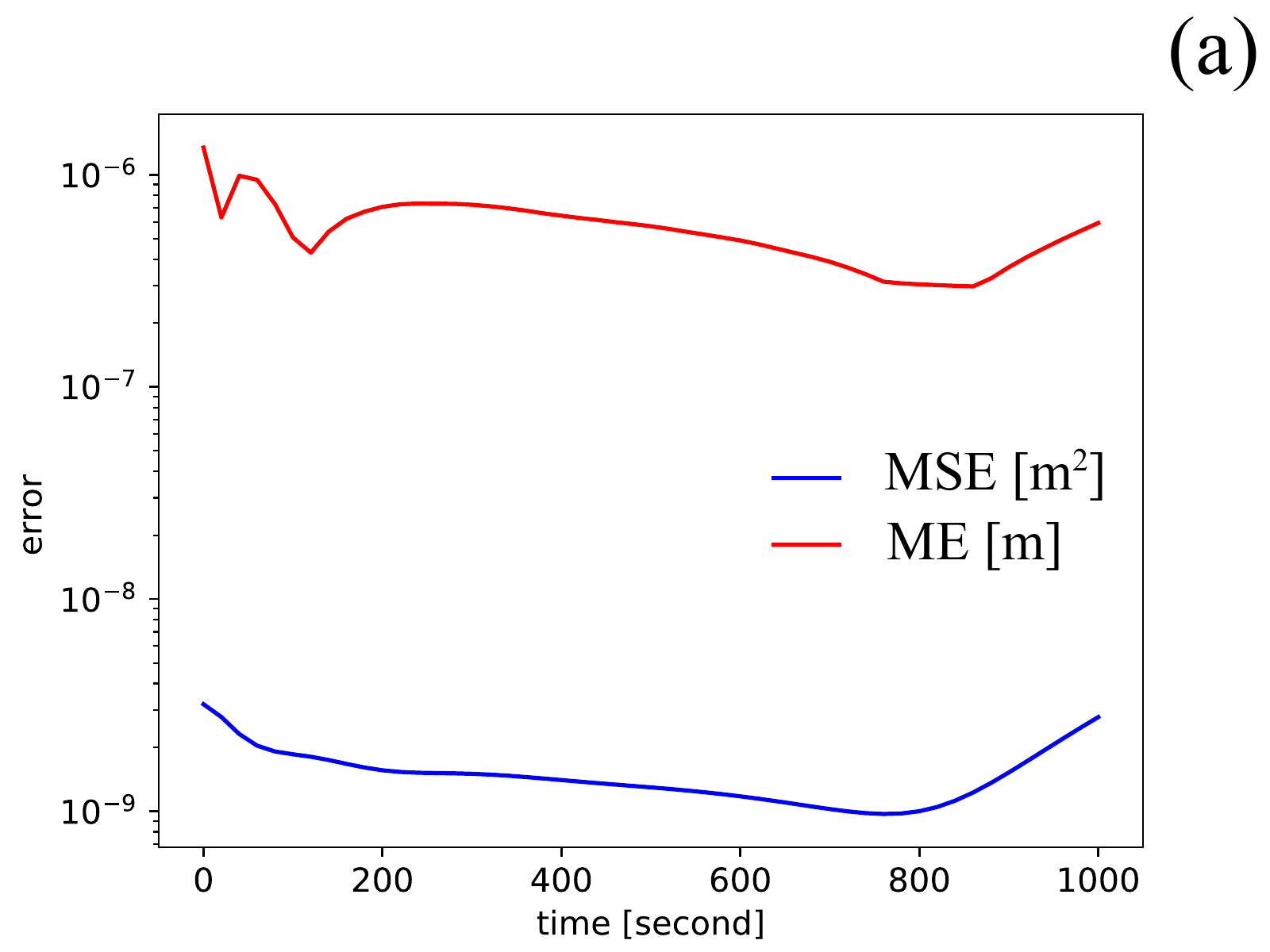}
         \includegraphics[keepaspectratio, height=6.0cm]{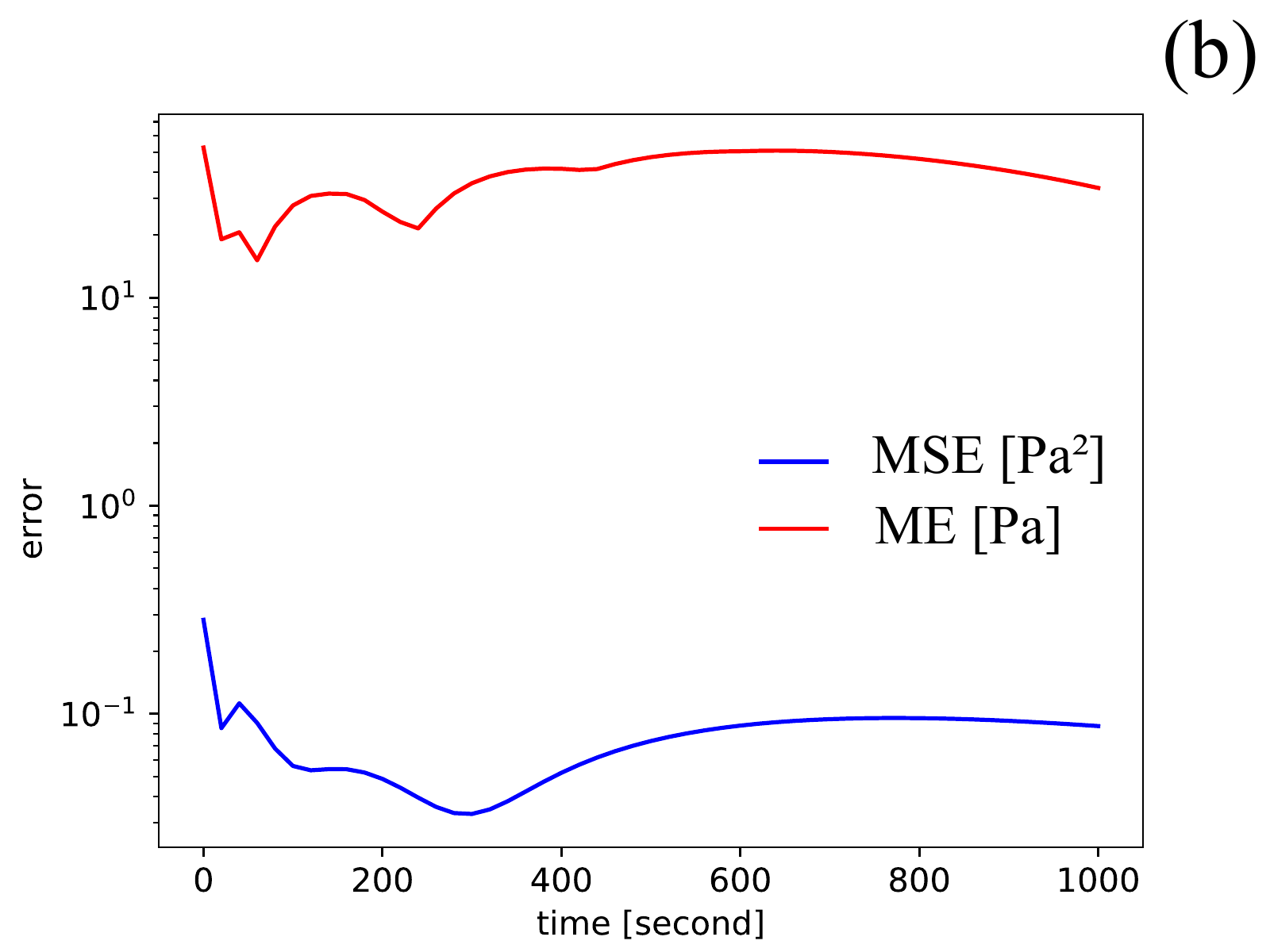}
          \includegraphics[keepaspectratio, height=6.0cm]{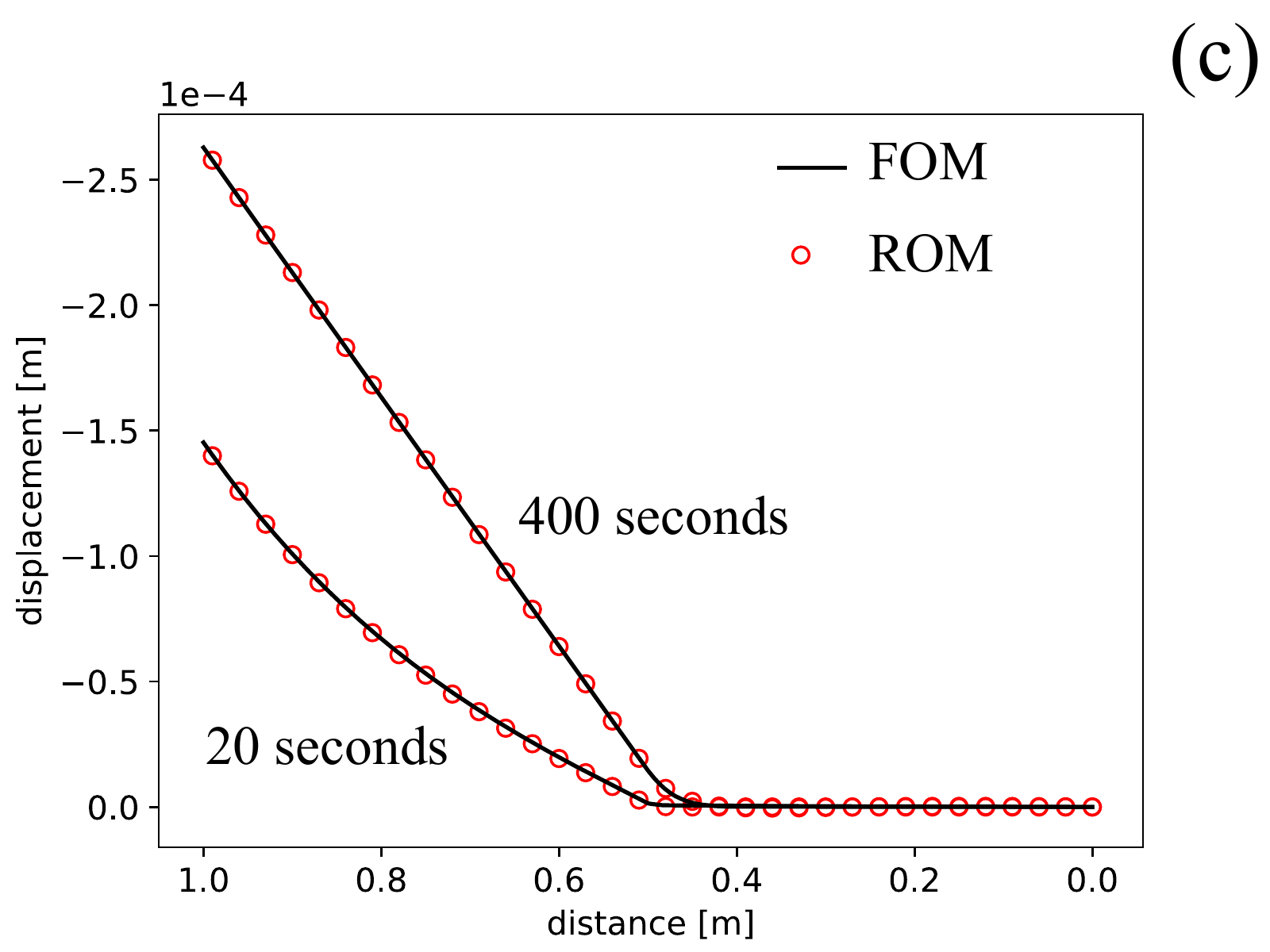}
         \includegraphics[keepaspectratio, height=6.0cm]{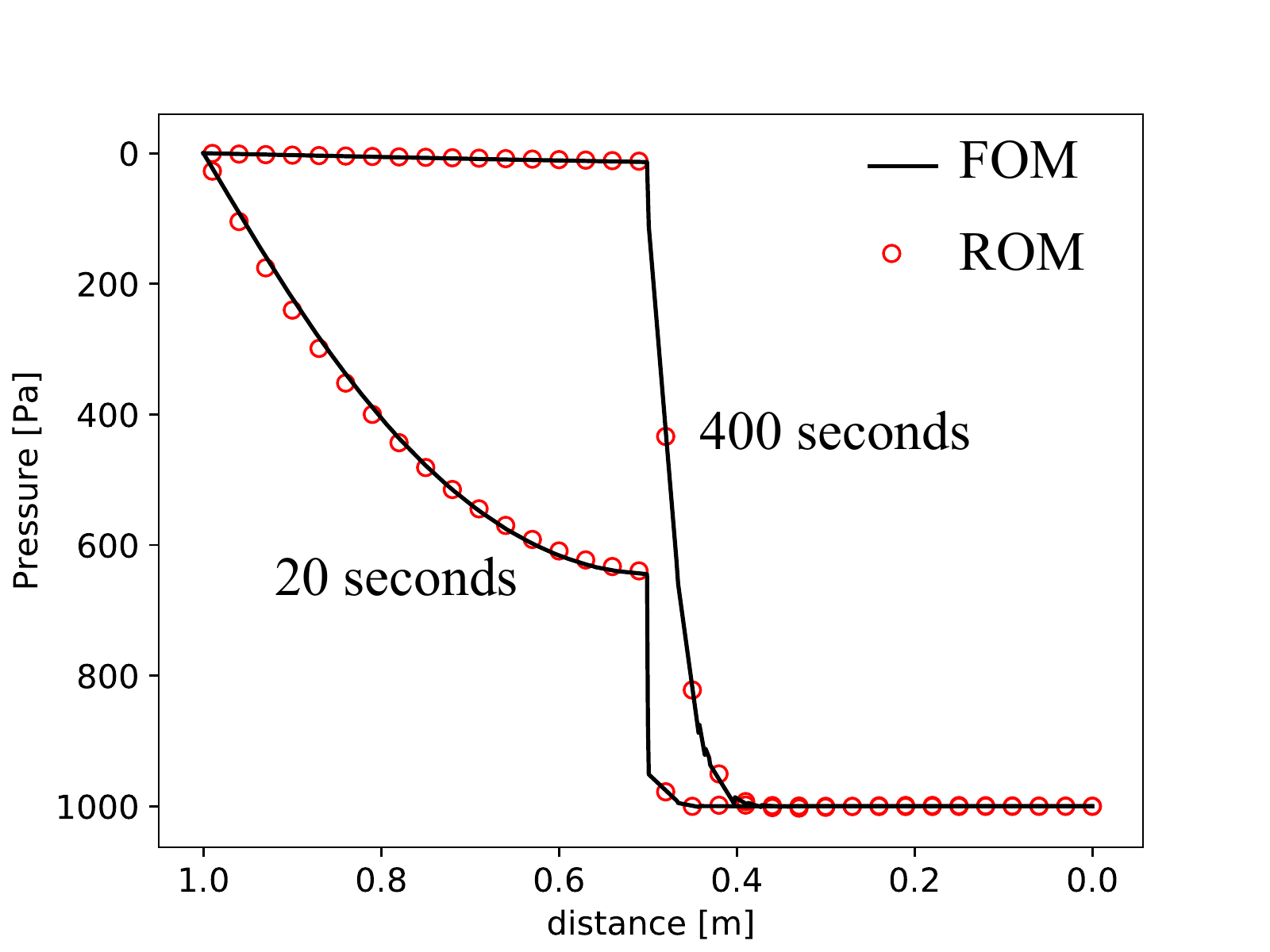}
   \caption{Example 3: mean squared error (MSE) and maximum error (ME) plots using ${\bm{\mu} = (\nu, k_{xx}) = (0.2, 5.0 \times 10^{-15})}$ - outside of the training snapshots: (a) displacement field ($\bm{u}$), (b) fluid pressure field ($p$), full order model (FOM) and reduced order model (ROM) results along the $x = 0.5$ line of (c) displacement field ($\bm{u}$), and (d) displacement field ($\bm{u}$).}
   \label{fig:ex3_err}
\end{figure}

\noindent
This test case has been used, e.g., in \cite{choo2018enriched,Kadeethum2019ARMA,kadeethum2020finite,kadeethum2020enriched} to underlying how, without a suitable stabilization, the solution may exhibit spurious oscillations at the interface between two layers. Since the DG method we employ in this work is oscillation-free, not only the FOM solutions will not have spurious oscillations at the interface (i.e., at $y = 0.5$), but also the ROM fulfill such desirable property, see Figures \ref{fig:ex3_err}c-d.
For what concerns the quantitative behavior of the error in time, we get a relative ME (the ratio between ME and the magnitude of the finite element solutions) of $O(10^{-2})$ for both primary variables. We note that $\mathrm{MSE}$, $\mathrm{ME}$ and relative ME values are higher than the ones presented in Sections \ref{sec:1d_1l} and \ref{sec:2d_ani} since the permeability is strongly heterogeneous, exhibiting a sharp contrast between the two material phases resulting in a discontinuity in the pressure field, which makes the problem significantly more challenging than the two previous test cases.


\subsection{Model analysis} \label{sec:analysis}

Following the verification of the developed ROM framework for a representative realization of the input parameters and for fixed values of the hyperparameters, we now perform a comprehensive analysis of the ROM using a more realistic example.

\subsubsection{Example 4: Consolidation problem with heterogeneous permeability} \label{sec:2d_het}

\begin{figure}[!ht]
   \centering
        \includegraphics[width=6.0cm,keepaspectratio]{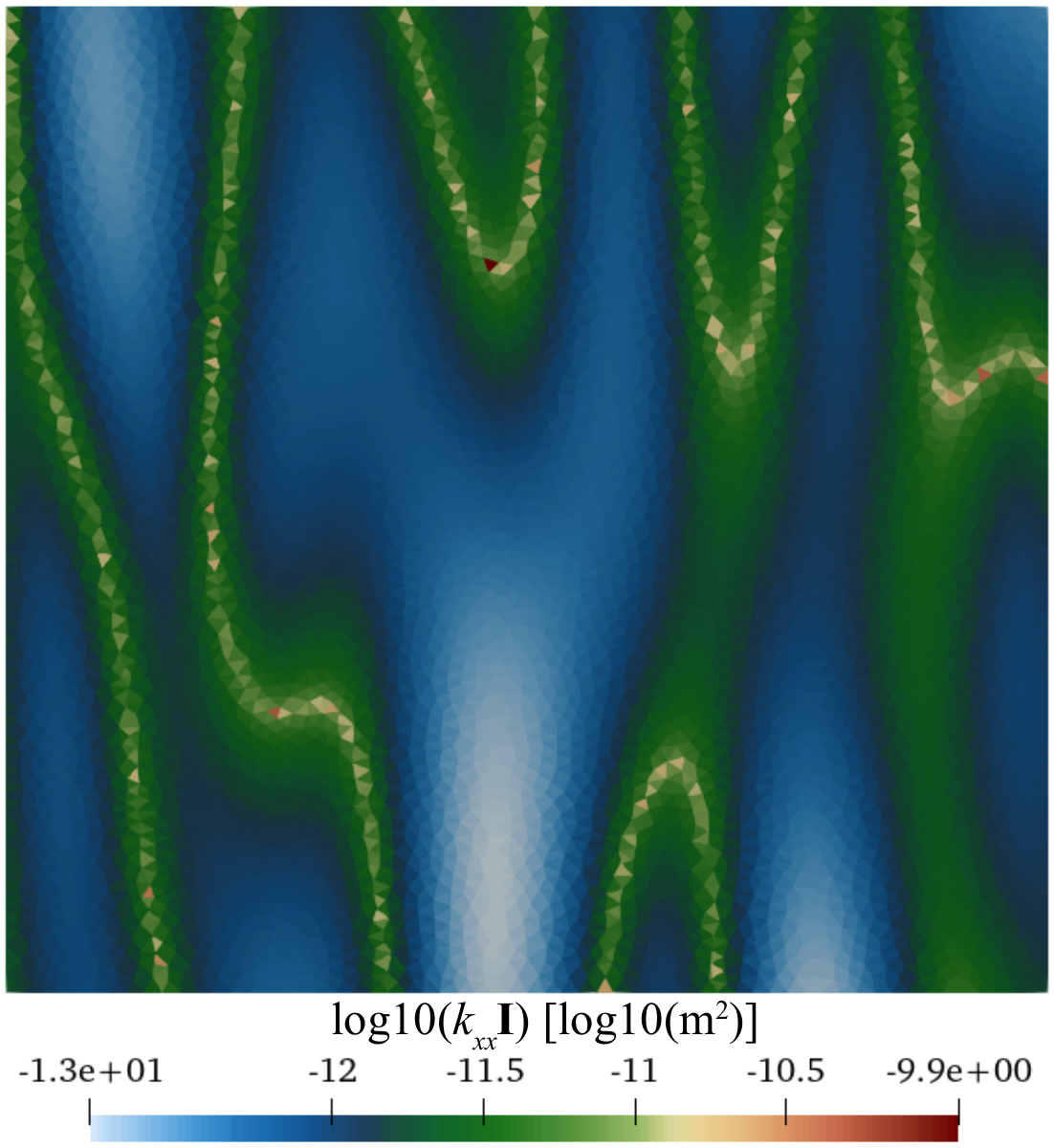}
   \caption{Example 4: A permeability field ($\bm{k_m} = k_{xx} \mathbf{I}$) generated as in \cite{muller2020}. The average of $k_{xx} = 1.77 \times 10^{-12}$ \si{m^2}. The variance of $k_{xx} = 5.53 \times 10^{-24}$ \si{m^4}. The Zinn \& Harvey transformation is applied at then end of the permeability field generation \cite{muller2020,zinn2003good}.}
   \label{fig:het_geo}
\end{figure}

In this example, the matrix permeability is heterogeneous as presented in Figure \ref{fig:het_geo}. This field is generated as in \cite{muller2020} with the average of $k_{xx} = 1.77 \times 10^{-12}$ \si{m^2}, and the variance of $k_{xx} = 5.53 \times 10^{-24}$ \si{m^4}. The field $\bm{k_m}$ is isotropic, which means $\bm{k_m} = k_{xx} \mathbf{I}$. The Zinn \& Harvey transformation is applied at then end of the permeability field generation \cite{muller2020,zinn2003good}.
Thus, in contrast to the previous test cases, the matrix permeability is fixed and not parametric. The input parameters are $\bm{\mu} = (\nu, \alpha) \in [0.1, 0.4] \times [0.4, 1.0]$. The remaining coefficients, as well as the boundary conditions, are as in Section \ref{sec:1d_1l}.
We note that the magnitude of the initial values for $\bm{u}_h$ and $p_h$ are $O(10^{-4})$ \si{m} and $O(10^{3})$ \si{Pa}, respectively throughout this example.

\begin{figure}[!ht]
   \centering
        \includegraphics[keepaspectratio, height=6.0cm]{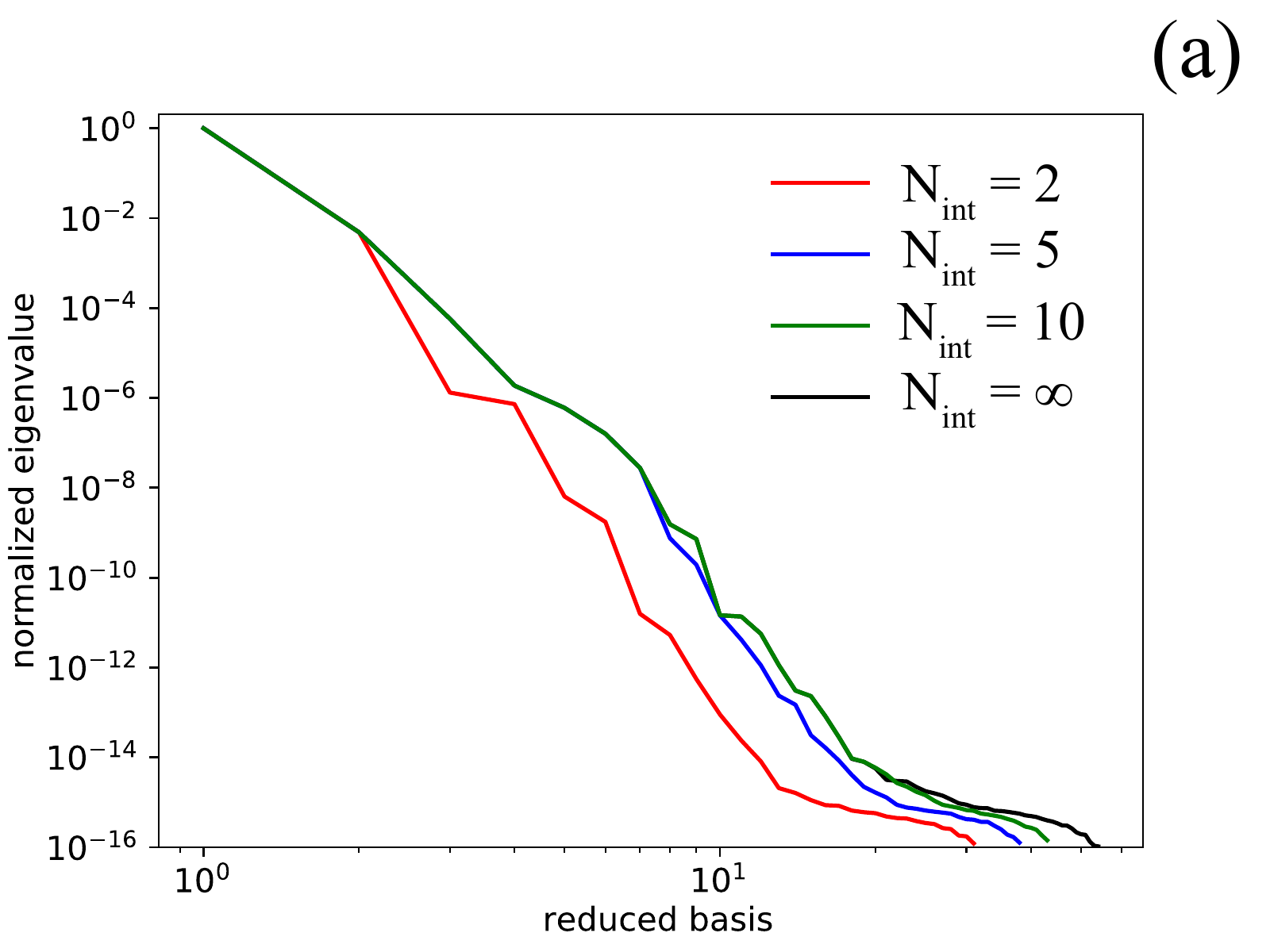}
         \includegraphics[keepaspectratio, height=6.0cm]{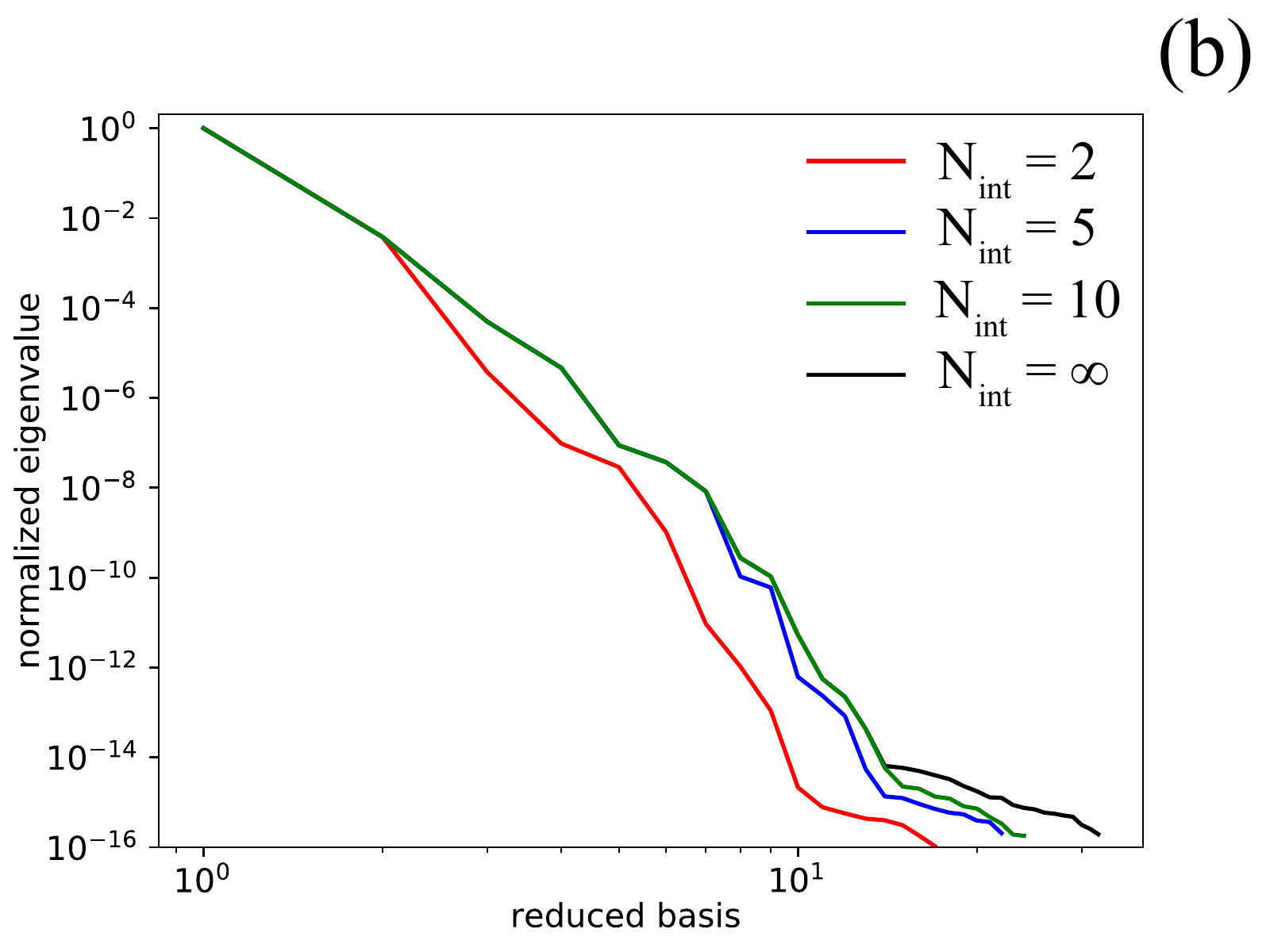}
   \caption{Example 4: normalized eigenvalue as a function of basis for (a) displacement field ($\bm{u}$) (b) fluid pressure field ($p$). $\mathrm{N_{int}} = \infty$ represents a case where we do not use the nested POD technique.}
   \label{fig:ex4_rb}
\end{figure}

The eigenvalue behavior obtained from the POD phase for both displacement and pressure fields is presented in Figure \ref{fig:ex4_rb}. We note that these eigenvalues are normalized by their maximum value for the sake of presentation. From this figure, we observe that by using 30 to 50 reduced bases, we could capture most of the information produced by FOM (i.e., the normalized eigenvalue reaches the machine precision.), regardless of the choice of $\mathrm{N_{int}}$. Besides, as the number of $\mathrm{N_{int}}$ increases, the behavior of eigenvalue becomes similar to the standard POD ($\mathrm{N_{int}} = \infty$) case (i.e., as the number of $\mathrm{N_{int}}$ increases, we could capture most of the information in the time domain; hence, there is no difference between the nested POD and standard POD.). By using $\mathrm{N_{int}} = 10$, the eigenvalue behavior is almost identical to the case where we use $\mathrm{N_{int}} = \infty$. In fact, the lines overlap for the most part of the plot, except for the trailing eigenvalues, which are below numerical precision.

The comparison of the wall time (seconds) used for SVD computations with a different number of snapshots ($\mathrm{M}$) and a number of intermediate reduced basis ($\mathrm{N_{int}}$) is presented in Table \ref{tab:time_num}. We observe that the $\mathrm{N_{int}} = \infty$ case consumes the longest wall time. Furthermore, the lower the number of  $\mathrm{N_{int}}$, the faster SVD computations are. In a general sense, we note that the nested POD technique could reduce the wall time required by the SVD computations significantly. For instance, the $\mathrm{N_{int}} = 10$ case provides a comparable eigenvalue behavior to the $\mathrm{N_{int}} = \infty$ case, but the generation of the reduced spaces is approximately ten times faster.

\begin{table}[!ht]
\centering
\caption{Example 4: Comparison of the wall time (seconds) used for proper orthogonal decomposition (POD) operation with different number of snapshots ($\mathrm{M}$) and number of intermediate reduced basis ($\mathrm{N_{int}}$). $\mathrm{N_{int}} = \infty$ represents a case where we do not use the nested POD technique.}
\begin{tabular}{|l|c|c|c|c|}
\hline
   &  \multicolumn{1}{l|}{$\mathrm{N_{int}} = 2$} & \multicolumn{1}{l|}{$\mathrm{N_{int}} = 5$} & \multicolumn{1}{l|}{$\mathrm{N_{int}} = 10$} &
   \multicolumn{1}{l|}{$\mathrm{N_{int}} = \infty$} \\ \hline
$\mathrm{M} = 100$                           & 100                             & 125                             & 170  & 1574                               \\ \hline
$\mathrm{M} = 400$                            & 437                             & 650                             & 1437 & 36705                             \\ \hline
$\mathrm{M} = 900$                          & 1168                             & 2319                             & 6475                        &     268754      \\ \hline
\end{tabular}
\label{tab:time_num}
\end{table}

\subsubsection{Sources of error} \label{sec:source_err}

Throughout this subsection, we study the model accuracy and the sources of error. Our goal is to differentiate the error arising due to the truncation of reduced bases (i.e., associated to a choice $\mathrm{N}$ that is less than $\mathrm{M} N^t$.) and the error introduced by the mapping between $(t, \bm{\mu})$ and the reduced order coefficients $\widehat{\bm{\theta}}^u(t, \bm{\mu})$ and $\widehat{\bm{\theta}}^p(t, \bm{\mu})$ provided by ANN. In particular, we will

\begin{figure}[!ht]
   \centering
        \includegraphics[keepaspectratio, height=6.0cm]{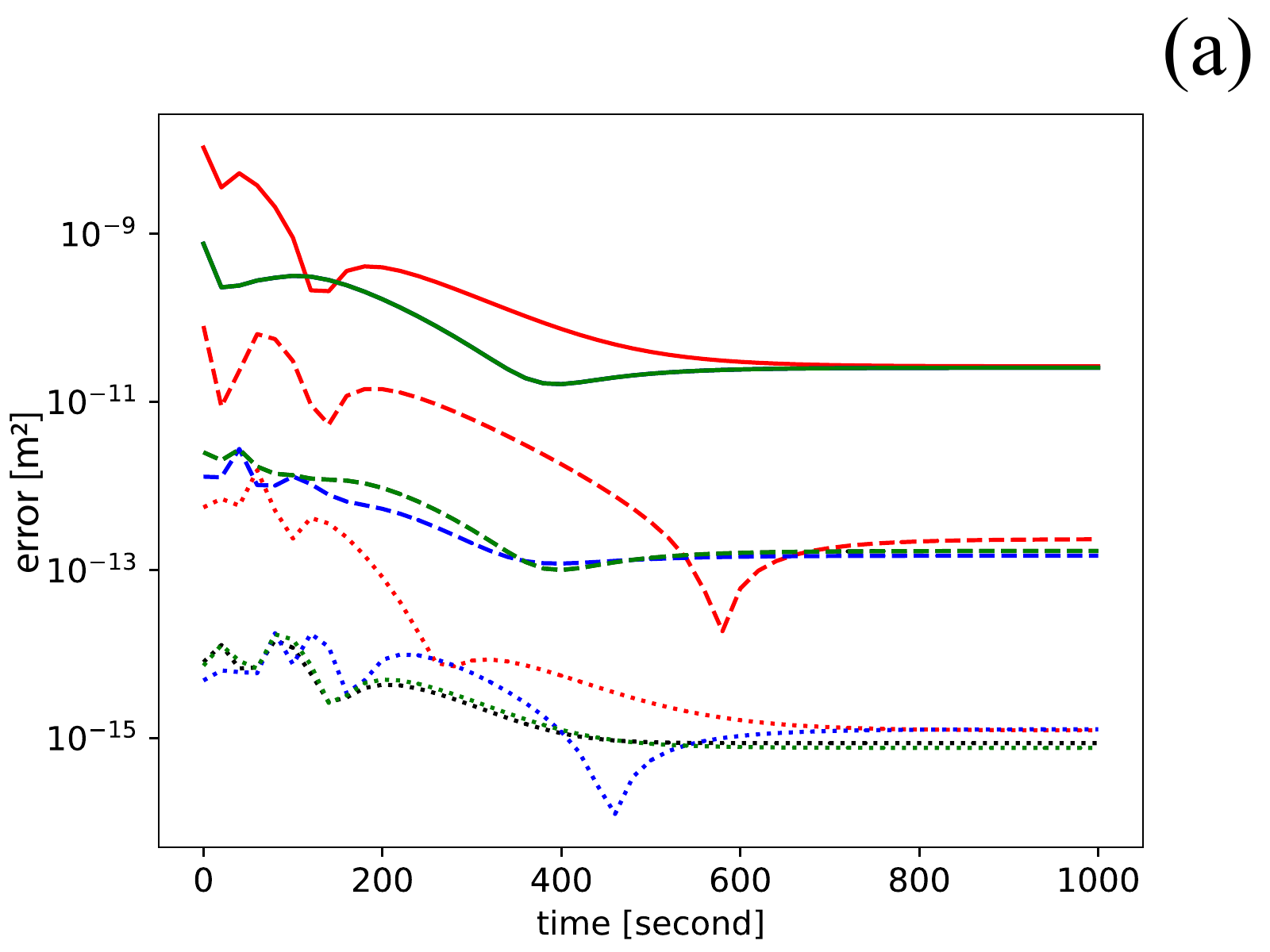}
         \includegraphics[keepaspectratio, height=6.0cm]{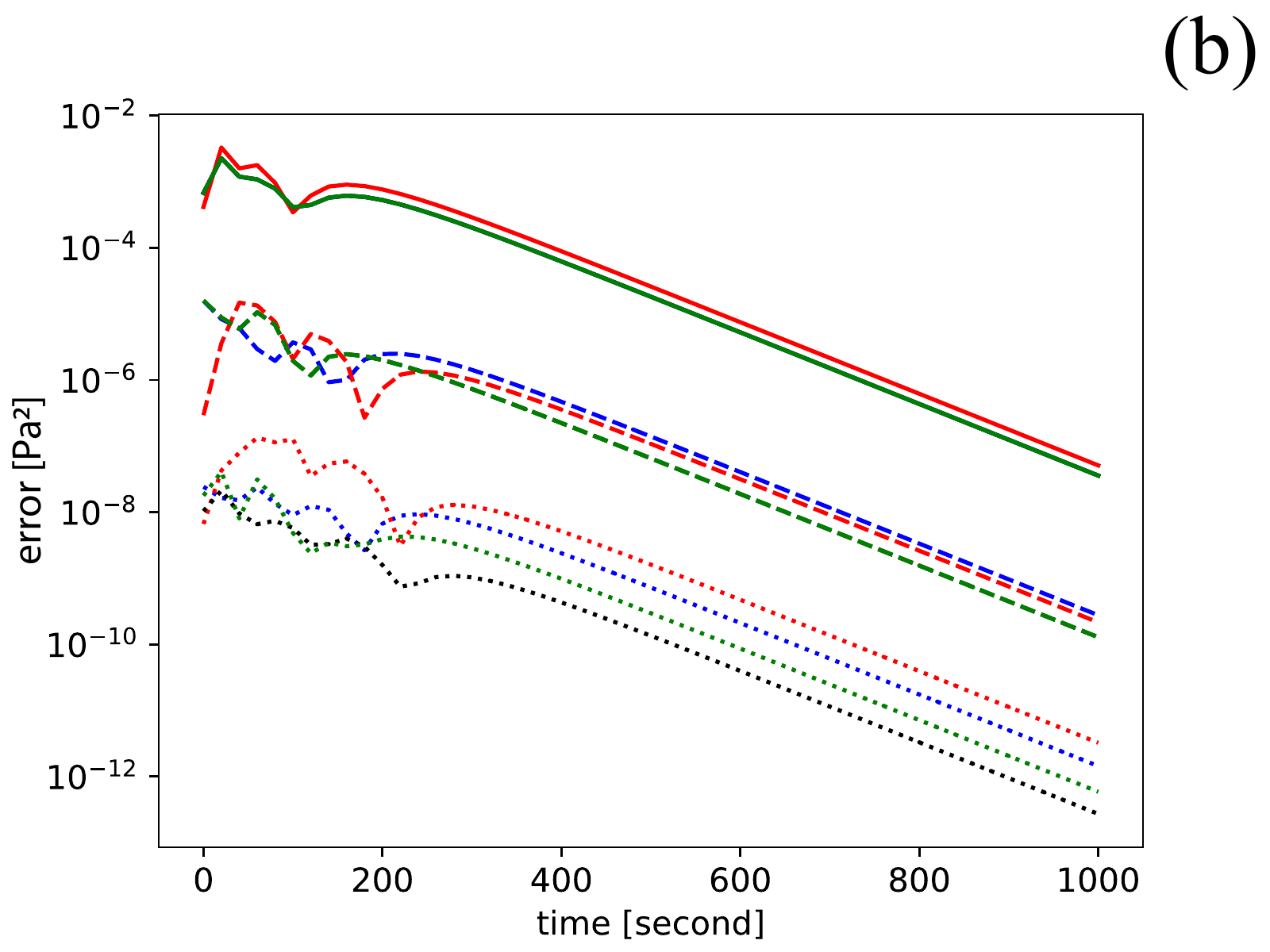}
         \includegraphics[keepaspectratio, height=6.0cm]{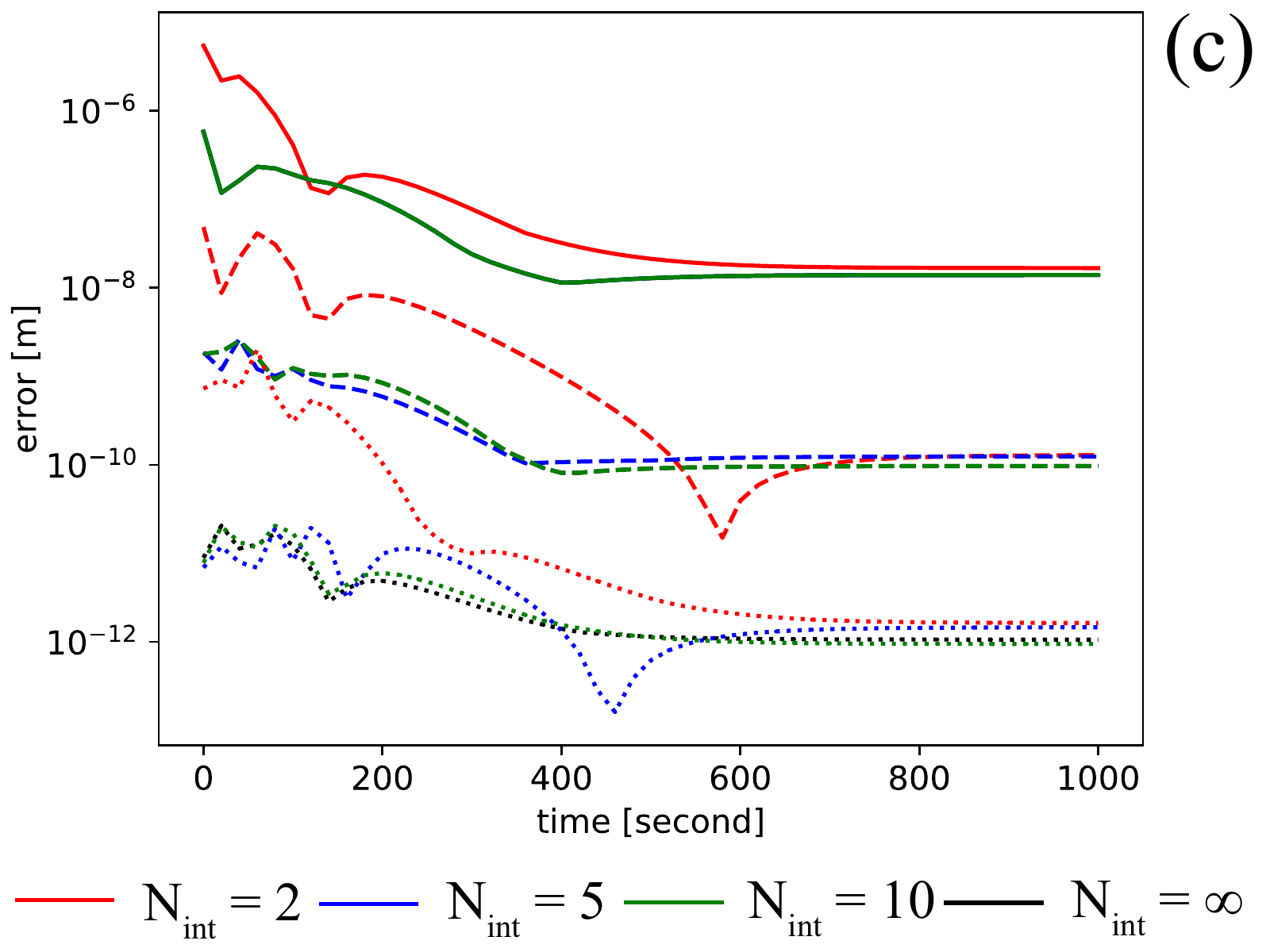}
         \includegraphics[keepaspectratio, height=6.0cm]{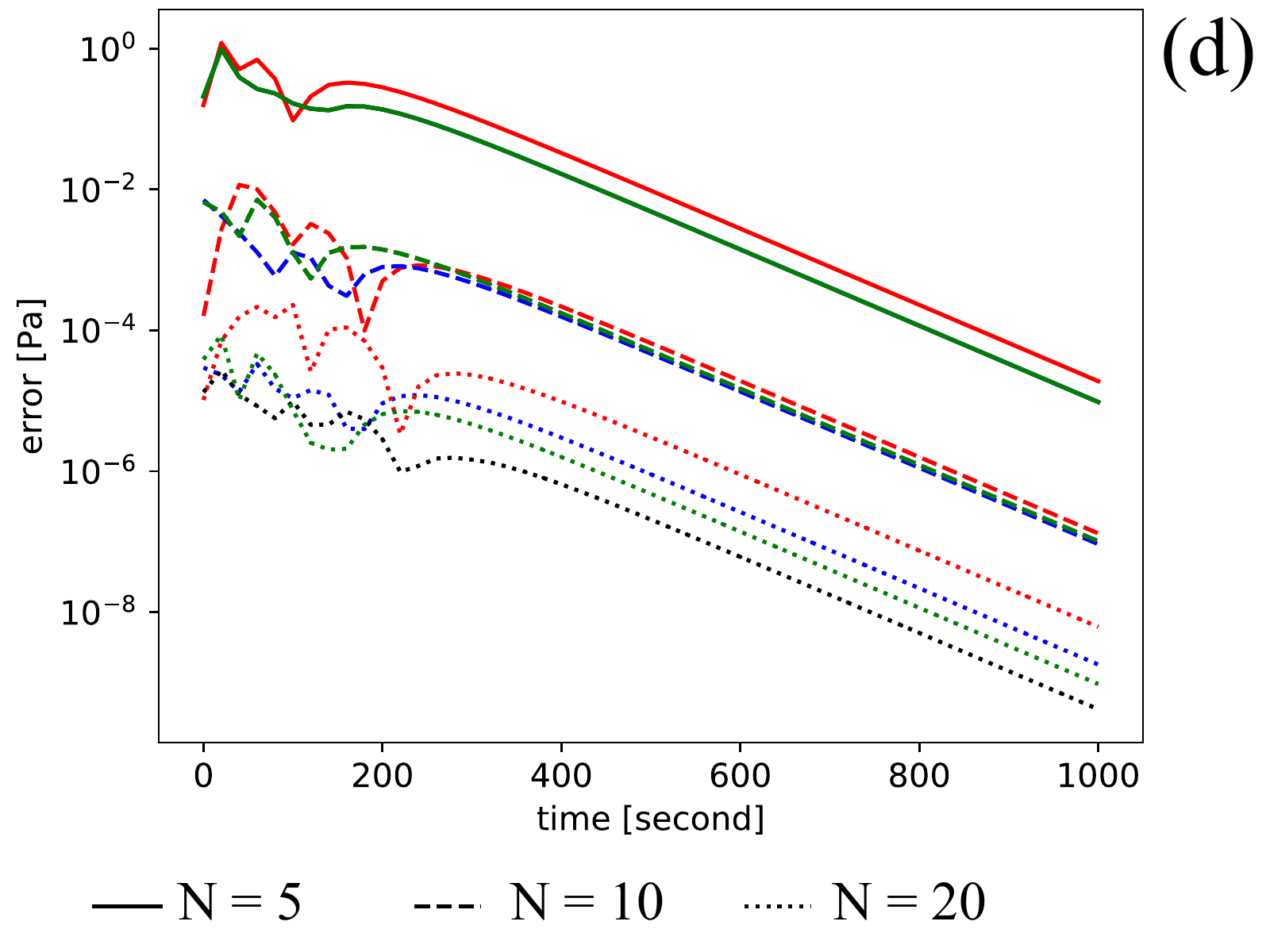}
   \caption{Example 4: Errors of reconstruction solutions using $\bm{\mu} = (\nu, \alpha) = (0.1, 0.8)$ - in the training snapshots and using ${\bm{\theta}}^u$, ${\bm{\theta}}^p$, using different numbers of reduced basis ($\mathrm{N}$): (a) mean squared error (MSE) of displacement field ($\bm{u}$), (b) mean squared error (MSE) of fluid pressure field ($p$), (c) maximum error (ME) of displacement field ($\bm{u}$) (d) maximum error (ME) of fluid pressure field ($p$). Colors correspond to increasing values of $\mathrm{N_{int}}$; solid, dashed, and dotted lines represent $\mathrm{N} = 5$, $\mathrm{N} = 10$, and $\mathrm{N} = 20$ cases,  respectively.}
   \label{fig:ex4_source_err_1}
\end{figure}

\begin{figure}[!ht]
   \centering
        \includegraphics[keepaspectratio, height=6.0cm]{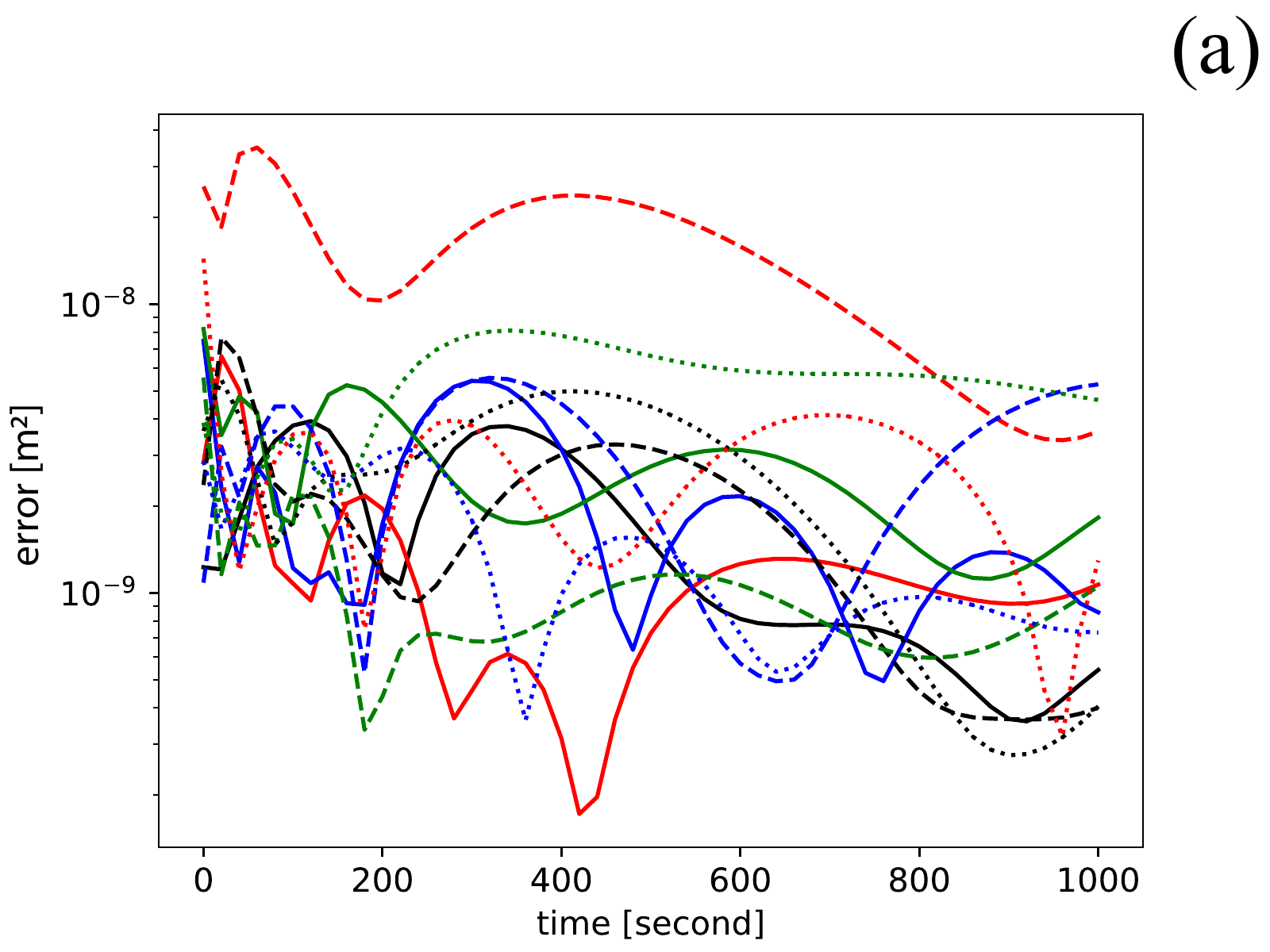}
         \includegraphics[keepaspectratio, height=6.0cm]{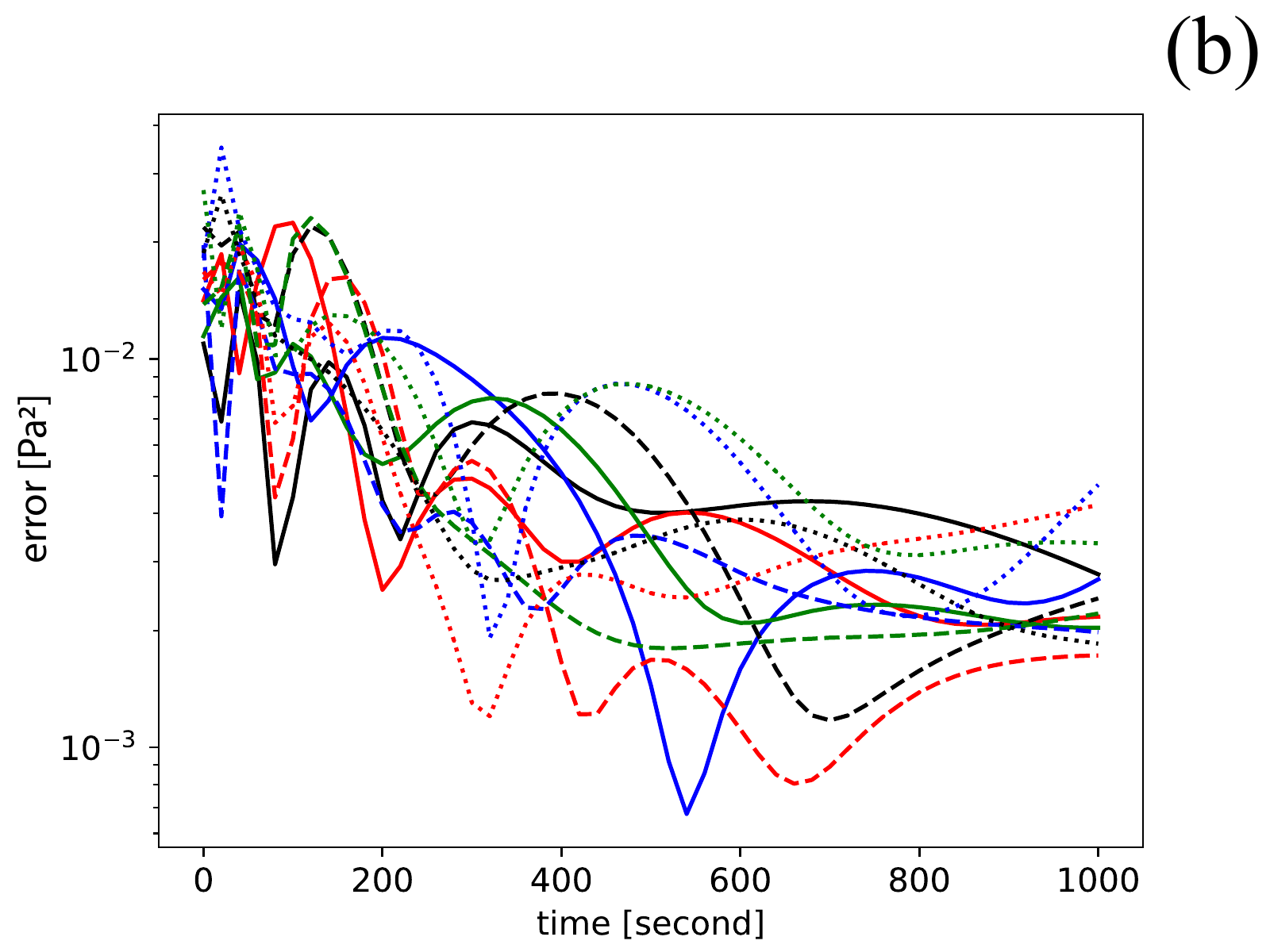}
         \includegraphics[keepaspectratio, height=6.0cm]{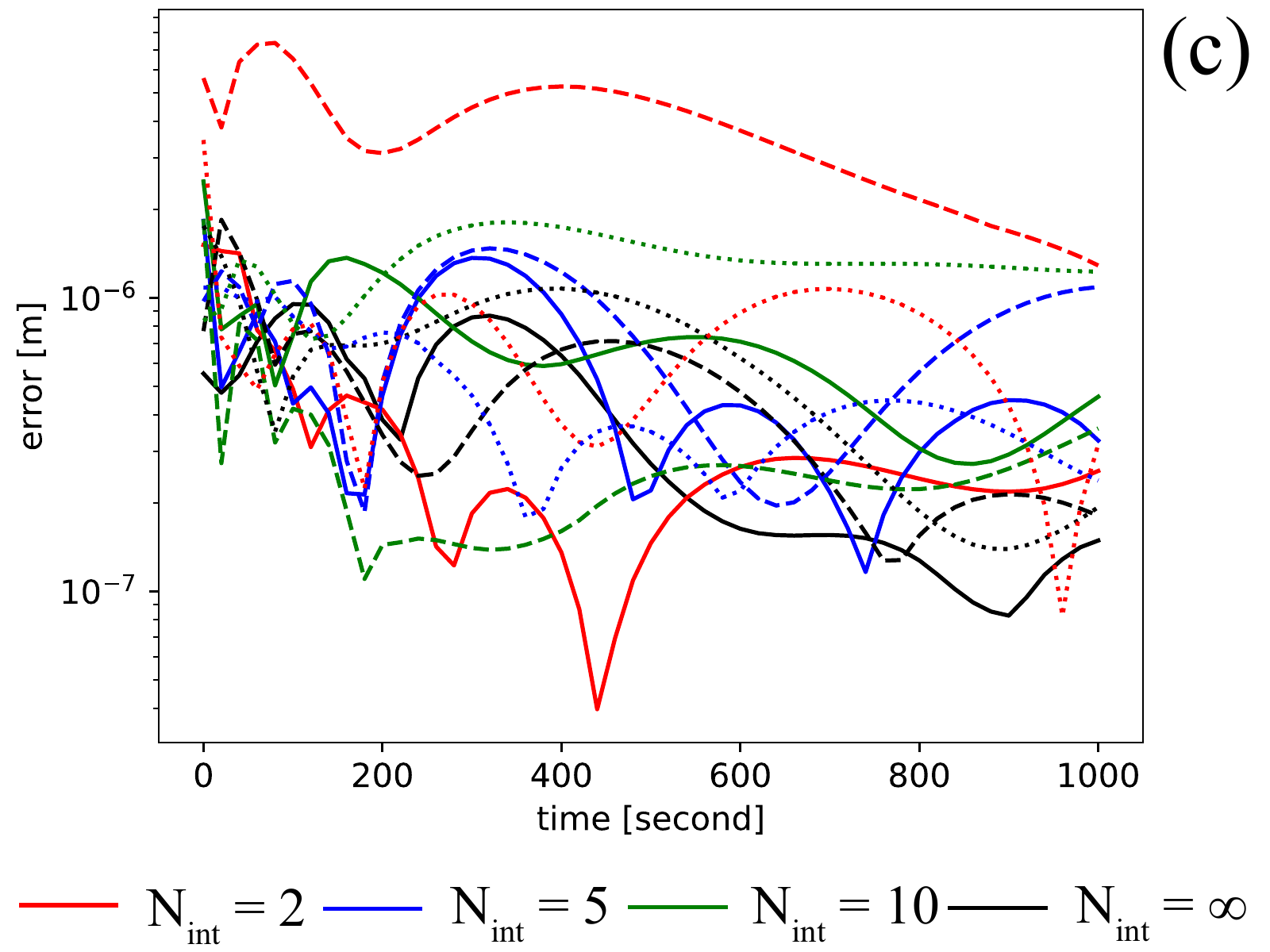}
         \includegraphics[keepaspectratio, height=6.0cm]{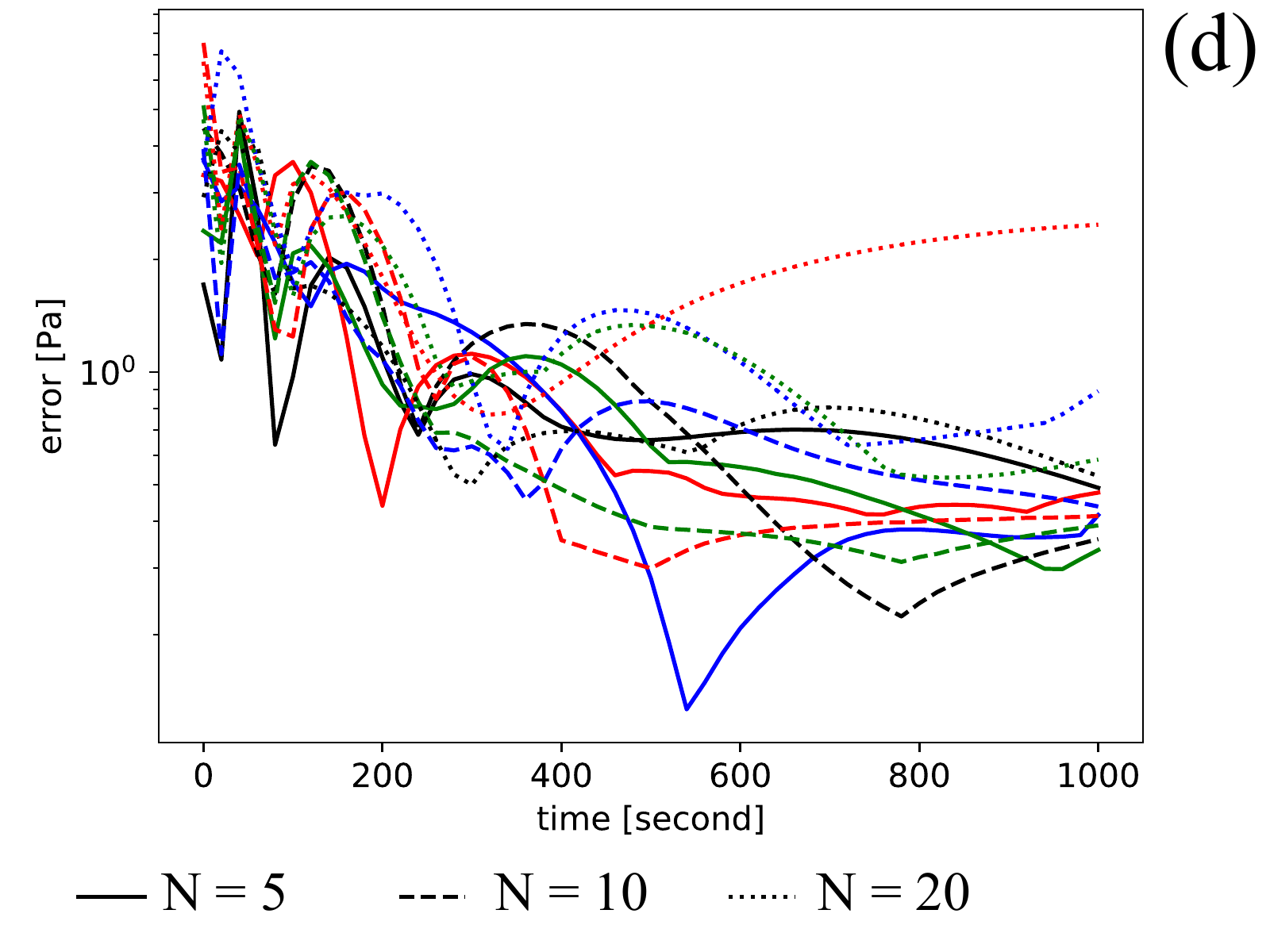}
   \caption{Example 4: Errors of reconstruction solutions using $\bm{\mu} = (\nu, \alpha) = (0.1, 0.8)$ - in the training snapshots and $\widehat{\bm{\theta}}^u$, $\widehat{\bm{\theta}}^p$, using different numbers of reduced basis ($\mathrm{N}$): (a) mean squared error (MSE) of displacement field ($\bm{u}$), (b) mean squared error (MSE) of fluid pressure field ($p$), (c) maximum error (ME) of displacement field ($\bm{u}$) (d) maximum error (ME) of fluid pressure field ($p$). Colors correspond to increasing values of $\mathrm{N_{int}}$; solid, dashed, and dotted lines represent $\mathrm{N} = 5$, $\mathrm{N} = 10$, and $\mathrm{N} = 20$ cases,  respectively.}
   \label{fig:ex4_source_err_2}
\end{figure}

\begin{figure}[!ht]
   \centering
         \includegraphics[keepaspectratio, height=6.0cm]{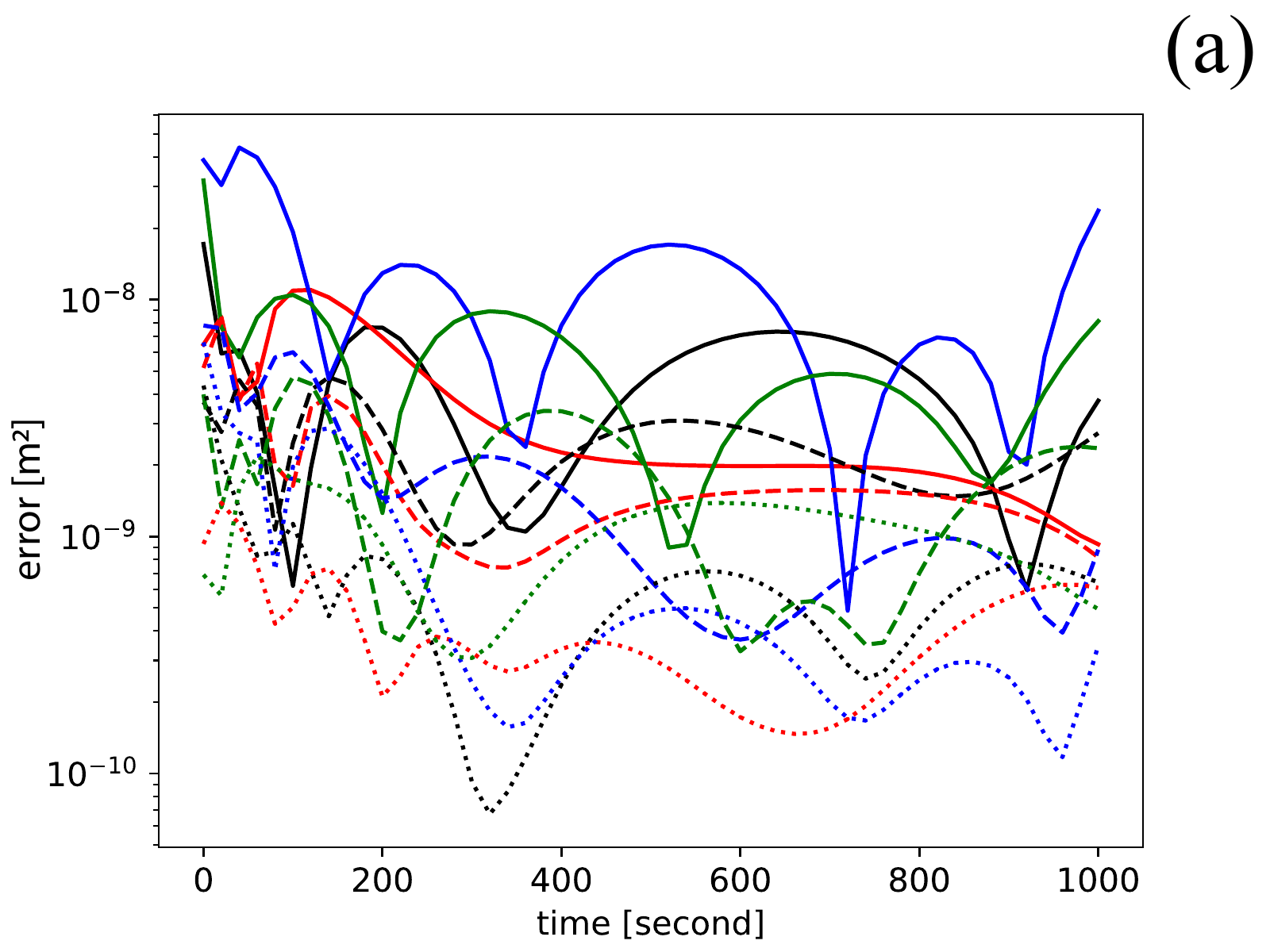}
         \includegraphics[keepaspectratio, height=6.0cm]{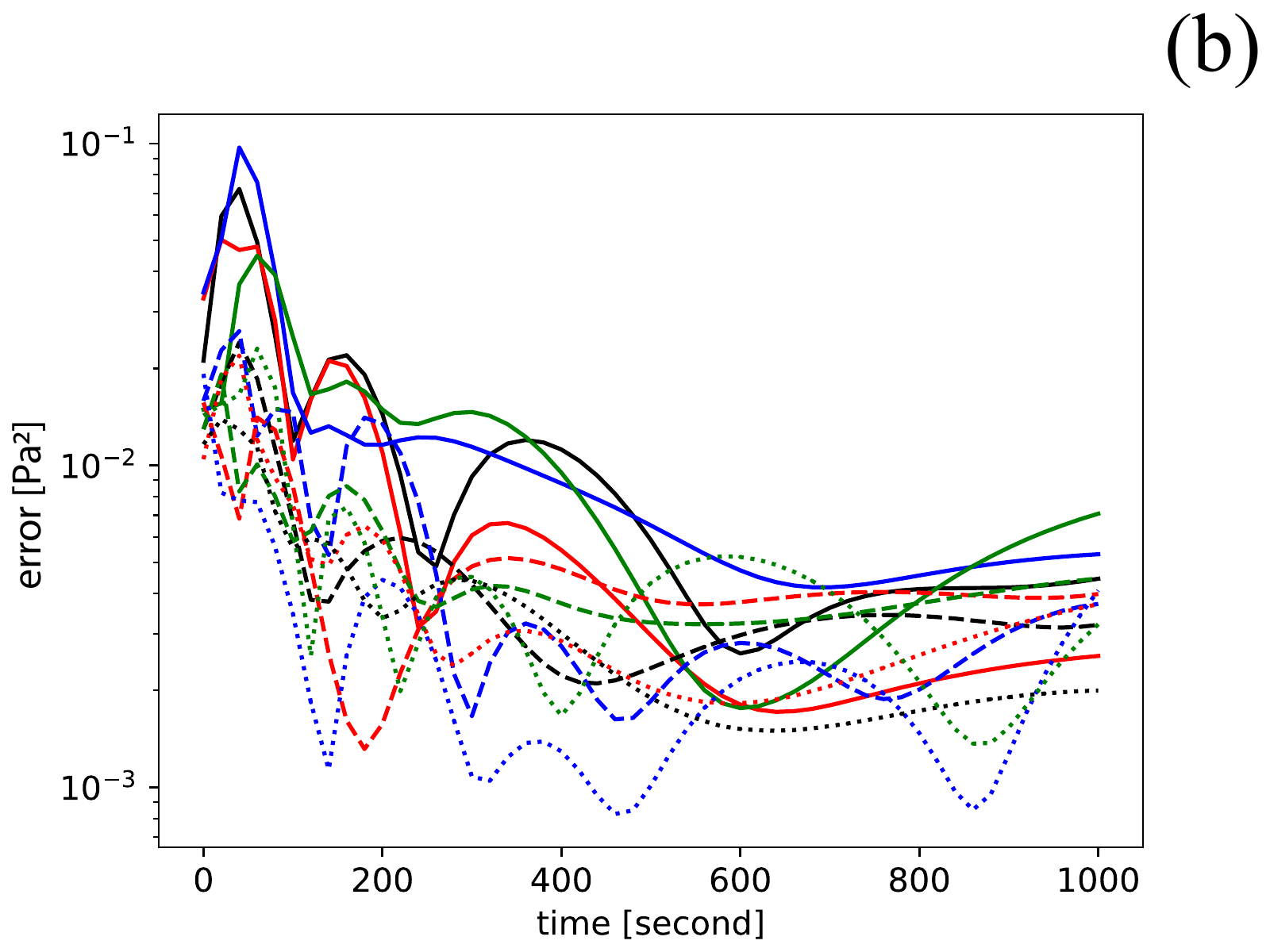}
         \includegraphics[keepaspectratio, height=6.0cm]{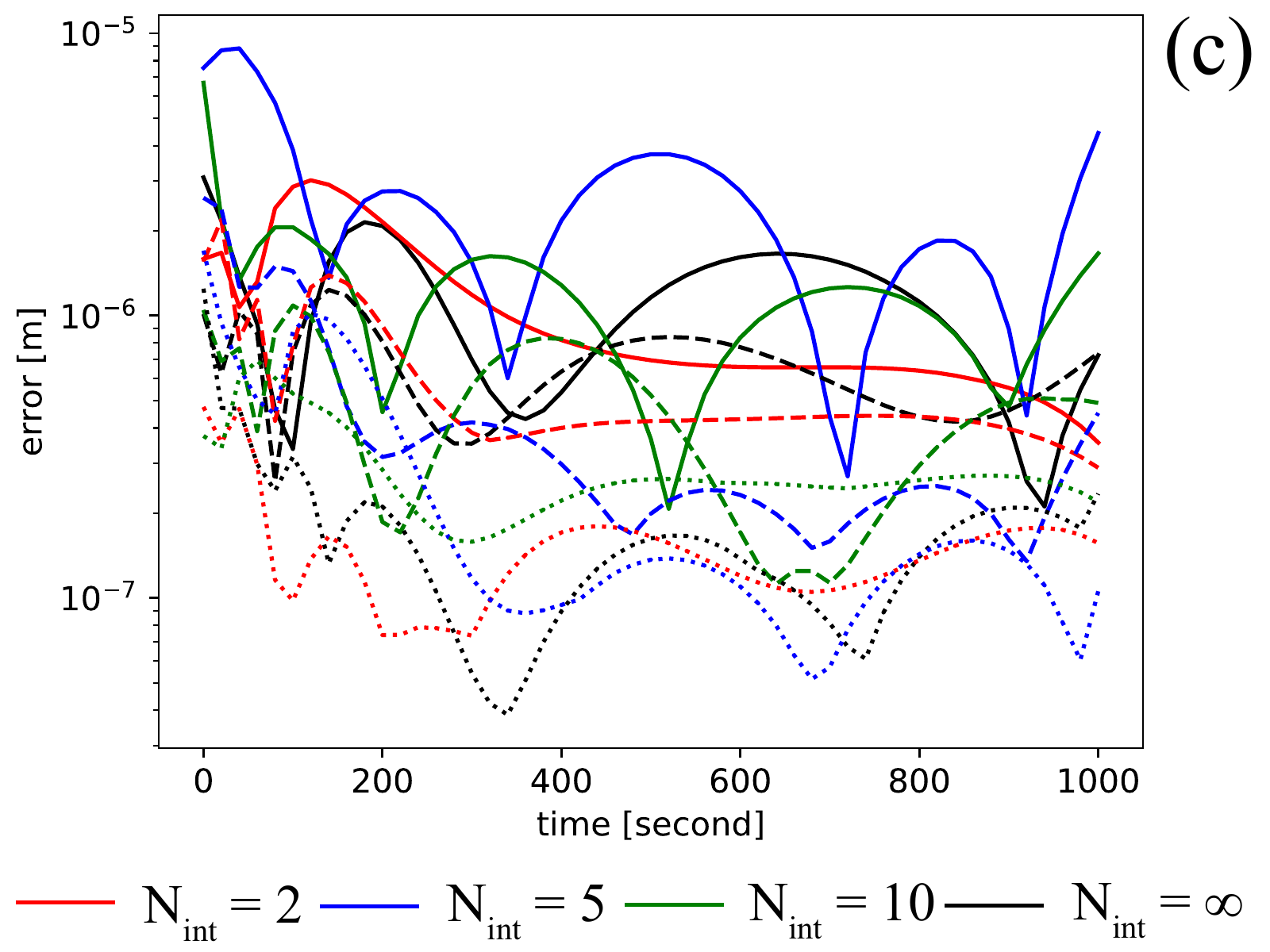}
         \includegraphics[keepaspectratio, height=6.0cm]{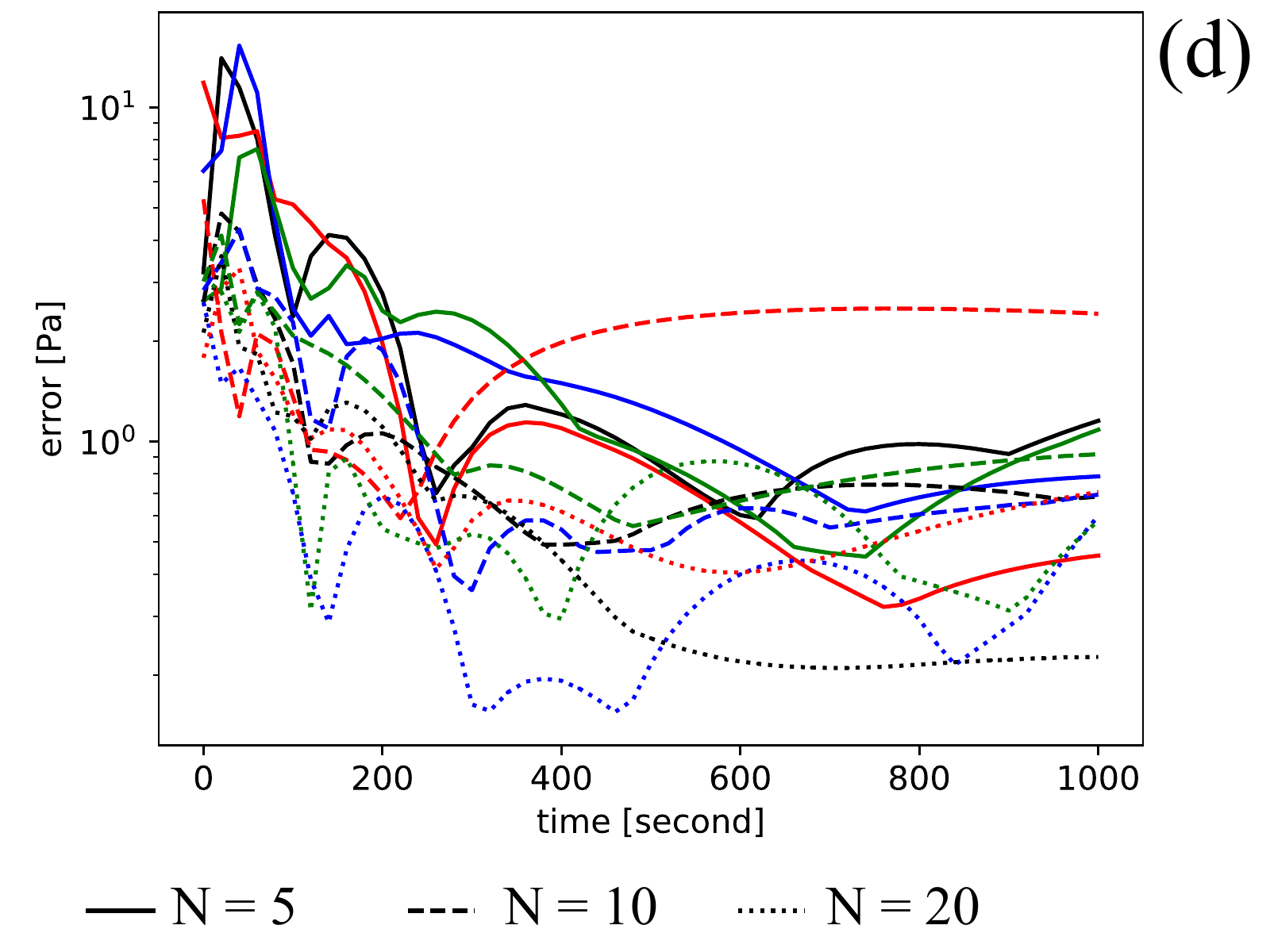}
   \caption{Example 4: Errors of reconstruction solutions using ${\bm{\mu} = (\nu, \alpha) = (0.2, 0.5)}$ - outside of the training snapshots and $\widehat{\bm{\theta}}^u$, $\widehat{\bm{\theta}}^p$, using different numbers of reduced basis ($\mathrm{N}$): (a) mean squared error (MSE) of displacement field ($\bm{u}$), (b) mean squared error (MSE) of fluid pressure field ($p$), (c) maximum error (ME) of displacement field ($\bm{u}$) (d) maximum error (ME) of fluid pressure field ($p$). Colors correspond to increasing values of $\mathrm{N_{int}}$; solid, dashed, and dotted lines represent $\mathrm{N} = 5$, $\mathrm{N} = 10$, and $\mathrm{N} = 20$ cases,  respectively.}
   \label{fig:ex4_source_err_3}
\end{figure}

\begin{enumerate}
    \item investigate MSE and ME results of ROM framework for a realization $\bm{\mu}$ \emph{in the training set}, and the coefficients ${\bm{\theta}}^u(t, \bm{\mu})$ and ${\bm{\theta}}^p(t, \bm{\mu})$ determined from the $L^2$ projection, rather than the ones from ANN, see Figure \ref{fig:ex4_source_err_1};

    \item investigate MSE and ME results of ROM framework for a realization $\bm{\mu}$ \emph{in the training set}, and $\widehat{\bm{\theta}}^u(t, \bm{\mu})$ and $\widehat{\bm{\theta}}^p(t, \bm{\mu})$ obtained by means of ANN, see Figure \ref{fig:ex4_source_err_2};

    \item investigate MSE and ME results of ROM framework for a realization $\bm{\mu}$ \emph{outside of the training set}, and $\widehat{\bm{\theta}}^u(t, \bm{\mu})$ and $\widehat{\bm{\theta}}^p(t, \bm{\mu})$ obtained by means of ANN, see Figure \ref{fig:ex4_source_err_3}.
\end{enumerate}

For what concerns the first goal, we select $\bm{\mu} = (\nu, \alpha) = (0.1, 0.8)$ in the training set and reuse coefficients ${\bm{\theta}}^u(t, \bm{\mu})$ and ${\bm{\theta}}^p(t, \bm{\mu})$ obtained by the $L^2$ projection. The corresponding results for MSE and ME indices, and for both primal variables, are presented in Figure \ref{fig:ex4_source_err_1}. Different colors correspond to different values of $\mathrm{N_{int}}$, including the label $\mathrm{N_{int}} = \infty$, which represents the use of the standard POD. Different line styles (solid, dashed, and dotted) correspond to increasing dimension $\mathrm{N}$ of the reduced basis spaces. As expected, the ROM accuracy increases as the number of $\mathrm{N}$ increase, following an exponential trend. As we increase the number of $\mathrm{N_{int}}$, the MSE and ME behaviors approach the ones of $\mathrm{N_{int}} = \infty$ case, which means that even a nested POD compression with moderate value of $\mathrm{N_{int}}$ is able to correctly capture the time evoluation. Besides, the values of MSE and ME initially decrease and remain constant over time since our problem reaches the steady-state solution.

For the second goal, the MSE and ME results using $\bm{\mu}$ in the training set and $\widehat{\bm{\theta}}^u$, $\widehat{\bm{\theta}}^p$ predicted by the ANN are shown in Figure \ref{fig:ex4_source_err_2}. For a fair comparison to \ref{fig:ex4_source_err_1} we use the same parameter instance $\bm{\mu} = (\nu, \alpha) = (0.1, 0.8)$ in the training set, to compare different approximation properties stemming from the use ${\bm{\theta}}^u$, ${\bm{\theta}}^p$ (Figure \ref{fig:ex4_source_err_1}) and $\widehat{\bm{\theta}}^u$, $\widehat{\bm{\theta}}^p$ (Figure \ref{fig:ex4_source_err_2}).

%

Compared to the previous results shown in Figure \ref{fig:ex4_source_err_1}, the MSE and ME values in Figure \ref{fig:ex4_source_err_2} are approximately three orders of magnitude higher. Moreover, there is no clear trend in how the increase of $\mathrm{N}$ and $\mathrm{N_{int}}$ affects the MSE and ME results. Indeed, as the number of $\mathrm{N_{int}}$ increase, we could not observe that the MSE and ME behaviors approach the ones of $\mathrm{N_{int}} = \infty$ case. The MSE and ME values, however, still decreases as the time domain progresses.
To better quantify the observations we obtain from Figure \ref{fig:ex4_source_err_2}, we take an average for all time steps of the MSE of the fluid pressure field ($p$) and present it in Table \ref{tab:l2_source_err_2}. We could see that the average MSE remains of the same order of magnitude for any $\mathrm{N}$ and $\mathrm{N_{int}}$ pair, which is a marked difference from the MSE results shown in Figure \ref{fig:ex4_source_err_1}. This indicates that to improve the accuracy, one should also consider various other ROM properties, such as the number of snapshots, hyperparameters, or network architecture. Each of these options will be explored in later subsections.

\begin{table}[!ht]
\centering
\caption{Example 4: Average for all time step of mean squared error (MSE) of fluid pressure field ($p$) presented in Figure \ref{fig:ex4_source_err_2}.}
\begin{tabular}{|l|c|c|c|c|}
\hline
   &  \multicolumn{1}{l|}{$\mathrm{N_{int}} = 2$} & \multicolumn{1}{l|}{$\mathrm{N_{int}} = 5$} & \multicolumn{1}{l|}{$\mathrm{N_{int}} = 10$} &
   \multicolumn{1}{l|}{$\mathrm{N_{int}} = \infty$} \\ \hline
$\mathrm{N} = 5$                           & 0.0053                             & 0.0057                             & 0.0050  & 0.0051                               \\ \hline
$\mathrm{N} = 10$                            & 0.0045                             & 0.0042                             & 0.0052 & 0.0066                             \\ \hline
$\mathrm{N} = 20$                          & 0.0049                             & 0.0073                             & 0.0073                        &     0.0052      \\ \hline
\end{tabular}
\label{tab:l2_source_err_2}
\end{table}

We then move to the third goal, where we present the MSE and ME results using $\bm{\mu}$ outside of the training set and $\widehat{\bm{\theta}}^u$, $\widehat{\bm{\theta}}^p$ predicted by the ANN, as shown in Figure \ref{fig:ex4_source_err_3}. These results are comparable to ones presented in Figure \ref{fig:ex4_source_err_2} as there is no clear relationship between the MSE or ME values and the numbers of $\mathrm{N}$ or $\mathrm{N_{int}}$. Furthermore, the MSE and ME values are approximately three orders of magnitude higher than those presented in Figure \ref{fig:ex4_source_err_1}.

The average for all time step of the MSE of the fluid pressure field ($p$) is presented in Table \ref{tab:l2_source_err_3}. Similar to Table \ref{tab:l2_source_err_2}, the average MSE is independent of the numbers of $\mathrm{N}$ and $\mathrm{N_{int}}$. Hence, from Figures \ref{fig:ex4_source_err_1}, \ref{fig:ex4_source_err_2}, and \ref{fig:ex4_source_err_3}, and Tables \ref{tab:l2_source_err_2} and \ref{tab:l2_source_err_3}, we could see that the errors introduced by POD and $L^2$ projection phases are negligible compared to the errors initiated from the ANN phase (prediction of $\widehat{\bm{\theta}}^u$, $\widehat{\bm{\theta}}^p$). Therefore we conclude that in practical applications one may want to choose moderate values for both $\mathrm{N}$ and $\mathrm{N_{int}}$. Keeping a moderate value for $\mathrm{N}$ guarantees the evaluation of a small network during the online phase; keeping a moderate value for $\mathrm{N_{int}}$ results in large computational savings during the offline phase. Once $\mathrm{N}$ and $\mathrm{N_{int}}$ are fixed, if a further increase in accuracy is desired one may then explore the possibility of changing the remain properties of the ROM framework, as we discuss in the following.

\begin{table}[!ht]
\centering
\caption{Example 4: Average for all time step of mean squared error (MSE) of fluid pressure field ($p$) presented in Figure \ref{fig:ex4_source_err_3}.}
\begin{tabular}{|l|c|c|c|c|}
\hline
   &  \multicolumn{1}{l|}{$\mathrm{N_{int}} = 2$} & \multicolumn{1}{l|}{$\mathrm{N_{int}} = 5$} & \multicolumn{1}{l|}{$\mathrm{N_{int}} = 10$} &
   \multicolumn{1}{l|}{$\mathrm{N_{int}} = \infty$} \\ \hline
$\mathrm{N} = 5$                           & 0.0042                             & 0.0035                             & 0.0049  & 0.0027                               \\ \hline
$\mathrm{N} = 10$                            & 0.0084                             & 0.0126                             & 0.0100 & 0.0110                             \\ \hline
$\mathrm{N} = 20$                          & 0.0048                             & 0.0054                             & 0.0049                        &     0.0048      \\ \hline
\end{tabular}
\label{tab:l2_source_err_3}
\end{table}

\subsubsection{Effect of number of snapshots} \label{sec:eff_snapshot}

The number of snapshots' ($\mathrm{M}$) effect is studied by comparing the MSE and ME results of cases with different $\mathrm{M}$. To reiterate, the higher $\mathrm{M}$, the longer the ROM framework will take to perform FOM solves and POD compressions. For instance, the wall time used to solve all FOM problems is 1780, 7120,  and 16020 seconds for $\mathrm{M} = 100$, $\mathrm{M} = 400$, and $\mathrm{M} = 900$, respectively. The wall time used corresponding to the cases with different $\mathrm{M}$ is presented in Table \ref{tab:time_num} for the POD compressions.

\begin{figure}[!ht]
   \centering

         \includegraphics[keepaspectratio, height=6.0cm]{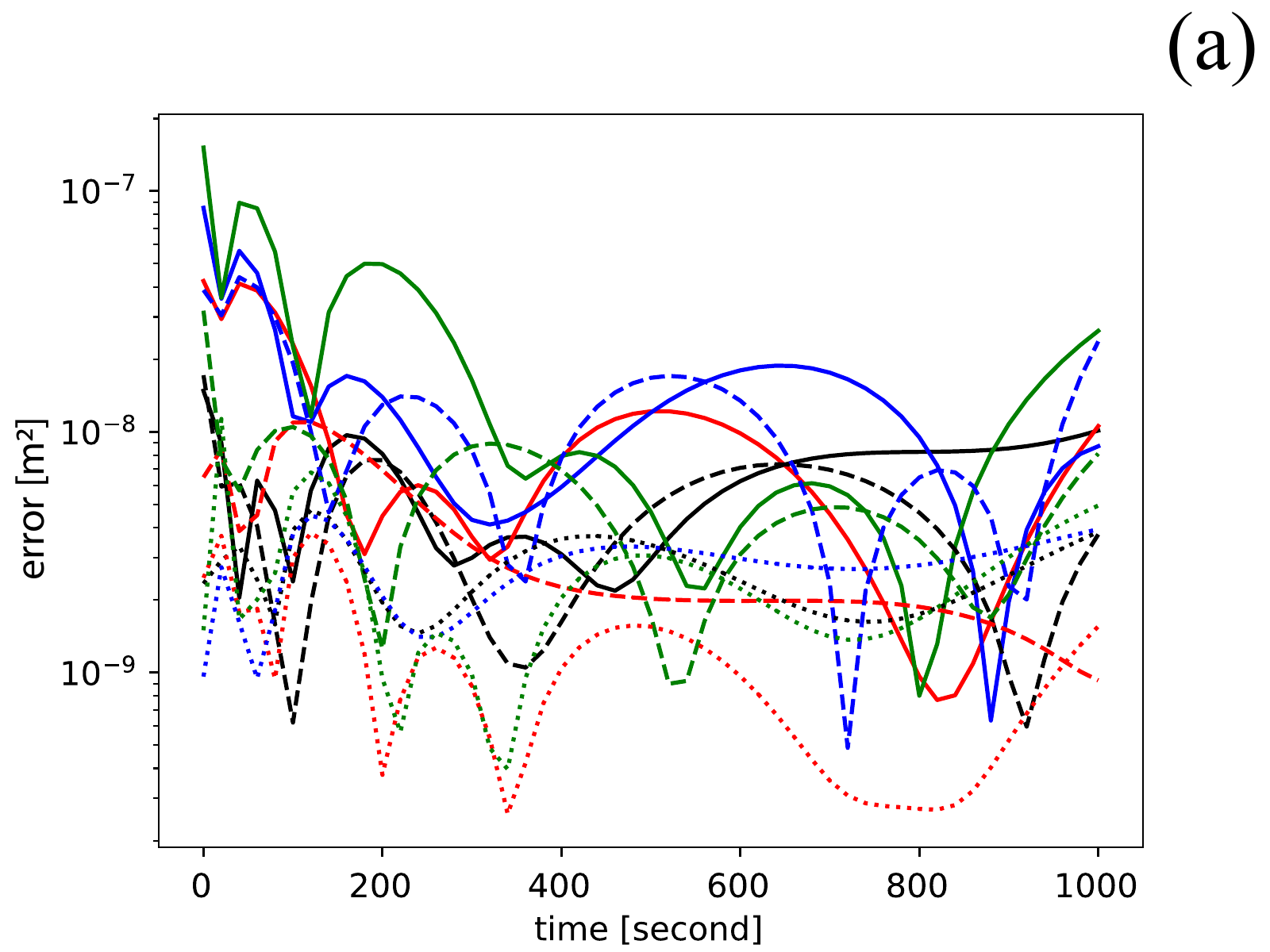}
         \includegraphics[keepaspectratio, height=6.0cm]{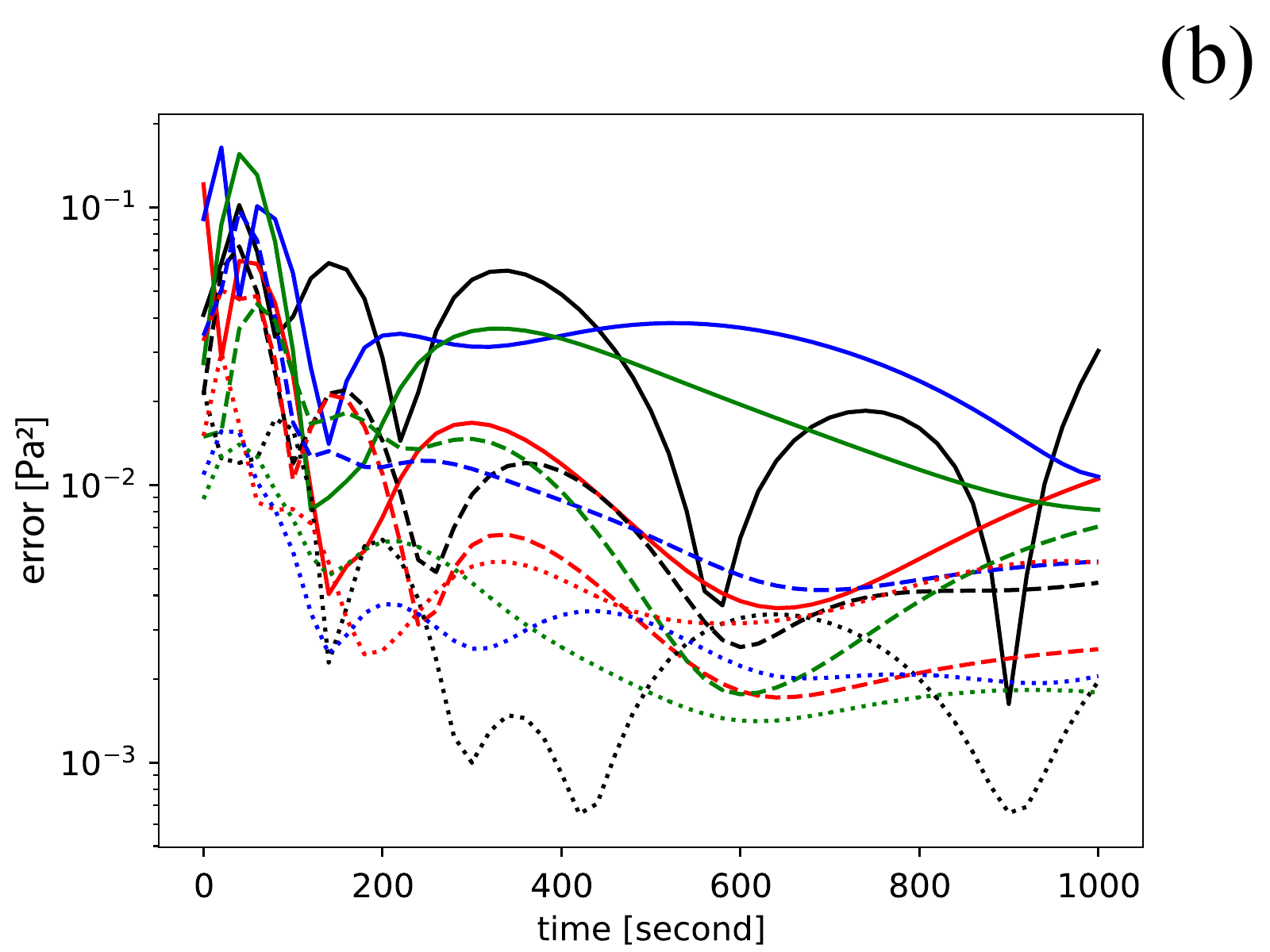}
         \includegraphics[keepaspectratio, height=6.0cm]{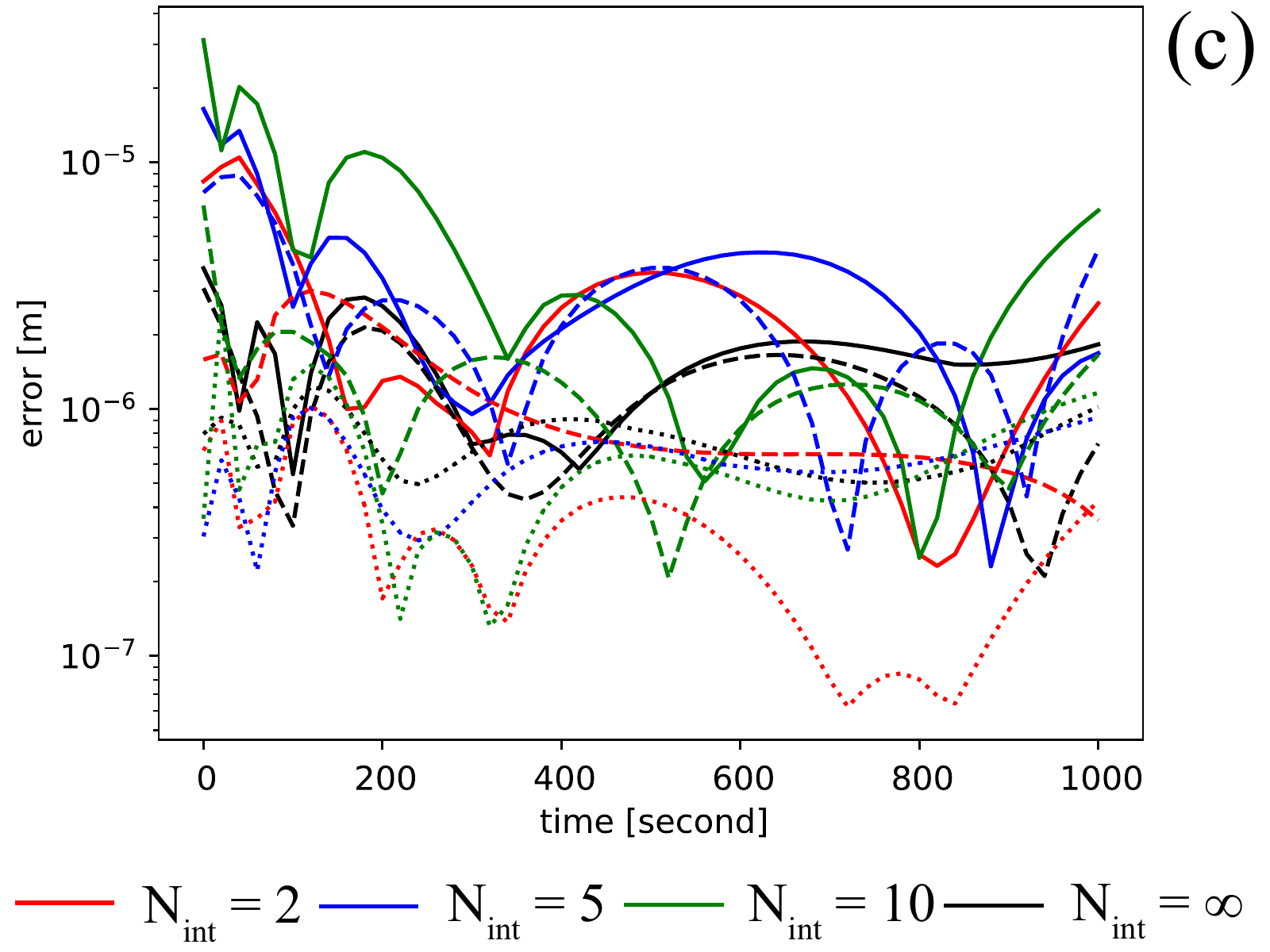}
         \includegraphics[keepaspectratio, height=6.0cm]{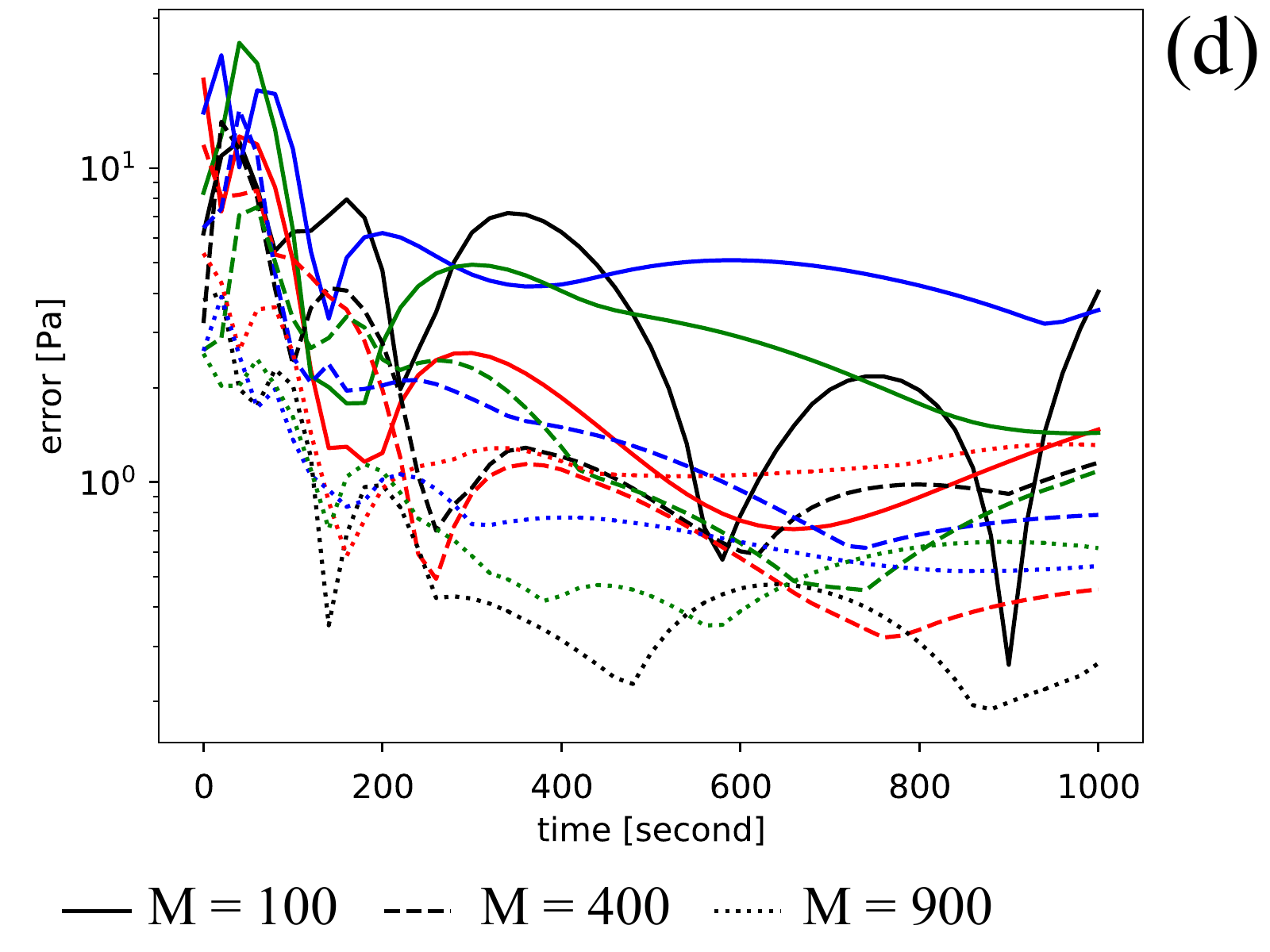}
   \caption{Example 4: Errors of reconstruction solutions using $\bm{\mu} = (\nu, \alpha) = (0.2, 0.5)$ - outside of the training snapshots and $\widehat{\bm{\theta}}^u$, $\widehat{\bm{\theta}}^p$, using different number of snapshots ($\mathrm{M}$): (a) mean squared error (MSE) of displacement field ($\bm{u}$), (b) mean squared error (MSE) of fluid pressure field ($p$), (c) maximum error (ME) of displacement field ($\bm{u}$) (d) maximum error (ME) of fluid pressure field ($p$). Colors correspond to increasing values of $\mathrm{N_{int}}$; solid, dashed, and dotted lines represent $\mathrm{M} = 100$, $\mathrm{M} = 400$, and $\mathrm{M} = 900$ cases,  respectively. Note that we fix $\mathrm{N} = 10$, $\mathrm{N_{hl}} = 3$, and $\mathrm{N_{nn}} = 7$.}
   \label{fig:ex4_num_snap}

\end{figure}

The MSE and ME results with different number of snapshots ($\mathrm{M}$) are presented in Figure \ref{fig:ex4_num_snap}. The MSE averages for all time step of $p$ is also presented in Table \ref{tab:l2_num}. We fix $\mathrm{N} = 10$, $\mathrm{N_{hl}} = 3$, and $\mathrm{N_{nn}} = 7$. From Figure \ref{fig:ex4_num_snap} and Table \ref{tab:l2_num}, we observe that the model with highest $\mathrm{M}$ provides the lowest MSE and ME results. However, there is no distinct different of the MSE and ME results among the models using different $\mathrm{N_{int}}$. Therefore an actionable way of increasing the ROM accuracy is to provide more input data to the training phase of the ANN; computational savings related to the basis generation during the offline stage can still be achieved by choosing a moderate value for $\mathrm{N_{int}}$.

\begin{table}[!ht]
\centering
\caption{Example 4: Average for all time step of mean squared error (MSE) of fluid pressure field ($p$) presented in Figure \ref{fig:ex4_num_snap}.}
\begin{tabular}{|l|c|c|c|c|}
\hline
   &  \multicolumn{1}{l|}{$\mathrm{N_{int}} = 2$} & \multicolumn{1}{l|}{$\mathrm{N_{int}} = 5$} & \multicolumn{1}{l|}{$\mathrm{N_{int}} = 10$} &
   \multicolumn{1}{l|}{$\mathrm{N_{int}} = \infty$} \\ \hline
$\mathrm{M} = 100$                           & 0.0139                             & 0.0361                             & 0.0271  & 0.0299                               \\ \hline
$\mathrm{M} = 400$                            & 0.0084                             & 0.0126                             & 0.0100 & 0.0110                             \\ \hline
$\mathrm{M} = 900$                          & 0.0054                             & 0.0036                             & 0.0037                        &     0.0039      \\ \hline
\end{tabular}
\label{tab:l2_num}
\end{table}

\subsubsection{Effect of network architecture} \label{sec:eff_nhl}

\begin{figure}[!ht]
   \centering

         \includegraphics[keepaspectratio, height=6.0cm]{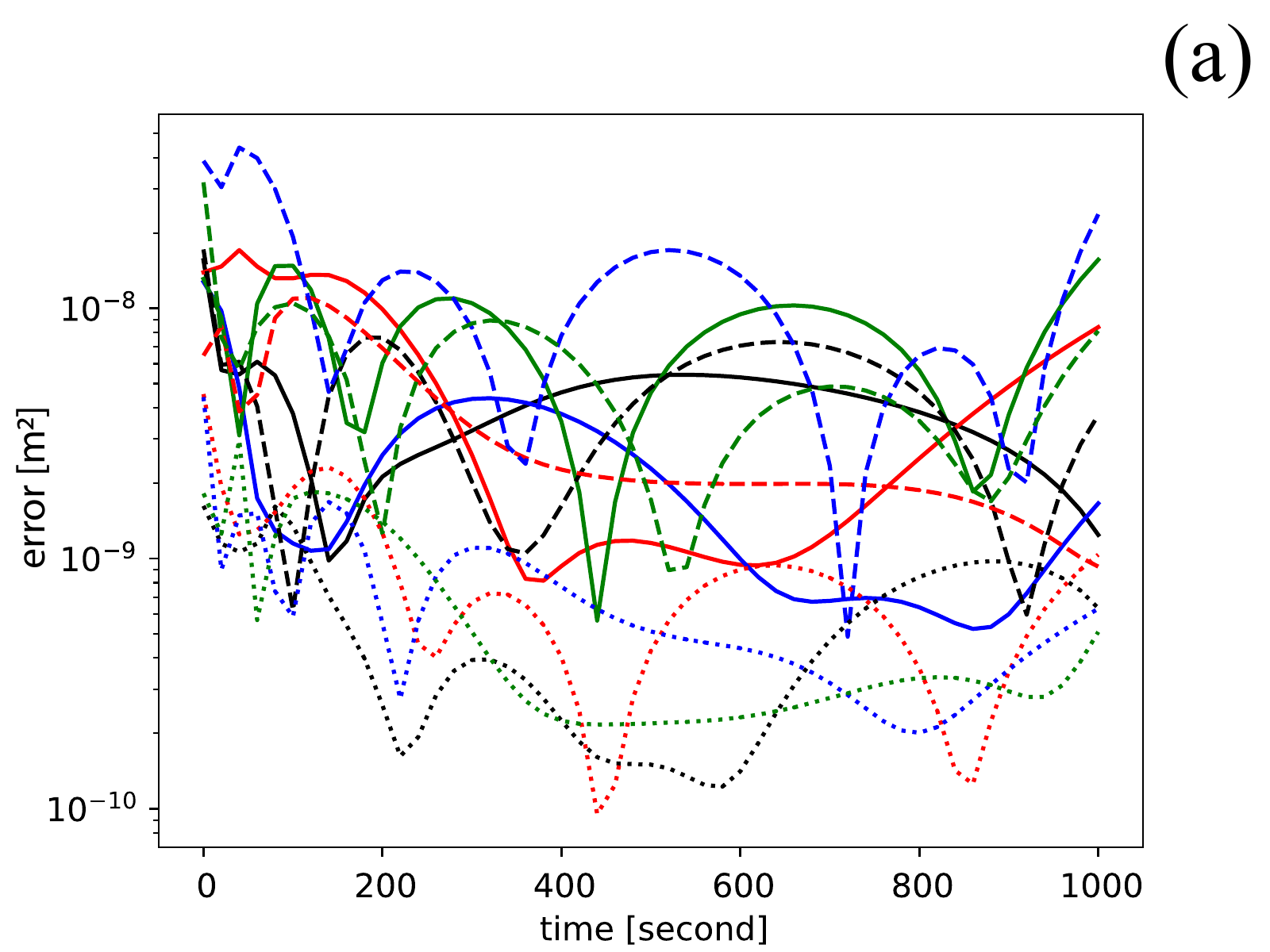}
         \includegraphics[keepaspectratio, height=6.0cm]{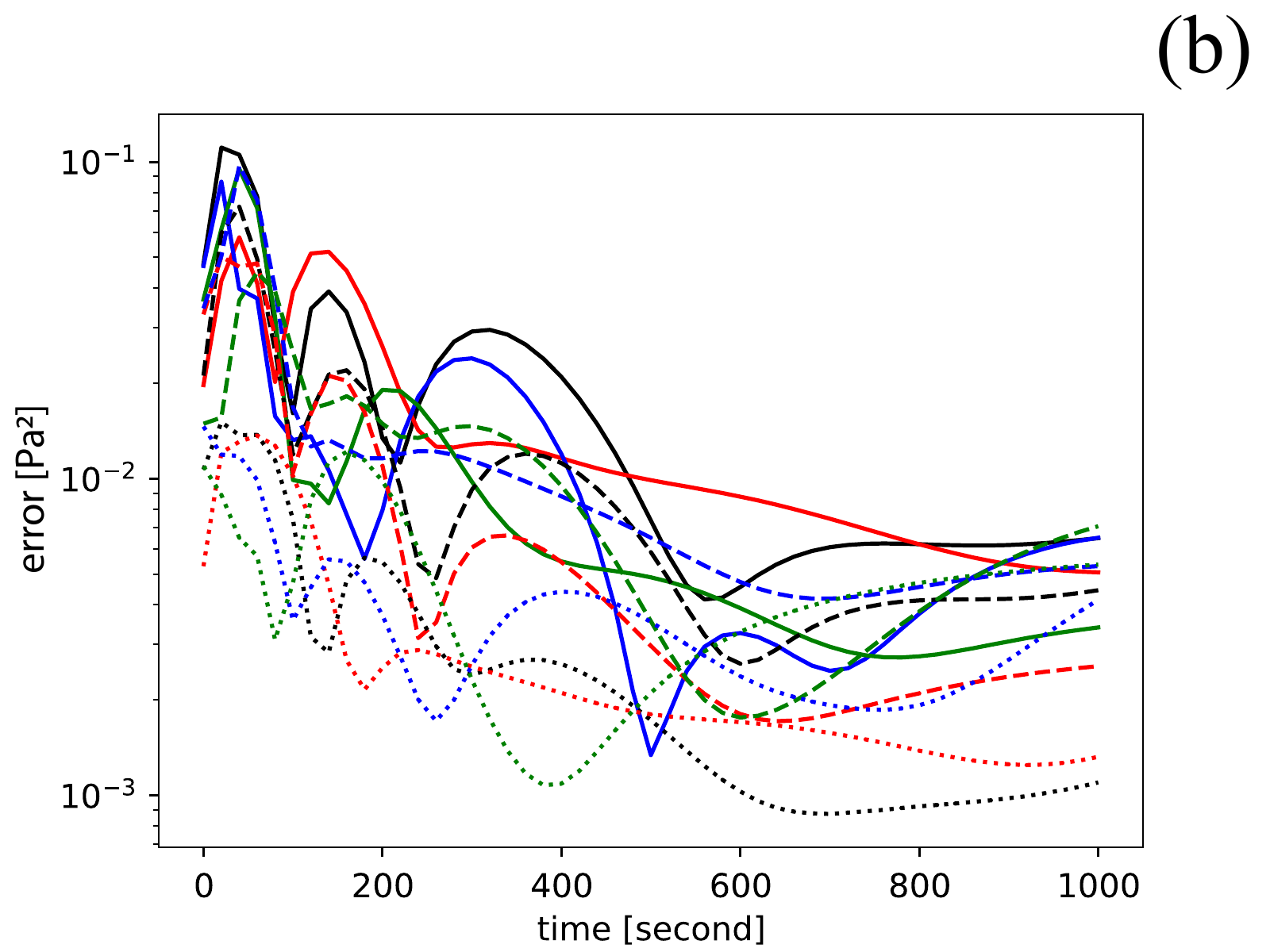}
         \includegraphics[keepaspectratio, height=6.0cm]{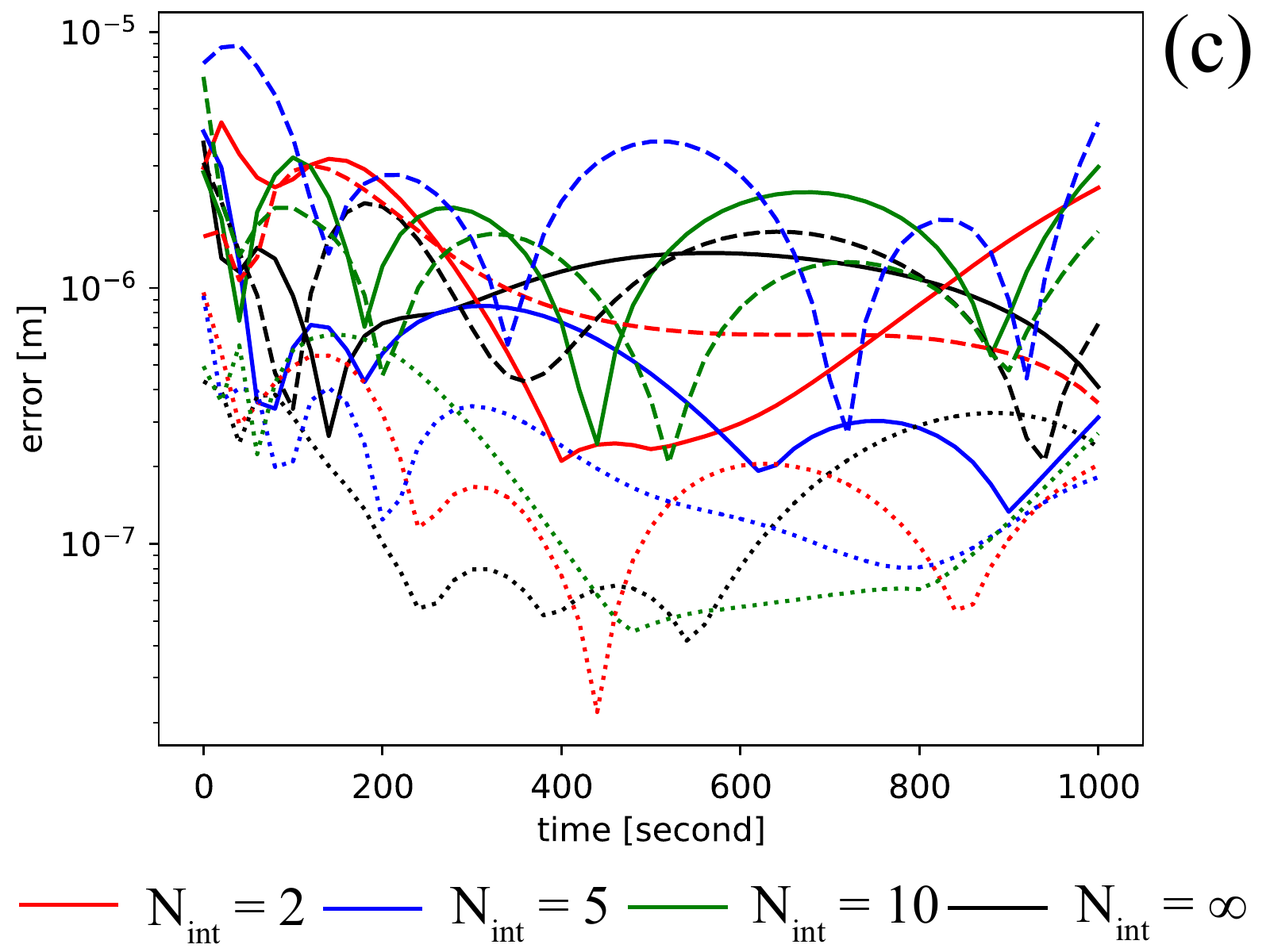}
         \includegraphics[keepaspectratio, height=6.0cm]{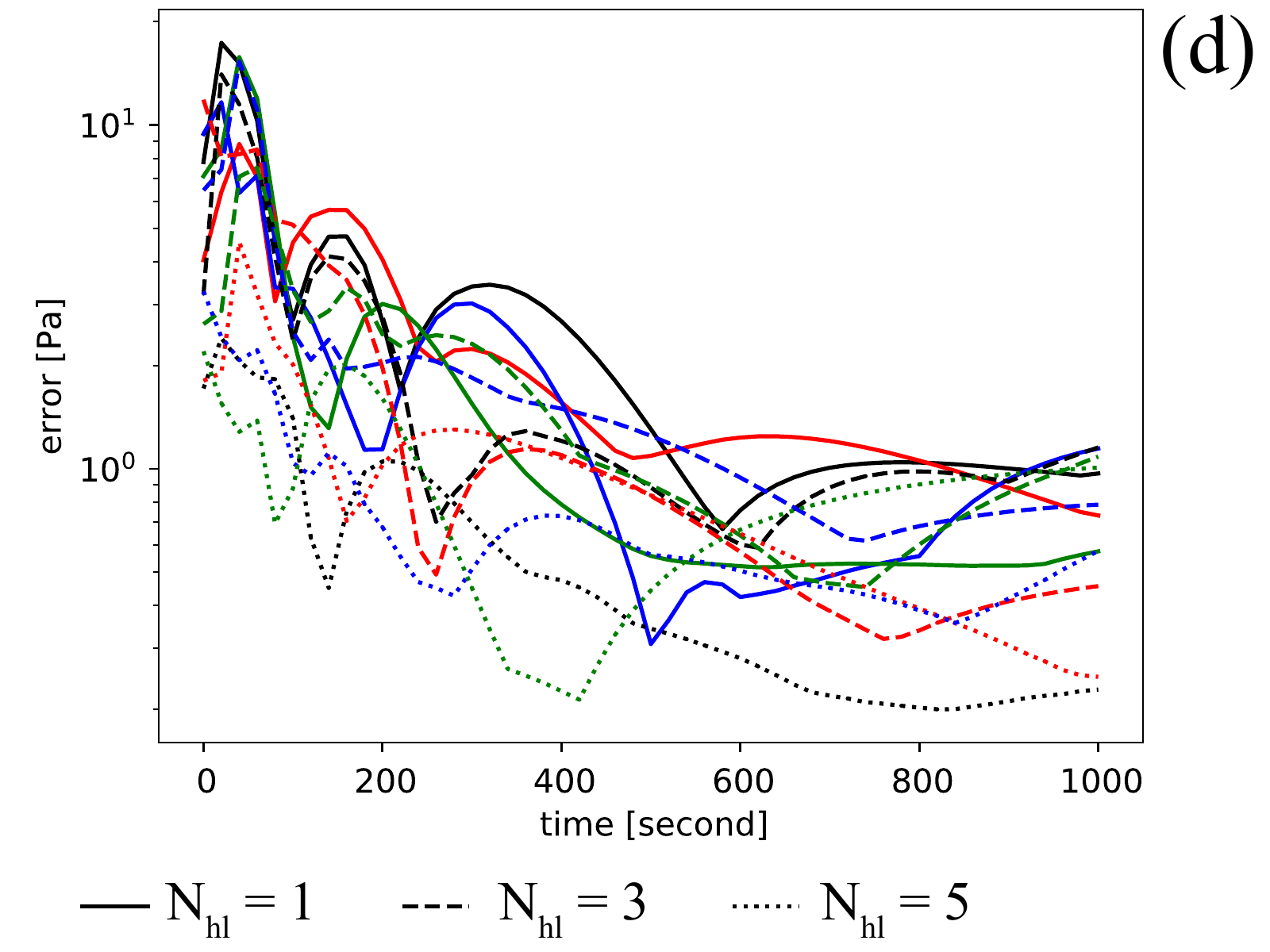}
   \caption{Example 4: Errors of reconstruction solutions using $\bm{\mu} = (\nu, \alpha) = (0.2, 0.5)$ - outside of the training snapshots and $\widehat{\bm{\theta}}^u$, $\widehat{\bm{\theta}}^p$, using different number of hidden layers ($\mathrm{N_{hl}}$): (a) mean squared error (MSE) of displacement field ($\bm{u}$), (b) mean squared error (MSE) of fluid pressure field ($p$), (c) maximum error (ME) of displacement field ($\bm{u}$) (d) maximum error (ME) of fluid pressure field ($p$). Colors correspond to increasing values of $\mathrm{N_{int}}$; solid, dashed, and dotted lines represent $\mathrm{N_{hl}} = 1$, $\mathrm{N_{hl}} = 3$, and $\mathrm{N_{hl}} = 5$ cases,  respectively. Note that we fix $\mathrm{N} = 10$ and $\mathrm{N_{nn}} = 7$.}
   \label{fig:ex4_num_nhl}

\end{figure}

We then examine the effect of network architecture (i.e., number of hidden layers ($\mathrm{N_{hl}}$) and number of neurons per hidden layer ($\mathrm{N_{nn}}$)). We begin with cases where we fix $\mathrm{N} = 10$ and $\mathrm{N_{nn}} = 7$, but we vary $\mathrm{N_{hl}}$ as presented in Figure \ref{fig:ex4_num_nhl} and Table \ref{tab:l2_hl}. From these MSE and ME results, we observe that the MSE and ME values decrease as $\mathrm{N_{hl}}$ increases. Similar to previous cases, the MSE and ME values, however, seems to be independent of a choice of $\mathrm{N_{int}}$.

\begin{table}[!ht]
\centering
\caption{Example 4: Average for all time step of mean squared error (MSE) of fluid pressure field ($p$) presented in Figure \ref{fig:ex4_num_nhl}.}
\begin{tabular}{|l|c|c|c|c|}
\hline
   &  \multicolumn{1}{l|}{$\mathrm{N_{int}} = 2$} & \multicolumn{1}{l|}{$\mathrm{N_{int}} = 5$} & \multicolumn{1}{l|}{$\mathrm{N_{int}} = 10$} &
   \multicolumn{1}{l|}{$\mathrm{N_{int}} = \infty$} \\ \hline
$\mathrm{N_{hl}} = 1$                           & 0.0155              & 0.0117                   &  0.0115                              & 0.0191                               \\ \hline
$\mathrm{N_{hl}} = 3$                            & 0.0084                             & 0.0126                             & 0.0100 & 0.0110                             \\ \hline
$\mathrm{N_{hl}} = 5$                          & 0.0031                             & 0.0038                             & 0.0047                        &     0.0031      \\ \hline
\end{tabular}
\label{tab:l2_hl}
\end{table}

The wall time as a function of $\mathrm{N_{hl}}$ and $\mathrm{N_{int}}$ is presented in Table \ref{tab:time_nhl}. As one expects, as the number of $\mathrm{N_{hl}}$ grows, the longer time the model takes to train during the ANN phase. Moreover, this table shows that the $\mathrm{N_{int}}$ does not affect the computational cost of the ANN phase. As mentioned in the methodology section, the wall time shown in Table \ref{tab:time_nhl} is a combination of the wall time used to perform $L^2$ projection and the wall time used to train the ANN. We combine these two operations because the wall time used to perform $L^2$ projection is relatively much smaller than the wall time used to train the ANN as well as other phases.

\begin{table}[!ht]
\centering
\caption{Example 4: Comparison of the wall time (seconds) used for training the neural networks with different number of hidden layers ($\mathrm{N_{hl}}$) and number of intermediate reduced basis ($\mathrm{N_{int}}$). $\mathrm{N_{int}} = \infty$ represents a case where we do not use the nested POD technique.}
\begin{tabular}{|l|c|c|c|c|}
\hline
   &  \multicolumn{1}{l|}{$\mathrm{N_{int}} = 2$} & \multicolumn{1}{l|}{$\mathrm{N_{int}} = 5$} & \multicolumn{1}{l|}{$\mathrm{N_{int}} = 10$} &
   \multicolumn{1}{l|}{$\mathrm{N_{int}} = \infty$} \\ \hline
$\mathrm{N_{hl}} = 1$                           & 6217             & 5963                   &  6065                              & 6039                               \\ \hline
$\mathrm{N_{hl}} = 3$                            & 7114                             & 7108                             & 7064 & 7103                             \\ \hline
$\mathrm{N_{hl}} = 5$                          & 8279                             & 8271                             & 8340                        &     8152      \\ \hline
\end{tabular}
\label{tab:time_nhl}
\end{table}

\begin{figure}[!ht]
   \centering

         \includegraphics[keepaspectratio, height=6.0cm]{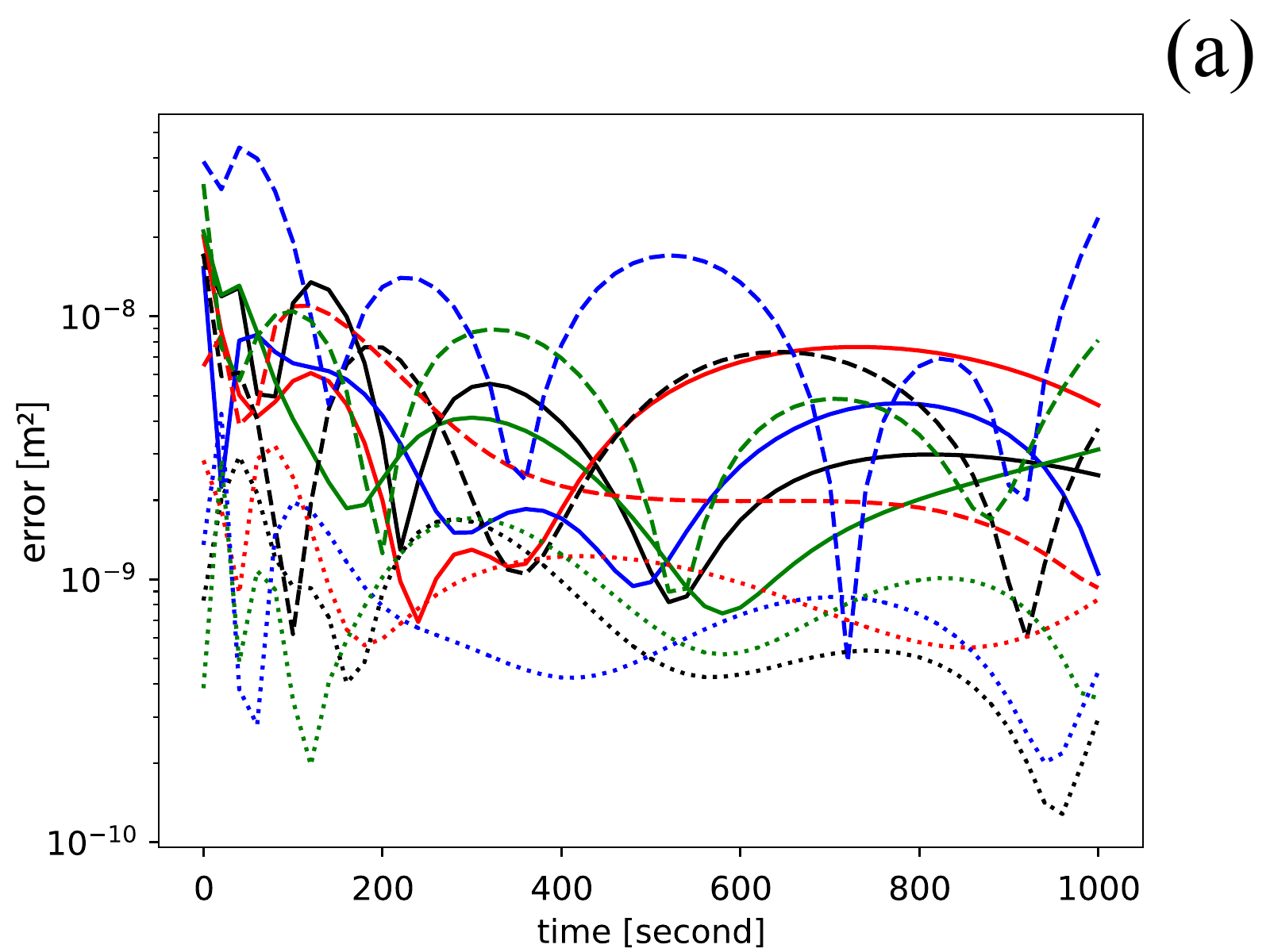}
         \includegraphics[keepaspectratio, height=6.0cm]{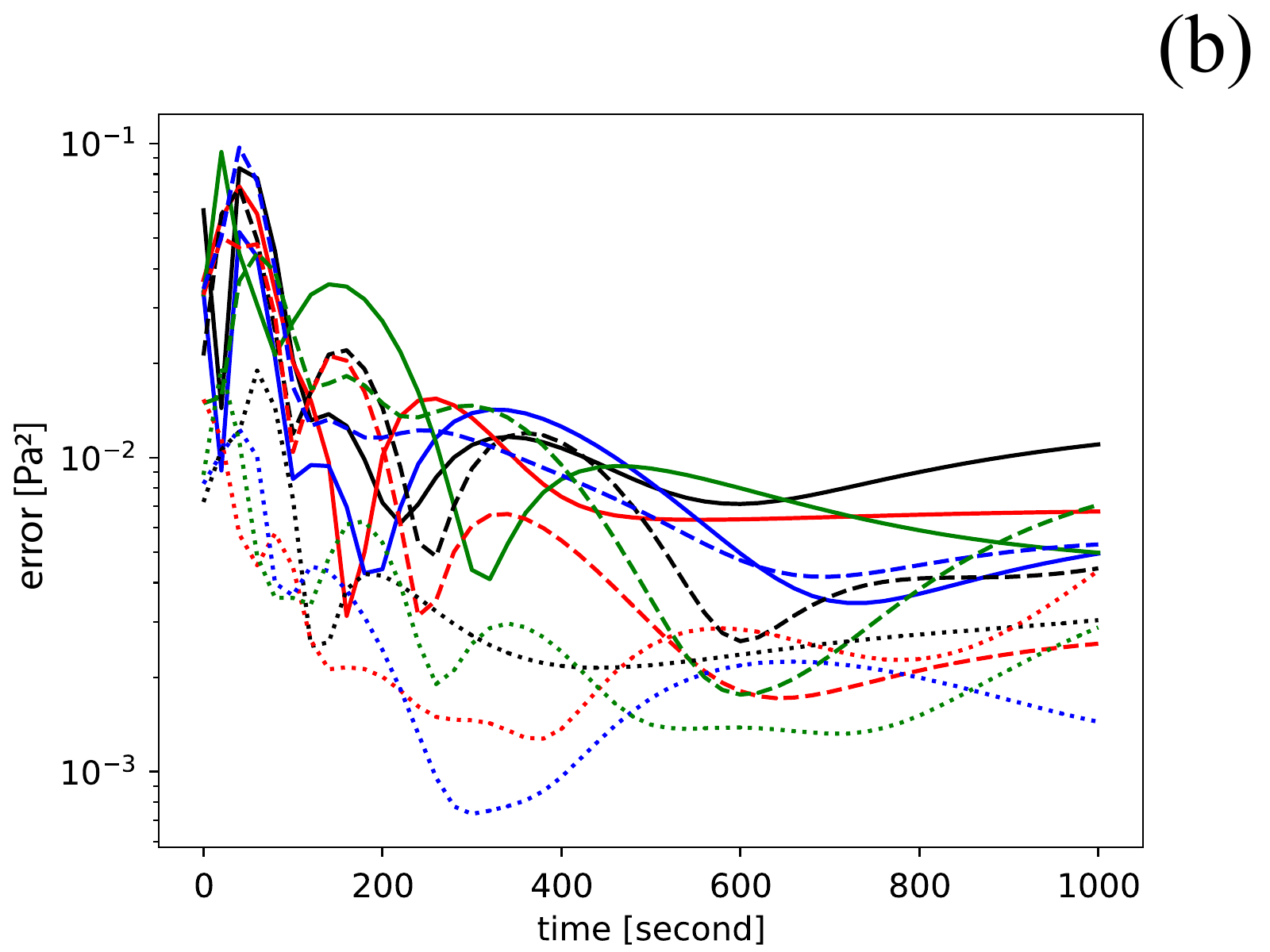}
         \includegraphics[keepaspectratio, height=6.0cm]{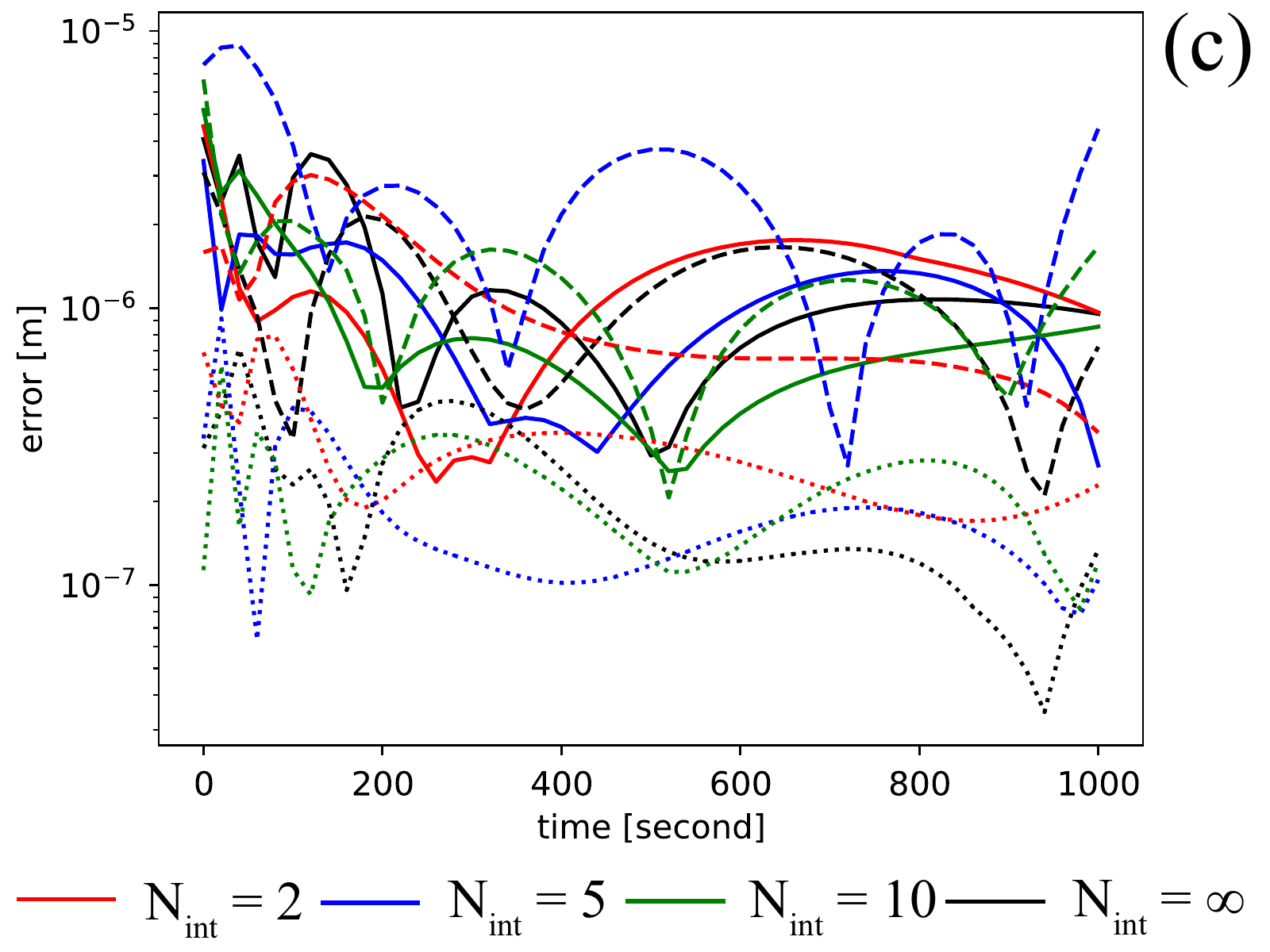}
         \includegraphics[keepaspectratio, height=6.0cm]{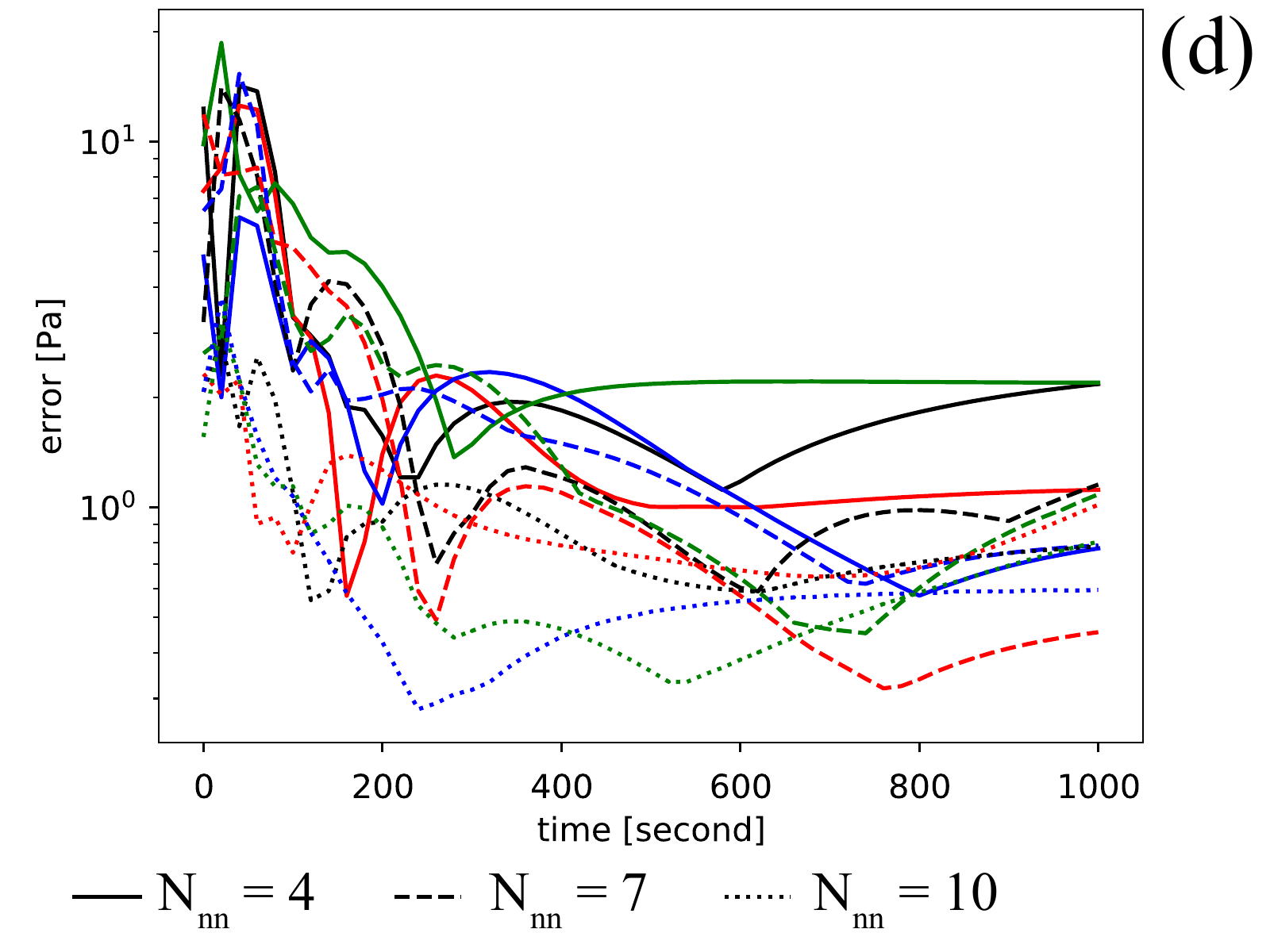}
   \caption{Example 4: Errors of reconstruction solutions using $\bm{\mu} = (\nu, \alpha) = (0.2, 0.5)$ - outside of the training snapshots and $\widehat{\bm{\theta}}^u$, $\widehat{\bm{\theta}}^p$, using different number of neuron per layer ($\mathrm{N_{nn}}$): (a) mean squared error (MSE) of displacement field ($\bm{u}$), (b) mean squared error (MSE) of fluid pressure field ($p$), (c) maximum error (ME) of displacement field ($\bm{u}$) (d) maximum error (ME) of fluid pressure field ($p$). Colors correspond to increasing values of $\mathrm{N_{int}}$; solid, dashed, and dotted lines represent $\mathrm{N_{nn}} = 4$, $\mathrm{N_{nn}} = 7$, and $\mathrm{N_{nn}} = 10$ cases, respectively. Note that we fix $\mathrm{N} = 10$ and $\mathrm{N_{hl}} = 3$.}
   \label{fig:ex4_num_nn}

\end{figure}

We then investigate the MSE and ME results of cases where we vary $\mathrm{N_{nn}}$, but we fix $\mathrm{N} = 10$ and $\mathrm{N_{hl}} = 3$ as presented in Figure \ref{fig:ex4_num_nn} and Table \ref{tab:l2_nn}. Again, we could not observe any clear relationships between the errors and $\mathrm{N_{int}}$. We, however, could see that the model accuracy is improved as we increase $\mathrm{N_{nn}}$.

\begin{table}[!ht]
\centering
\caption{Example 4: Average for all time step of mean squared error (MSE) of fluid pressure field ($p$) presented in Figure \ref{fig:ex4_num_nn}.}
\begin{tabular}{|l|c|c|c|c|}
\hline
   &  \multicolumn{1}{l|}{$\mathrm{N_{int}} = 2$} & \multicolumn{1}{l|}{$\mathrm{N_{int}} = 5$} & \multicolumn{1}{l|}{$\mathrm{N_{int}} = 10$} &
   \multicolumn{1}{l|}{$\mathrm{N_{int}} = \infty$} \\ \hline
$\mathrm{N_{nn}} = 4$                           & 0.0126             & 0.0096                   &  0.0140                             & 0.0142                               \\ \hline
$\mathrm{N_{nn}} = 7$                            & 0.0084                             & 0.0126                            & 0.0100 & 0.0110                            \\ \hline
$\mathrm{N_{nn}} = 10$                          & 0.0030                             & 0.0026                             & 0.0031                        &     0.0038      \\ \hline
\end{tabular}
\label{tab:l2_nn}
\end{table}

The comparison of the wall time used for the $L^2$ projection phase and training the neural networks with a different number of $\mathrm{N_{nn}}$ and number of $\mathrm{N_{int}}$ is presented in Table \ref{tab:time_nn}. Similar to Table \ref{tab:time_nhl}, the number of $\mathrm{N_{int}}$ does not significantly affect the wall time. The $\mathrm{N_{nn}}$, on the other hand, could influence the wall time used to train the ANN.

\begin{table}[!ht]
\centering
\caption{Example 4: Comparison of the wall time (seconds) used for training the neural networks with different number of hidden layers ($\mathrm{N_{nn}}$) and number of intermediate reduced basis ($\mathrm{N_{int}}$). $\mathrm{N_{int}} = \infty$ represents a case where we do not use the nested POD technique.}
\begin{tabular}{|l|c|c|c|c|}
\hline
   &  \multicolumn{1}{l|}{$\mathrm{N_{int}} = 2$} & \multicolumn{1}{l|}{$\mathrm{N_{int}} = 5$} & \multicolumn{1}{l|}{$\mathrm{N_{int}} = 10$} &
   \multicolumn{1}{l|}{$\mathrm{N_{int}} = \infty$} \\ \hline
$\mathrm{N_{nn}} = 4$                           & 6956             & 6932                   &  6923                              & 6983                               \\ \hline
$\mathrm{N_{nn}} = 7$                            & 7114                             & 7108                             & 7064 & 7103                             \\ \hline
$\mathrm{N_{nn}} = 10$                          & 7607                             & 7552                             & 7657                        &     7538      \\ \hline
\end{tabular}
\label{tab:time_nn}
\end{table}

\subsubsection{Sensitivity analysis} \label{sec:sensitivity}

So far, we have seen the ROM framework's performance with only a single instance of $\bm{\mu}$. This section aims to present the model's performance when it is utilized as a sensitivity analysis tool. We have 1000 test cases randomly selected from the parameter range $\mathbb{P} = [0.1, 0.4] \times [0.4, 1.0] \ni (\nu, \alpha) = \bm{\mu}$. We employ two model settings

\begin{figure}[!ht]
   \centering

         \includegraphics[keepaspectratio, height=4.8cm]{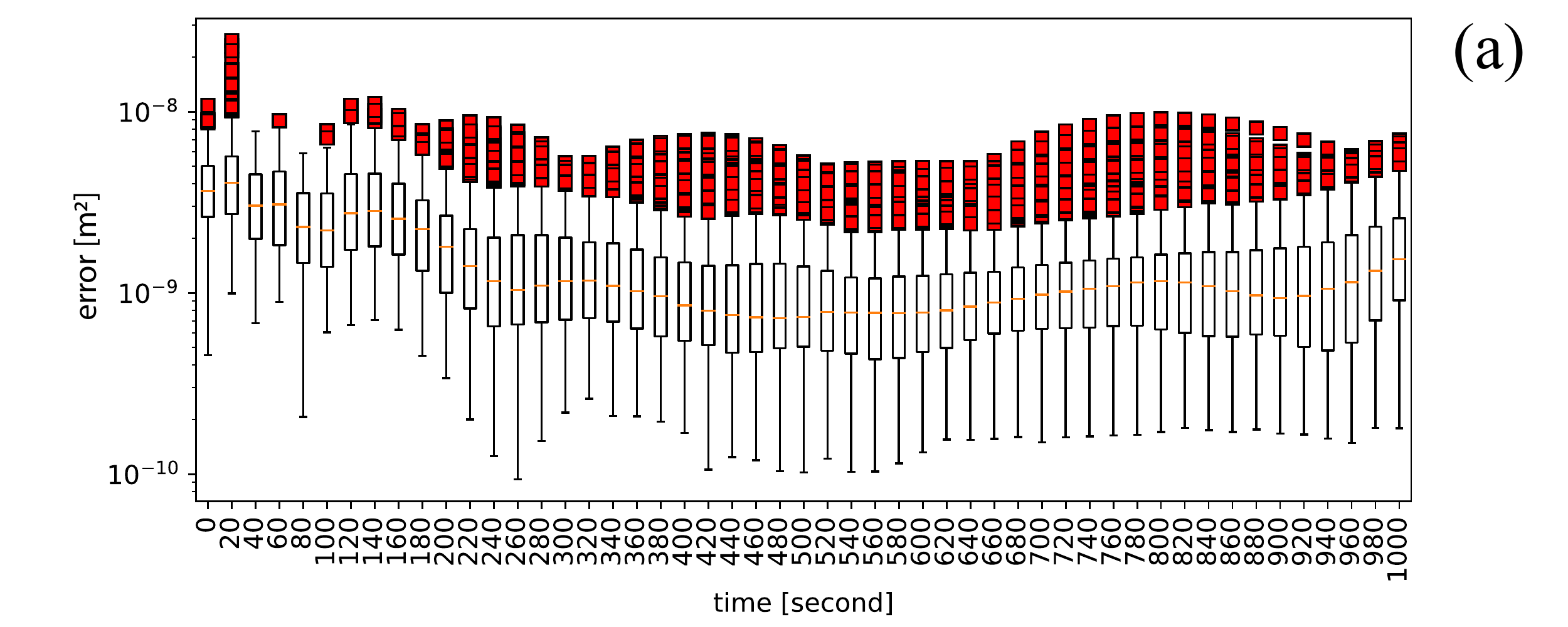}
         \includegraphics[keepaspectratio, height=4.8cm]{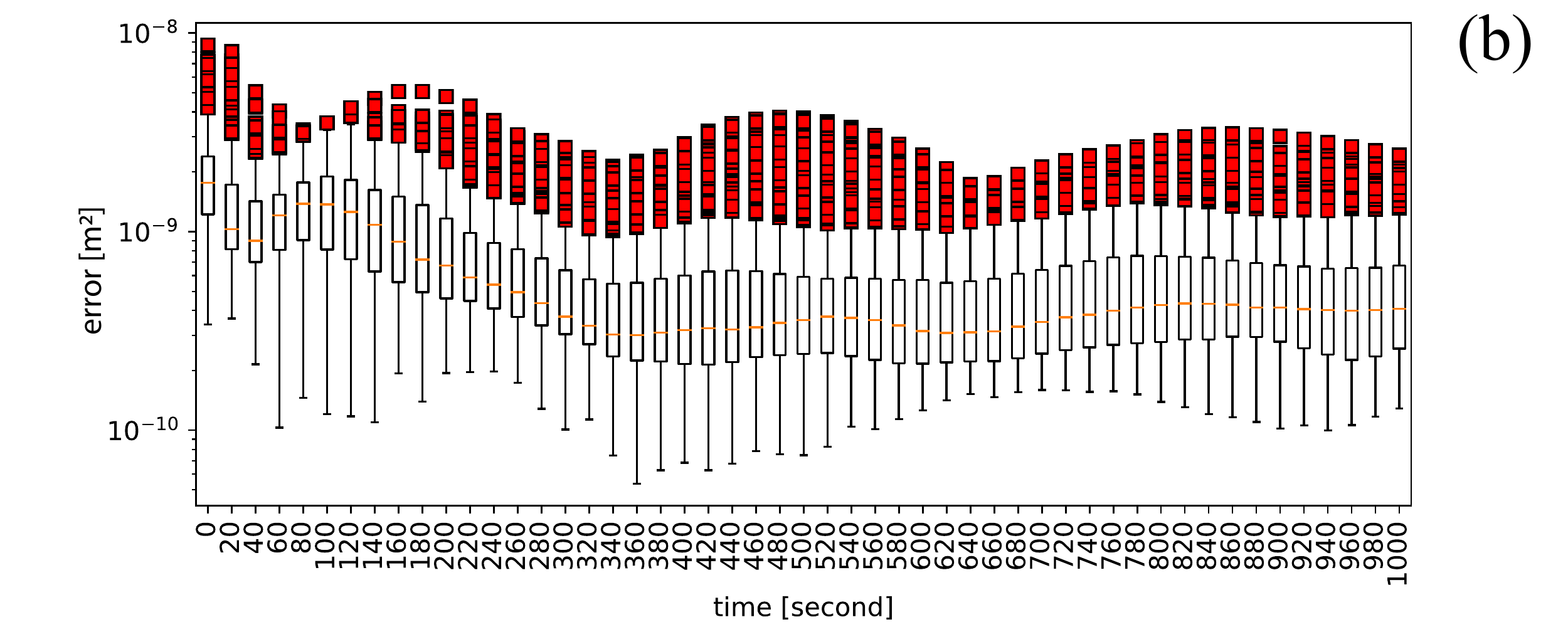}
         \includegraphics[keepaspectratio, height=4.8cm]{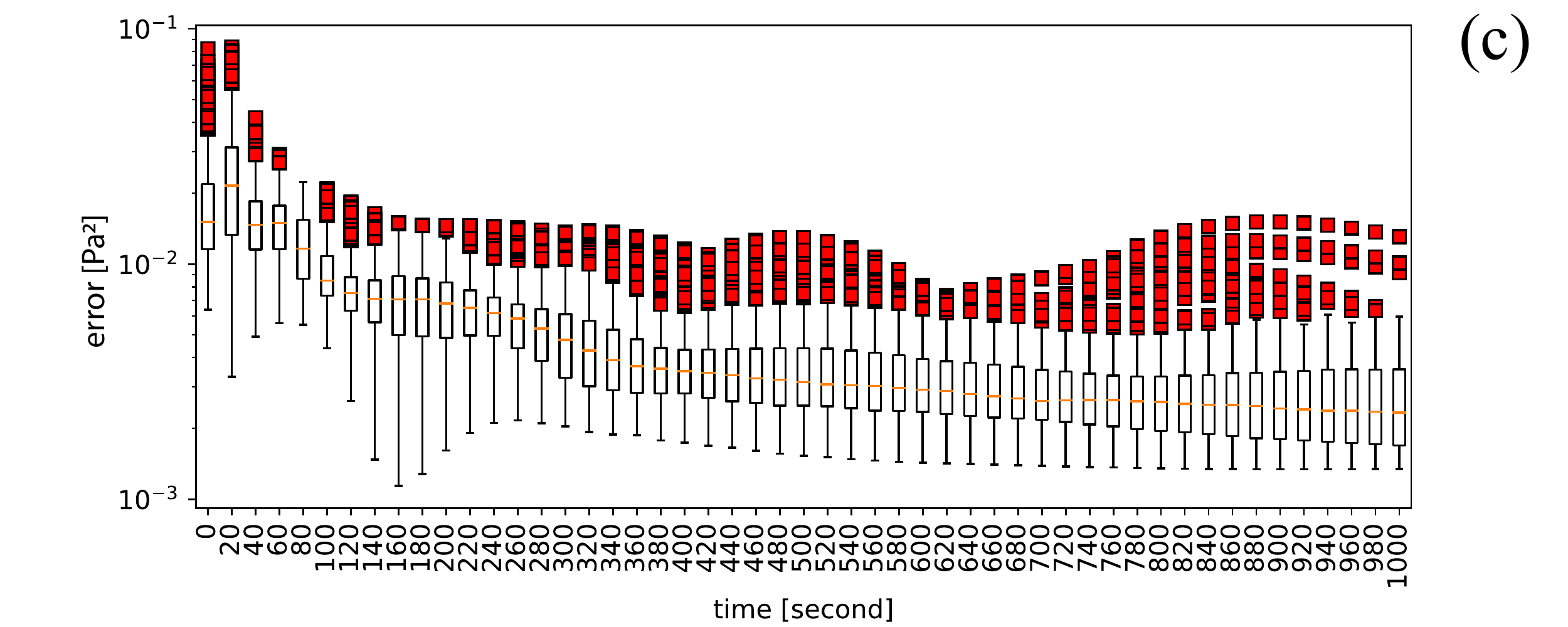}
         \includegraphics[keepaspectratio, height=4.8cm]{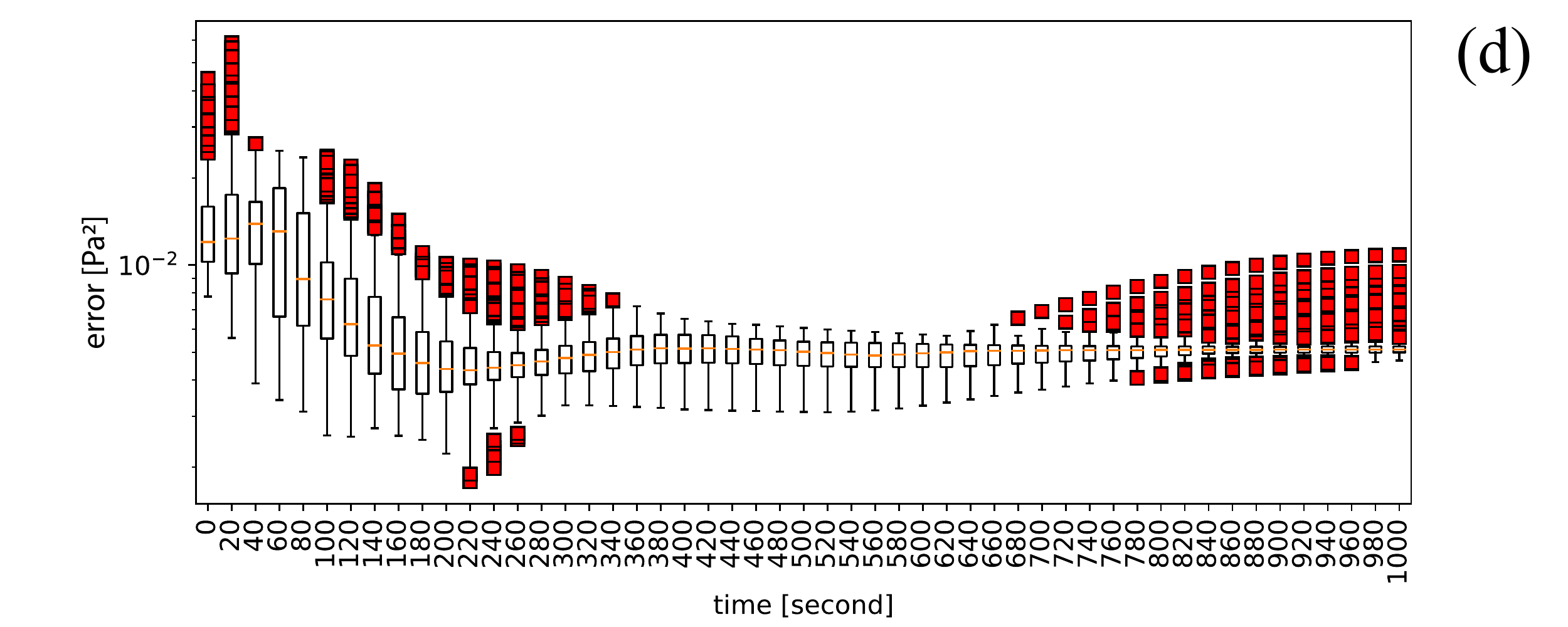}
   \caption{Example 4: Sensitivity analysis - Errors of reconstruction solutions using 1000 testing $\bm{\mu}$: (a) mean squared error (MSE) of displacement field ($\bm{u}$) with $\mathrm{M} = 400$, $\mathrm{N} = 10$, $\mathrm{N_{hl}} = 3$, and $\mathrm{N_{nn}} = 7$, (b) mean squared error (MSE) of displacement field ($\bm{u}$) with $\mathrm{M} = 900$, $\mathrm{N} = 20$, $\mathrm{N_{hl}} = 5$, and $\mathrm{N_{nn}} = 10$, (c) mean squared error (MSE) of fluid pressure field ($p$) with $\mathrm{M} = 400$, $\mathrm{N} = 10$, $\mathrm{N_{hl}} = 3$, and $\mathrm{N_{nn}} = 7$, and (d) mean squared error (MSE) of fluid pressure field ($p$) with $\mathrm{M} = 900$, $\mathrm{N} = 20$, $\mathrm{N_{hl}} = 5$, and $\mathrm{N_{nn}} = 10$. Note that we fix $\mathrm{N_{int}} = 10$. We note that the red squares represent outliers, and the box plot covers the interval from the 25th percentile to 75th percentile, highlighting the mean (50th percentile) with an orange line.}
   \label{fig:ex4_sen}

\end{figure}

\begin{enumerate}
    \item model 1: $\mathrm{M} = 400$, $\mathrm{N_{int}} = 10$, $\mathrm{N} = 10$, $\mathrm{N_{hl}} = 3$, and $\mathrm{N_{nn}} = 7$,
    \item model 2: $\mathrm{M} = 900$, $\mathrm{N_{int}} = 10$, $\mathrm{N} = 20$, $\mathrm{N_{hl}} = 5$, and $\mathrm{N_{nn}} = 10$,
\end{enumerate}

\noindent
and present the range of MSE values of $\bm{u}$ and $p$ fields in Figure \ref{fig:ex4_sen}. We note that the red squares represent outliers, and the box plot covers the interval from the 25th percentile to 75th percentile, highlighting the mean (50th percentile) with an orange line. We note that most of the outliers are located above the 75th percentile, i.e. are all cases in which the error is larger than average. From this figure, it is clear that the second model has approximately one order of MSE magnitude less than the first model. Besides, the range of uncertainties is reduced significantly (see outliers and length of the box plots). Again, the MSE values tend to decrease with time since the solutions approach the steady-state solutions.

Next, we compare the wall time used for ROM and FOM to perform sensitivity analysis (test set of 1000 members) as shown in Table \ref{tab:ex4_sen}. We note that we do not present the wall time used in the initialization of $\bm{\mu}$ (see blue box in Figure \ref{fig:pod_nn_explain}) because it is insignificant compared to the other phases. The second model (i.e., the one with better accuracy) require higher wall time for all operations than the first model. Using ROM, we could speed-up by six-folds, as one could see from the last row in Table \ref{tab:ex4_sen}. A more comprehensive discussion on the effectiveness of the ROM framework will be provided in the following section.


\begin{table}[!ht]
\centering
\caption{Example 4: Comparison of the wall time (seconds) used for sensitivity analysis}
\begin{tabular}{|l|c|c|c|}
\hline
                             & M = 400 (model 1) & M = 900 (model 2) & FOM   \\ \hline
Train FOM snapshots                          & 7160    & 16020   & -     \\ \hline
Perform POD                          & 1437    & 6475    & -     \\ \hline
Train ANN                          & 7064    & 18492   & -     \\ \hline
Prediction - 1000 testing $\bm{\mu}$ & 2895    & 3160    & 17790 \\ \hline
Prediction - per testing $\bm{\mu}$  & 2.9     & 3.2     & 17.8  \\ \hline
\end{tabular}
\label{tab:ex4_sen}
\end{table}

\section{Discussion}

The numerical observations that the benchmark cases in the previous section have highlighted can be summarized by three main points of discussion. First, in Section \ref{sec:source_err} we investigated the sources of the ROM error, and we observed that the main error contribution comes from the prediction of $\widehat{\bm{\theta}}^u$, $\widehat{\bm{\theta}}^p$ by ANN with a given $t$ and $\bm{\mu}$. The main evidence lies in the comparisons of MSE and ME values between Figures \ref{fig:ex4_source_err_1} and \ref{fig:ex4_source_err_2}. We can see that with the same parameter $\bm{\mu}$, the MSE and ME values resulted from $\widehat{\bm{\theta}}^u$, $\widehat{\bm{\theta}}^p$ are about three orders of magnitude higher than those of the ones obtained from ${\bm{\theta}}^u$, ${\bm{\theta}}^p$. Therefore, for future works, we will focus on improving the ANN model's accuracy by using different types of network architecture (recurrent neural networks) or regularization (physics-guided machine learning).

Second, throughout Section \ref{sec:analysis}, we could not observe any clear relationships between the ROM's accuracy and $\mathrm{N_{int}}$ (see Figures \ref{fig:ex4_source_err_2}, \ref{fig:ex4_source_err_3}, \ref{fig:ex4_num_snap}, \ref{fig:ex4_num_nhl}, and \ref{fig:ex4_num_nn} and Tables \ref{tab:l2_source_err_2}, \ref{tab:l2_source_err_3}, \ref{tab:l2_num}, \ref{tab:l2_hl}, and \ref{tab:l2_nn}). As discussed in the previous paragraph and previous sections, the errors introduced by the POD and $L^2$ projection phases are much less than the errors stemming from the ANN phase. Consequently, the errors introduced by the truncation of $\mathrm{N_{int}}$ could not be observed in the final results. This observation implies that we could utilize the nested POD technique to save computational time (see Table \ref{tab:time_num}) without any observable losses in the model's accuracy.

Third, according to Table \ref{tab:ex4_sen}, we could see that the ROM framework is approximately six times faster than the FOM or the finite element model during the online phase. Moreover, the ROM framework errors are very small (see Figure \ref{fig:ex4_sen}), especially relative to the magnitude of the FOM solutions. The ROM framework, however, has an extra cost of training (i.e., initialization, FOM, POD, $L^2$ projection, and ANN phases, see Figure \ref{fig:pod_nn_explain}). The training time (wall time) of model 1 and model 2 in Section \ref{sec:sensitivity} are 15661 and 40987 seconds, respectively. Taking the training time into account, we need to perform at least 1050 and 2850 inquiries (online phase) to have a break-even point for model 1 and model 2, respectively. To this end, before one wants to build a ROM framework, one should consider how many inquiries are expected to have.

\section{Conclusion}

A non-intrusive reduced order model (ROM) has been developed for linear poroelasticity problems in heterogeneous media. We employ the discontinuous Galerkin (DG) finite element framework as a full order model (FOM); the DG solutions are thus employed as snapshots to train and test the ROM. During the offline phase, this framework utilizes one of two variants of the proper orthogonal decomposition (POD) to define a reduced basis space, namely a standard POD and a nested POD, and artificial neural networks (ANN) to construct an inexpensive map from a time and parameter pair to coefficients associated with each reduced basis. We validate the framework through a series of benchmark problems.  Our results show that the framework could provide reasonable approximations of the FOM results, but it is significantly faster. Moreover, the reduced order framework can capture both displacement and pressure fields' sharp discontinuities resulting from the heterogeneity in the media's conductivity. We then present the error sources and show that the error inherited from the ANN model trumps the error associated with the POD operation. Consequently, we illustrate that the nested POD technique, in which time and uncertain parameter domains are compressed consecutively, could provide comparable accuracy to the classical POD method, in which all domains are compressed simultaneously, but at a fraction of the offline computational cost. Finally, we emphasize in which circumstances the ROM framework is more suitable than the FOM. Further developments could be to consider different ANN architectures, as well as coupled problems involving poroelasticity.

\section{Acknowledgements}
 The computational results in this work have been produced by the RBniCS project \cite{ballarin2015rbnics} (a reduced order modeling library built upon FEniCS \cite{AlnaesBlechta2015a}), the multiphenics library~\cite{Ballarin2019} (an extension of FEniCS for multiphysics problems), and PyTorch \cite{NEURIPS2019_9015}. We acknowledge the developers of and contributors to these libraries.
FB thanks Horizon 2020 Program for Grant H2020 ERC CoG 2015 AROMA-CFD project 681447 that supported the development of RBniCS and multiphenics. NB acknowledges startup  support from the Sibley School of Mechanical and Aerospace Engineering, Cornell University.

\bibliographystyle{unsrt}
\bibliography{references.bib}

\end{document}